\documentclass[10pt]{memoir} 
\usepackage{amsmath,amsthm,amsfonts,amscd, amssymb, epsf}
\usepackage{palatino}
\usepackage{graphicx}
\usepackage{geometry}
\title{The Selberg Trace Formula and Selberg Zeta-Function for Cofinite Kleinian Groups with Finite Dimensional Unitary Representations \\ {\small Stony Brook University PhD Thesis}}
\author{Joshua S. Friedman}
\date{May 2005} 

\numberwithin{section}{chapter}
 \theoremstyle{plain}    
 \newtheorem{thm}{Theorem}[section]
 
 \numberwithin{equation}{section} 
 \numberwithin{figure}{section} 
 \theoremstyle{plain}
 \theoremstyle{plain}    
 \newtheorem{lem}[thm]{Lemma} 
 \newtheorem{rem}[thm]{Remark}
 \newtheorem{as}[thm]{Assumption} 
 \newtheorem{cor}[thm]{Corollary}
 \newtheorem*{cor*}{Corollary}
 \newtheorem*{conj*}{Conjecture}
 \theoremstyle{definition}
 \newtheorem{defn}[thm]{Definition}
 \newtheorem{nota}[thm]{Notation}

 \theoremstyle{plain}    
 \newtheorem{prop}[thm]{Proposition}
  
\newcommand{\bl}{\begin{lem}}
\newcommand{\el}{\end{lem}}

\DeclareMathOperator{\sz}{\mathnormal{Z}(\mathnormal{s},\Gamma,\chi)}

\newcommand{\bnot}{\begin{nota}}
\newcommand{\enot}{\end{nota}}

\newcommand{\ben}{\begin{enumerate}}
\newcommand{\een}{\end{enumerate}}

\newcommand{\bml}{\begin{multline}}
\newcommand{\eml}{\end{multline}}

\newcommand{\beq}{\begin{equation}}
\newcommand{\eeq}{\end{equation}}
\newcommand{\bp}{\begin{prop}}
\newcommand{\ep}{\end{prop}}
\newcommand{\bd}{\begin{defn}}
\newcommand{\ed}{\end{defn}}
\newcommand{\pf}{\begin{proof}}
\newcommand{\epf}{\end{proof}}

\newcommand{\field}[1]{\ensuremath{\mathbb{#1}}}
\newcommand{\CC}{\field{C}}

\newcommand{\NN}{\field{N}}
\newcommand{\QQ}{\field{Q}}

\newcommand{\HHHH}{\field{H}}
\newcommand{\hh}{\field{H} \,}

\DeclareMathOperator{\HC}{\HHHH^3}

\DeclareMathOperator*{\lto}{\mathnormal{o}(1)}

\DeclareMathOperator{\df}{\equiv}
\DeclareMathOperator{\lsum}{\mathnormal{\sum}}
\DeclareMathOperator*{\psum}{ \lsum^\prime}
\DeclareMathOperator{\gip}{\mathnormal{\Gamma ^ \prime _\infty   }}
\DeclareMathOperator{\gi}{\mathnormal{\Gamma_\infty   }}
\DeclareMathOperator{\HH}{\HHHH^3}

\newcommand{\PP}{\field{P}}

\newcommand{\BB}{\field{B}}
\newcommand{\RR}{\field{R}}

\newcommand{\ZZ}{\field{Z}}

\newcommand{\bh}{\textbf{H}}

\newcommand{\py}{\varphi}
\newcommand{\F}{\mathcal{F}}
\newcommand{\E}{\mathcal{E}}
\newcommand{\scz}{\mathcal{S}}
\newcommand{\pg}{\mathcal{P}}
\newcommand{\K}{\mathcal{K}}
\newcommand{\D}{\mathcal{D}}
\newcommand{\hil}{\mathcal{H}}

\newcommand{\mC}{\mathcal{C}}

\newcommand{\smat}{\mathfrak{S}}

\DeclareMathOperator{\rep}{Rep(\Gamma,\mathnormal{V})}
\DeclareMathOperator{\rp}{Rep(\Gamma,\mathnormal{V})}
\DeclareMathOperator{\PSL}{PSL}

\DeclareMathOperator{\SL}{SL}

\DeclareMathOperator{\vol}{vol}

\DeclareMathOperator{\R}{Re}
\DeclareMathOperator{\I}{Im}
\DeclareMathOperator{\pc}{PSL(2,\CC)}
\DeclareMathOperator{\GL}{GL}

\DeclareMathOperator{\ID}{id}
\DeclareMathOperator{\CUSP}{cusp}
\DeclareMathOperator{\PAR}{par}
\DeclareMathOperator{\CE}{ce}
\DeclareMathOperator{\NCE}{nce}
\DeclareMathOperator{\LOX}{lox}
\DeclareMathOperator{\cuspi}{\mathcal{CE}}

\DeclareMathOperator{\lds}{\mathnormal{ \frac{\phi^{\prime}}{\phi}}}

\DeclareMathOperator{\tr}{tr}

\DeclareMathOperator{\en}{\mathnormal{\mathcal{E}(T)}}
\DeclareMathOperator{\oen}{\mathnormal{\left|\mathcal{E}(T) \right|}}

\DeclareMathOperator{\ren}{\mathnormal{\mathcal{E}(R)}}

\DeclareMathOperator{\cinf}{\mathnormal{\PP}}

\DeclareMathOperator{\hs}{\mathnormal{\hil(\Gamma,\chi)}}
\DeclareMathOperator{\lp}{\mathnormal{\Delta}}

\DeclareMathOperator{\rv}{\mathnormal{F_{\chi}(P,Q,s)}}

\newcommand{\Z}{\mathbb{Z}}

\renewcommand{\Re}{\operatorname{Re}}

\newcommand{\<}{\langle}
\renewcommand{\>}{\rangle}

\newcommand{\ra}{\rightarrow}

\newcommand{\union}{\cup}

\newcommand{\tdv}{two-dimensional vector }
\newcommand{\td}{two-dimensional }
\newcommand{\thd}{three-dimensional  }
\newcommand{\fd}{finite-dimensional  }
\newcommand{\thds}{three-dimensional scalar }
\newcommand{\thdv}{three-dimensional vector }
\newcommand{\ckus}{cofinite Kleinian groups with finite-dimensional unitary representations}
 \usepackage{setspace} 
\begin{document}
\maketitle

\begin{abstract}
For cofinite Kleinian groups,  with finite-dimensional unitary representations, we derive the Selberg trace formula.  As an application we define the corresponding Selberg zeta-function and compute its divisor, thus generalizing results of Elstrodt, Grunewald and Mennicke  to non-trivial unitary representations.  We show that the presence of cuspidal elliptic elements sometimes adds ramification point to the zeta function.  In fact, if \mbox{$\mathcal{O} = \ZZ[-\frac{1}{2}+\frac{\sqrt{-3}}{2}]$} is the ring of Eisenstein integers, then the Selberg zeta-function of $\PSL(2,\mathcal{O})$ contains ramification points.
\end{abstract}
\newpage
\tableofcontents* 

\section*{Acknowledgment}
I would like to thank my thesis advisor Professor Leon Takhtajan for the years he has spent guiding and teaching me.  I would also like to thank Professor J\"{u}rgen Elstrodt for reading over the results of this thesis, and for many useful suggestions.  Special thanks  are due to Jay Jorgenson, Irwin Kra, Lee-Peng Teo, Alexei Venkov, and Peter Zograf for useful comments and suggestions.

I would also like to thank the Stony Brook Mathematics Department for supporting me while I carried out this research.

\chapter{Introduction}
The Selberg theory (Selberg trace formulas, Selberg zeta-functions, and related applications) has been well studied in both the two-dimensional scalar case (\cite{Iwaniec}) and the two-dimensional vector case\footnote{The works \cite{Roelcke}, \cite{Hejhal}, \cite{Fischer}, contain not only the two-dimensional vector case, but also its generalization, the case of unitary multiplier 
systems of arbitrary real weight. \\ } (\cite{Roelcke}, \cite{Venkov}, \cite{Hejhal}, \cite{Fischer}).  By ``Two-dimensional vector case'' we mean: cofinite Fuchsian groups with finite-dimensional unitary representations, and the ``Scalar case'' refers to the case with the trivial  representation.  Elstrodt, Grunewald and Mennicke extended the Selberg theory to the three-dimensional scalar case in \cite{Elstrodt}.  By the ``Three-dimensional case'' we mean: cofinite Kleinian groups.  The main goal of this thesis is to extend the Selberg theory to the three-dimensional vector case.

In this thesis we derive the Selberg trace formula for cofinite Kleinian groups\footnote{A  Kleinian groups is referred to in some texts as a discrete group of isometries acting on hyperbolic three-space, or a discrete subgroup of $\pc.$ \\ },  with finite-dimensional unitary representations.  As an application we define the corresponding Selberg zeta-function and compute its divisor, thus generalizing results of Elstrodt, Grunewald and Mennicke \cite{Elstrodt} to non-trivial unitary representations. 

Much of the \tdv and  three-dimensional scalar cases extends in a straight forward manner to the \thdv case.  However, the extension of several parts of the Selberg theory are more subtle in the three-dimensional case, especially in the vector case.  One reason for this is because the set of finite-dimensional unitary representations of a fixed cofinite Kleinian groups is not well understood.  Another reason is related to the structure of the stabilizer subgroup of a cusp.  In the \td case the stabilizer subgroup of a cusp is a purely parabolic group that is isomorphic to a rank-one lattice,  while in three dimensions the stabilizer subgroup of a cusp is a non-abelian group that contains elliptic elements, with a finite-index purely parabolic subgroup that is  isomorphic to a rank-two lattice. The presence of elliptic elements in the stabilizer subgroup introduces some subtleties to the \thdv case, particularly in the computation of the divisor of the Selberg zeta-function.  In addition, the fact that the stabilizer subgroup (in the \thd case) contains a rank-two parabolic subgroup forced us to prove some additional estimates involving two-dimensional lattice sums.


A  Kleinian group  is a discrete subgroup of $\PSL(2,\CC) =   \SL(2,\CC) / \{\pm I\}. $ Each element of $\PSL(2,\CC) $ is identified with a M\"{o}bius transformation, and has a well-known action on hyperbolic three-space  $\HH$  and on its boundary at infinity$-$ the Riemann sphere $\PP^1$ (see \cite[Section 1.1]{Elstrodt})  . A Kleinian  group is \emph{cofinite} iff it has a fundamental domain $\F \subset \HH $ of finite hyperbolic volume. 

We use the following coordinate system for hyperbolic three-space, 
$$\HH \df   \{(x,y,r)\in\RR^{3}~|~r>0 \} \df  \{ (z,r) ~| z \in \CC, ~r > 0 \} \df    \{z + rj\in\RR^{3}~|~r>0 \}, $$ with the hyperbolic metric 
$$ ds^{2} \df \frac{dx^{2}+dy^{2}+dr^{2}}{r^{2}}, $$ and volume form 
$$ dv \df \frac{dx\, dy\, dz}{r^{3}}. $$
The Laplace-Beltrami operator is defined by
$$ \lp \df -r^{2}(\frac{\partial^{2}}{\partial x^{2}}+
\frac{\partial^{2}}{\partial y^{2}}+ \frac{\partial^{2}}{\partial r^{2}})+ r\frac{\partial}{\partial r}, $$ and it acts on the space of  smooth functions $f:\HH \mapsto V, $ where $V$ is a \fd complex vector space with inner-product $\langle~,~\rangle_V.$  

Suppose that $\Gamma $ is a cofinite Kleinian group and  $\chi \in \rep$ ($\rep$ is the space of \fd unitary representations of $\Gamma$ in $V$). Then the Hilbert space of \emph{$\chi-$automorphic} functions is defined by  
\begin{multline*}
\hs  \df  \{ f: \HH \ra V ~|~ f(\gamma P) = \chi(\gamma) f(P)~\forall \gamma \in \Gamma, \\ P \in \HH, $ and $\left<f,f \right> \df \int_{\F} \left<f(P),f(P)\right>_V\,dv(P) < \infty \}. 
\end{multline*}
Here $\F$ is a fundamental domain for $\Gamma $ in $\HH$,  and $\left<~,~\right>_V $ is the inner product on $V.$ Finally, let   $ \lp = \lp(\Gamma,\chi) $ be the corresponding positive self-adjoint Laplace-Beltrami operator on $\hs.$

Our first result is the spectral decomposition of $\lp$  on $\hs$ (see Theorem \ref{T:spectral}).  Except for one important point, the proof of the spectral decomposition theorem is analogous to the \tdv and \thds cases.  The one important point being, \emph{singularity at a cusp}.  To the best of the author's knowledge, prior to this thesis, the notion of singularity at a cusp was only defined for cofinite Fuchsian groups \cite{Selberg1} \cite{Roelcke} \cite{Venkov} \cite{Hejhal}.  In \S\ref{secSingular} we extend the notion of singularity to cofinite Kleinian groups.  

In  \S\ref{sectionSelberg} we give an explicit form of the Selberg trace formula for cofinite Kleinian groups with finite-dimensional unitary representations (see Theorem \ref{T:Selberg}).
The new feature in the trace formula is a term of the form,
$$ \sum _{\alpha=1}^{\kappa}\frac{g(0)}{|\Gamma _{\alpha}:\Gamma '_{\alpha}|}\sum _{k=l_{\alpha}+1}^{\dim_\CC V}L(\Lambda_{\alpha},\psi_{k \alpha}). $$
The above term comes about from  \emph{regularity} at a cusp, and its value is computed using Kronecker's second limit formula ( see  \S\ref{secDoubleSum}).
   
As an application of the spectral decomposition theorem we derive an identity  involving conjugacy relations of cuspidal elliptic elements.  This identity is used in the proof of the Selberg trace formula and to show that under certain conditions, the Selberg zeta function admits a meromorphic continuation (see Lemma~\ref{lemCuspElip} for the identity).

For $\R(s)>1 $ the Selberg zeta-function $Z(s,\Gamma,\chi)$ is defined by the following product 
$$
Z(s,\Gamma,\chi) \df \prod_{ \{T_0 \} \in \mathcal{R}} ~ \prod_{j=1}^{ \dim_\CC V} \prod_{  \substack{ l,k \geq 0 \\  c(T,j,l,k)=1   } } \left( 1-\mathfrak{t}_{j} a(T_0)^{-2k} \overline{ a(T_0) ^{-2l}} N(T_0)^{-s - 1}    \right).  
$$
In \S\ref{secSZF} we introduce the various definitions and notations that are needed in order to define the zeta function, and   meromorphically\footnote{We give the meromorphic continuation for certain cases, and for others show that the zeta function is a rational root of a meromorphic function. \\}  continue    $Z(s,\Gamma,\chi)$ to $\R(s) \leq 1$ while computing its divisor. The main new difficulty is handling the contribution of the cuspidal elliptic elements to the topological (or trivial) zeros and poles of $Z(s,\Gamma,\chi).$  We show the following in \S\ref{secTopZeros}:
\begin{cor*}
Let $\Gamma = \PSL(2,\ZZ[-\frac{1}{2}+\frac{\sqrt{-3}}{2}]),$ and $\chi \equiv 1$ (the trivial representation). Then $Z(s,\Gamma,\chi)$ is not a meromorphic function\footnote{This is the first example that the author is aware of where the Selberg zeta-function is not meromorphic.\\ } (it is the 6-th root of a meromorphic function).
\end{cor*}
In addition, the methods of \S\ref{secTopZeros} imply that the Selberg zeta-function of the Picard group is meromorphic:
\begin{cor*}
Let $\Gamma = \PSL(2,\ZZ[\sqrt{-1}]),$ and let $\chi \in \rep.$ Then $Z(s,\Gamma,\chi)$ is a meromorphic function.
\end{cor*}

\chapter{Prerequisite Material}
In order to make our presentation self-contained, we present some well-known results concerning hyperbolic three-space, the Laplace-Beltrami operator, and cofinite Kleinian groups.  For more details see  \cite{Elstrodt}.  

\section{Hyperbolic three-space $\HH$}
Let $\HH$ denote the upper half space (of $\RR^{3}$) model of  hyperbolic three-space.  The space $\HH$ is parametrized by the following coordinates:  
$$ \HH  \df  \{(x,y,r)\in\RR^{3}~|~r>0 \}  \df  \{(z,r)~| ~z \in\CC,  r>0 \}. $$ A point $P \in \HH$ will be denoted by $(x,y,r),~(z,r),$ or $z+rj.$ 
The standard hyperbolic metric\footnote{The metric with -1 sectional curvature.} hyperbolic metric and volume form are written respectively as
$$ 
ds^{2} \df \frac{dx^{2}+dy^{2}+dr^{2}}{r^{2}}~\text{and}~dv \df \frac{dx\, dy\, dz}{r^{3}}. $$ 

\begin{rem} \label{remRoundHyp} An alternate model of hyperbolic three-space is the open three-ball, $~\BB^3 \df \{ (x,y,z)\in \RR^3~|~x^2 + y^2 + z^2 < 1~  \}. $ When equipped with the metric $$ 4 \frac{dx^2 + dy^2 + dz^2}{ \left( 1-x^2 -y^2 - z^2 \right)^2},$$ $\BB^3$ is isometric to $\HH.$  
\end{rem}

The boundary (at infinity) of $\HH $ can be realized (by Remark \ref{remRoundHyp}) as the Riemann sphere $\PP = \CC \cup \infty. $

The Laplace-Beltrami $\lp$ operator lies at the heart of this thesis.  In our coordinates it can be written explicitly as the following differential operator\footnote{In our notation $\lp$ is a positive self-adjoint operator \\}:
$$ \lp \df -r^{2}(\frac{\partial^{2}}{\partial x^{2}}+
\frac{\partial^{2}}{\partial y^{2}}+
\frac{\partial^{2}}{\partial r^{2}})+
r\frac{\partial}{\partial r}. $$ 

The orientation preserving isometry group of  $\HH $ can be identified with the group $\PSL(2,\CC)  = \SL(2,\CC) / \{ \pm I \} $. Each element 
$$M=
\left(\begin{array}{cc}
a & b\\
c & d\end{array}\right)
\in\PSL(2,\CC)$$ acts on $\HH$ as follows:  $M(z+rj)=w+tj, $ where 
$$
w=\frac{(az+b)(\bar{c}\bar{z}+\bar{d})+ a\bar{c}r^{2}}{|cz+d|^{2}+|c|^{2}r^{2}} ~~\text{and}~~ t=\frac{r}{|cz+d|^{2}+|c|^{2}r^{2}}.  $$
The element $M$ also acts on the  boundary at infinity of $\HH$ via  standard M\"{o}bius action on $\PP,$  
$$ M\zeta = \frac{a\zeta + b}{c\zeta + d} $$ for $\zeta \in \PP. $ 

\section{Harmonic Analysis on $\HH$}  \label{secFreCas}
For $P=z+rj,P'=z'+r'j \in \HH$ denote by  $d(P,P')$  the (hyperbolic) distance in $\HH$ between $P$ and $P',$ and let $\delta(P,P')$ be defined by 
$$\delta(P,P') \df  \frac{|z-z'|^{2}+r^{2}+ r'^{2}}{2rr'}. $$
It follows that $\cosh(d(P,P')) = \delta(P,P'),$ and that $\delta$ is a point-pair invariant\footnote{A point-pair invariant is a map  $f:\HC\times\HC\rightarrow\CC$ defined almost everywhere  satisfying $
f(MP,MQ)=f(P,Q)  $ for all $P,Q\in\HC$, $M\in\PSL(2,\CC).$ \\}.   

We can use the concept of a point-pair invariant to construct the resolvent kernel for $\lp.$ For  $s \in \CC, \, t>1,$ set  
$$
\py_{s}(t) \df \frac{\left(t+\sqrt{t^{2}-1}\right)^{-s}}{\sqrt{t^{2}-1}}. $$ 
\bl \label{lemScalRes} \cite[Lemma 4.2.2]{Elstrodt}
Let  $u \in C_c^2(\HH), $
$s\in\CC$ and  $\lambda=1-s^{2}.$ Then for all $Q\in\HC$
$$
u(Q)=\frac{1}{4\pi}\int_{\HH}\varphi_{s}(\delta(P,Q))(\Delta-\lambda)u(P)\, dv(P). $$

(2) The point-pair invariant $\py_{s}\circ \delta$ is the resolvent kernel for $\lp$ on the Hilbert space of square integrable functions on $\HH.$
\el

Let $ \scz \df \scz[1,\infty)  $ denote the Schwartz space of smooth functions \linebreak \mbox{$k:[1,\infty) \rightarrow \CC$} that satisfy $ \lim_{x \ra \infty} x^n k^{(m)}(x) = 0 $ for all $n,m \in \NN_{\geq 0} $.  For each $k \in \scz ,$  $K \df k \circ \delta $ is a point-pair invariant and the kernel of an operator $ \mathcal{K}:L^2_{\text{loc}}(\HH) \mapsto L^2_{\text{loc}}(\HH)  $ defined by
$$
\mathcal{K}f(P) = \int_{\HH} K(P,Q) f(Q) ~dv(Q).
$$
We have
\bl \cite[Lemma 3.5.3]{Elstrodt} \label{E:Selberg Transform}
Let $k \in \scz $, $K(P,Q)=k(\delta (P,Q))$ for $P,Q \in \HH$,   $f:\HH \ra \CC $ be a solution of $\lp f = \lambda f$,  $ \lambda = 1-s^2 $, and 
\beq
h(\lambda)\df h(1-s^2)\df \frac{\pi}{s}
\int_{1}^{\infty}k\left(\frac{1}{2}\left(t+\frac{1}{t}\right)\right)
(t^{s}-t^{-s})\left(t-\frac{1}{t}\right)\,\frac{dt}{t}.
\eeq
Then \eqref{E:Selberg Transform} converges absolutely and 
$$ 
\int_{\HH}K(P,Q)f(Q)\, dv(Q)=h(\lambda)f(P). 
 $$
\el

The lemma above says that if $f$ is an eigenfunction\footnote{The function $f$ need not be in any Hilbert space. \\} of $\lp $ with eigenvalue $\lambda $  then $f $ is also an eigenfunction of $\mathcal{K} $ with eigenvalue $h(\lambda)$ depending only on the eigenvalue $\lambda $ and not on the particular eigenfunction.

The map sending $ k \in \scz $ to holomorphic function $h$ above is called the \mbox{Selberg-Harish-Chandra} transform of $k.$

\section{Kleinian Groups} \label{secKG}
A subgroup $\Gamma<\PSL(2,\CC)$ is called a Kleinian group if for each $P\in\HH$ the orbit $\Gamma P$ has no accumulation points in $\HH$. An equivalent formulation is that  $\Gamma$ is a discrete subset of $\PSL(2,\CC)$ in the topology induced from $\CC^4.$

\noindent A closed subset $\F \subset \HH$ is called a fundamental domain of $\Gamma$ if 
\begin{itemize}
\item $\F$ meets each $\Gamma-$orbit at least once, 
\item the interior $\F^{o}$ meets each $\Gamma-$orbit at most once, 
\item the boundary of $\F$ has Lebesgue measure zero.
\end{itemize}
\noindent For each $Q \in \HH$  the set $$
\mathcal{P}_{Q}(\Gamma) \df \{ P\in\HH\,|\, d(P,Q)\leq d(\gamma P,Q)\,\,\,\,\forall\,\gamma\in\Gamma\} $$ is a fundamental domain for $\Gamma$ that is centered at the point $Q.$ 

We say that $\Gamma$ is \emph{cocompact} if it has a fundamental domain $ \F $ that is compact, and  \emph{cofinite}\footnote{Note that cocompact groups are also cofinite. \\} if it has a fundamental domain $ \F $ with
$$ \vol(\Gamma) \df \int_{\F}\, dv < \infty $$ ($dv$ is the volume  form of $ \HH $).   
\section{M\"{o}bius Transformations}
Each element $$\gamma=\left(\begin{array}{cc}
a & b\\
c & d\end{array}\right)\in\PSL(2,\CC) $$ falls into exactly\footnote{The identity element is an exception. Its trace is $\pm 2$ but it is not usually thought of as a parabolic element. \\} one of the following categories:
\begin{itemize} 
\item  \emph{parabolic} if $|\tr(\gamma)|=2$ and $\tr\gamma\in\RR,$
\item \emph{elliptic} if $0\leq|\tr(\gamma)|<2$ and $\tr\gamma\in\RR,$
\item \emph{loxodromic} if it is neither \emph{elliptic} nor \emph{parabolic}.
\end{itemize}

There is a useful geometric characterization of the above notions.
An element $\gamma \in \pc $ is 
\begin{itemize} 
\item  parabolic\footnote{See the previous footnote. \\} iff it has exactly one fixed point in $\cinf$,
\item  elliptic iff it has two fixed points in $\cinf$ and fixes the geodesic line in $\HH$ connecting the two points.
\item loxodromic iff it has two fixed points in $\cinf$ and has no fixed points in $\HH.$ \end{itemize}

\section{Stabilizer Subgroups and Cusps} \label{secStaCus}
Let $\Gamma $ be a Kleinian group. For each $Q \in\HH\cup\cinf$ the stabilizer subgroup of $Q$ is denoted by 
$ \Gamma_{Q} \df \{\,\gamma\in \Gamma \,|\,\gamma Q=Q\,\} $
and 
for $ \zeta \in \cinf, $  $$
\Gamma_{\zeta}^{\prime} \df \{ \gamma \in \Gamma_{\zeta} \, | \, \gamma \, \text{ is parabolic or the identity} \,  \}.  
$$
Set 
$$  B(\CC) \df \left\{ \,\left.\left(\begin{array}{cc}
a & b\\
0 & a^{-1}\end{array}\right)\,\right|\,0\neq a\in\CC,\, b\in\CC\,\right\} /\{\pm I\}<\pc, $$
 $$ N(\CC) \df \left\{ \,\left.\left(\begin{array}{cc}
1 & b\\
0 & 1\end{array}\right)\,\right|\,\, b\in\CC\,\right\} /\{\pm I\}<\pc. $$

The group $B(\CC) $ is the stabilizer subgroup of $\pc$ that fixes $\infty$ and $N(\CC)$ is the  maximal parabolic subgroup of $B(\CC). $ It contains all parabolic elements fixing $\infty,$ and the identity element. 

A point $\zeta \in \cinf $  is called a \emph{cusp} of $\Gamma$ if $ \Gamma_{\zeta}^{\prime} $ is a free abelian group of rank two.  The set of cusps  is denoted by $C_{\Gamma}$. Two cusps, $\alpha ,\ \text{and} \, \beta $ of $\Gamma$ called equivalent or $\Gamma-$equivalent if\footnote{$\Gamma \beta$ denotes the orbit of the point $\beta \in \PP.$ \\} $\alpha \in \Gamma \beta.$  The equivalence class of cusps is denoted by $\Gamma\setminus C_{\Gamma}.$

The following is well known (see \cite[Chapter 2]{Elstrodt}).
\begin{lem}\label{lemFiniteGeo}
\noindent Let $\Gamma$ be a Kleinian group. 
\begin{enumerate}
\item If $\Gamma$ contains a parabolic element then $\Gamma$ is not
cocompact.
\item If $\Gamma$ is cofinite and is not cocompact then $\Gamma$ contains a parabolic element.
\item If $\Gamma$ is cofinite, $\zeta\in\cinf$, and $\Gamma_{\zeta}$
contains a parabolic element then $\zeta$ is a cusp of $\Gamma.$
\item If $\Gamma$ is cofinite then $\Gamma$ has only finitely many
$\Gamma-$equivalent classes of cusps.
\end{enumerate} 
\end{lem}
\noindent Next we study the structure\footnote{The following structure theorem is extremely important, and without it we would be unable to extend the Selberg theory to cofinite Kleinian groups. \\} of the stabilizer subgroup of a cusp for cofinite Kleinian groups.
\begin{lem} \cite[Theorem 2.1.8 part (3)]{Elstrodt}
\label{lemCuspStructI}Let $\Gamma$ be cofinite with a cusp 
at $\infty.$ Then $\gip$ is a lattice in $B(\CC) \approx \CC$ and one of the following three holds.
\begin{enumerate}
\item  $\Gamma_{\infty}=\Gamma'_{\infty}$
\item $\Gamma_{\infty}$ is conjugate in $B(\CC)$ to a group of the
form \[
\left\{ \,\left.\left(\begin{array}{cc}
\epsilon & \epsilon b\\
0 & \epsilon^{-1}\end{array}\right)\,\right|\,0\neq b\in\Lambda,\,\epsilon\in\{1,i\}\,\right\} /\{\pm I\}\]
 where $\Lambda<\CC$ is an arbitrary lattice. As an abstract group
$\Gamma_{\infty}$ is isomorphic to $\ZZ^2 \rtimes \ZZ/2 \ZZ$ where the
nontrivial element of $\ZZ/2 \ZZ$ acts by multiplication by $-1.$
\item $\gip$ is conjugate in $B(\CC)$ to a group of the form \[
\Gamma(n,t)=\left\{ \,\left.\left(\begin{array}{cc}
\epsilon & \epsilon b\\
0 & \epsilon^{-1}\end{array}\right)\,\right|\,\begin{array}{l}
b\in\mathcal{O}_{n}\\
\epsilon=\exp\left(\frac{\pi ivt}{n}\right)\, for\,\,1\leq v\leq2n\end{array}\,\right\} /\{\pm I\}\]
 where $n=4$ or $n=6$ and $t|n$ and where $\mathcal{O}_{n}$ is
the ring of integers in the quadratic number field $\QQ\left(\exp\left(\frac{2\pi i}{n}\right)\right).$
Hence , as an abstract group $\gi$ is isomorphic to the group $\ZZ^2\rtimes\ZZ/m\ZZ$
for some $m\in\{1,2,3,4,6\}.$ An element \[
\epsilon\in\ZZ/m\ZZ\cong\left\{ \exp\left(\frac{\pi iv}{n}\right)\,|v\in\ZZ\,\right\} \]
 acts on $\ZZ^{2}\cong\mathcal{O}_{m'}^{+}$ by multiplication with
$\epsilon^{2}$, where $m'=4$ in case \linebreak $m\in\{1,2,4\}$ and $m'=6$
otherwise. 
\end{enumerate}
\end{lem}
\section{Cofinite Kleinian Groups} \label{secCofKleGro}
Let $\Gamma $ be a cofinite  Kleinian group. By Lemma~\ref{lemFiniteGeo}   $\Gamma $ has finitely many equivalence classes of cusps.  
\bnot
Unless otherwise noted $\Gamma$ is a cofinite Kleinian group with $$ \kappa \df \left| \Gamma \setminus C_{\Gamma}  \right|, $$ and a maximal set $\{\zeta_\alpha \}_{\alpha = 1}^{\kappa} $ of representatives for the equivalence classes of cusps.  We set 
$ \Gamma_\alpha \df \Gamma_{\zeta_\alpha}$ and 
$ \Gamma_\alpha^\prime \df \Gamma_{\zeta_\alpha}^\prime. $
\enot
\noindent The following is elementary.
\bl \label{lemConjCusp}
Let $\zeta_\alpha$ be a cusp of $\Gamma.$ Then there exists $B_\alpha \in \pc $ and a lattice $\Lambda_\alpha = \ZZ \oplus \ZZ \tau_\alpha,~\I(\tau_\alpha) > 0 $ satisfying the following.
\ben
\item $\zeta_{\alpha} = B_{\alpha}^{-1} \infty,$
\item $B_{\alpha}\Gamma_{\alpha}B_{\alpha}^{-1}  $  acts discontinuously on $\cinf \setminus \{\infty\} = \CC$.
\item $$  B_{\alpha}\Gamma_{\alpha}^\prime B_{\alpha}^{-1} =  
\left\{ \,\left.\left(\begin{array}{cc}
1 &  b\\
0 &   1\end{array}\right)\,\right|\,0\neq b \in \Lambda_\alpha \,\right\}.  $$
\een
\el
\begin{nota}
For each cusp class $\alpha = 1\dots \kappa $ we fix $B_\alpha $ and $\Lambda_\alpha $ from Lemma~\ref{lemConjCusp} and set $ \Gamma_{\alpha_\infty} \df B_{\alpha}\Gamma_{\alpha}B_{\alpha}^{-1},  $ and $ \Gamma_{\alpha_\infty}^\prime \df B_{\alpha}\Gamma_{\alpha}^\prime B_{\alpha}^{-1}.$
\end{nota}
While the action of  $ \Gamma_{\alpha_\infty}^\prime $ on $\CC$ is exactly the action of the (additive) lattice $\Lambda_\alpha$ on\footnote{For $z\in \CC,$ we have $ \left(\begin{array}{cc}
1 &  b\\
0 &   1\end{array}\right)z = z+b.$ \\ } $\CC, $ the action of $ \Gamma_{\alpha_\infty} $ is a combination of the lattice action of  $\Lambda_\alpha$ and possibly  finite ordered euclidean rotations of $\CC.$

The fundamental domain for a cofinite Kleinian group can be realized as the union of a compact (hyperbolic) polyhedron and $\kappa$ \emph{cusp sectors}. Let $\pg_{\alpha} \subset \CC $ be a fundamental domain\footnote{The domain $\pg_{\alpha} $ is a euclidean polygon. \\} for the action of $\Gamma_{\alpha_\infty}$ on $\CC,$ and let $\pg_{\alpha}^\prime $ be the fundamental parallelogram with base point at the origin of the lattice $\Lambda_\alpha. $ For $Y>0$ set $$
 \widetilde{\F}_{\alpha}(Y) \df \{\, z+rj\,|\, z\in\pg_{\alpha},\, r\geq Y\,\}$$
and define the cusp sector,  $
\F_{\alpha}(Y) \df B_{\alpha}^{-1} \widetilde{\F}_{\alpha}(Y). $
\begin{lem}
\label{lemFunDom} \cite[Prop. 2.3.9]{Elstrodt} Let $\Gamma < \pc$ be cofinite with $\kappa = \left|\Gamma\setminus C_{\Gamma}\right|.$
Then there exist $Y>0$ and a compact set $\F_{Y}\subset\HH$
such that \[
\F \df \F_{Y}\cup\F_{1}(Y)\cup\cdots\cup\F_{\kappa}(Y)\]
 is a fundamental domain for $\Gamma.$ The compact set $\F_{Y}$
can be chosen such that  \mbox{$\F_{Y}\cap\F_{\alpha}(Y_{\alpha})$} 
all are contained in the boundary of  $\F_{Y}$  and hence have Lebesgue
measure 0.  Also, $\F_{\alpha}(Y_{\alpha})\cap\F_{\beta}(Y_{\beta})=\emptyset$
if $\alpha \neq \beta.$
\end{lem}

\section{Cuspidal elliptic elements}
Throughout this section $\Gamma $ is a cofinite Kleinian group.
\bl Let  $\Gamma$ be cofinite.  Then $\Gamma $  has only finitely many elliptic conjugacy classes in $\Gamma$. 
\el
\pf
Assume not.  Then there is an infinite sequence of elliptic $\Gamma-$conjugacy classes $ \{[e_n] \}. $  Next for each $n$ choose representative $e_n$ which fixes a point $P_n$ on the boundary of $\F,$ the fundamental domain of $\Gamma$ given in Lemma~\ref{lemFunDom}.  Since $\Gamma$ is discrete, the points \emph{must} accumulate at least at one cusp.  After conjugating $\Gamma$  and passing to a subsequence of $\{P_n \}$ we may assume that  $P_{n} \ra \infty.$ An application of \cite[Corollary 2.3.3]{Elstrodt} implies that $e_n \in \Gamma_\infty.$ In short we have constructed infinitely many elliptic $\Gamma-$conjugacy classes with a representative fixing the cusp $\infty.$  An elementary computation using  Lemma~\ref{lemCuspStructI} shows that there are only finitely many elliptic $\Gamma_\infty-$conjugacy classes, a contradiction.
\epf

\bnot For each $ \alpha \in \{ 1 \dots \kappa \} $ set 
$$ \Upsilon_{\alpha} \df  \left\{ \left(\begin{array}{cc}
\epsilon & 0\\
0 & \epsilon^{-1}\end{array}\right)~ \right. \left| ~ \left(\begin{array}{cc}
\epsilon & 0\\
0 & \epsilon^{-1}\end{array}\right) \in     \Gamma_{\alpha_\infty}  \right\}. $$ 
\enot 
If $$ \left(\begin{array}{cc}
\epsilon & 0\\
0 & \epsilon^{-1}\end{array}\right) \in \Upsilon_\alpha, $$ then for some non-zero natural number $N,$  $\epsilon^N  = 1. $ 
We have the following application of Lemma~\ref{lemCuspStructI}.
\bl
The set $\Upsilon_{\alpha}$ is isomorphic to a finite subgroup of the unit circle $ S^{1}$,  cyclic,  and isomorphic to $\Gamma_{\infty \alpha}/\Gamma'_{\alpha_\infty}.$ Any element $\gamma \in \Gamma_{\alpha_\infty}$ can be written uniquely in the form $
\gamma = \alpha \beta, $
 where $\alpha \in \Upsilon_{\alpha}$ and $\beta \in \Gamma_{\alpha_\infty}^{\prime}.$ \el

\bd An elliptic element $\gamma \in \Gamma$, is said to be 
\emph{cuspidal elliptic} if at least one\footnote{   Actually by \cite[Cor 2.3.11]{Elstrodt} \label{footCuspEllip},  if one fixed point is a cusp, then the other fixed point is also a cusp. } of its fixed points in 
$\cinf$ is a cusp of $\Gamma$. Otherwise it is called a  \emph{non-cuspidal elliptic}  element. The set of cuspidal elliptic elements is denoted by $\Gamma^{\CE}.$ 
\ed

\bnot
For $\alpha \in \{1\dots \kappa \} $ define $ \cuspi_\alpha $ to be the set of elements of $\Gamma$ 
which are $\Gamma$-conjugate to an element of $ \Gamma_\alpha
\setminus \Gamma_\alpha^\prime. $  
We fix representatives of $\cuspi_\alpha ~:
g_{1}^{\alpha}, \dots , g_{d_\alpha}^{\alpha} $
and define  $
q_{i}^{\alpha} \df  B_{\alpha} g_{i}^{\alpha} B_{\alpha}^{-1} $
\enot
Since $q_{i}^{\alpha}$ fixes $\infty, $  
\beq q_{i}^{\alpha} =  \left(\begin{array}{cc}
 \epsilon^{\alpha}_{i} & \epsilon^{\alpha}_{i}w^{\alpha}_{i}\\
 0 & \left(\epsilon^{\alpha}_{i}\right)^{-1}\end{array}
\right),
\eeq
where $\epsilon^{\alpha}_{i}$ is a finite-ordered root of unity and $w^{\alpha}_{i} \in \Lambda_\alpha.$

\chapter{The Spectral Decomposition Theorem}
In this chapter we prove the spectral decomposition theorem for \ckus. While much of our work is analogous to the \tdv and \thds cases\footnote{We encourage the reader to first read \cite[Sections 4.1,4.2,4.3,6.1,6.2,6.3,6.4]{Elstrodt} and \cite[Chapter 1,2,3]{Venkov} \\} there are some new complications that arise because of cuspidal elliptic elements\footnote{Fuchsian groups do not have cuspidal elliptic elements.}.  

\section{Unitary Representations} \label{secSingular}
Throughout this section $\Gamma $ is a cofinite Kleinian group with $\kappa = \left|\Gamma\setminus C_{\Gamma}\right| $ and cusp representatives  $ \{ \zeta_\alpha \}_{\alpha=1}^\kappa.$

Let $V$ be an $n$-dimensional complex vector space, with an inner product $\left< \cdot,\cdot \right>_V$ 
linear in the first argument and anti-linear in the second argument.
For $v \in V$ its norm $|v|_{V}$  is given by $\sqrt{\left< v,v \right>_V}$ and the norm of a linear operator $L:V\rightarrow V$ is defined as 
$$
|L|_{V} \df \sup_{v\in V}\left(\frac{|Lv|_{V}}{|v|_{V}}\right). $$
\bd
Define $\rp$ to be the set of pairs $(\chi,V)$ where $V$ a finite
dimensional complex inner product space and $\chi$ is a unitary representation
of $\Gamma$ in $\GL(V)$.
\ed   
We will abuse notation and identify $\chi $ with $(\chi,V). $

\bd
\label{def:Fixed V}
Let $\Gamma $ be cofinite, $\kappa = \left|\Gamma\setminus C_{\Gamma}\right| $  and $\chi\in\rep.$  
For each $ \alpha \in \{1,\cdots,\kappa\}$ define 
$$V_{\alpha}\df \{ v \in V \, | \, \chi(\gamma)v=v \, \, \, \forall \gamma \in \Gamma_{\alpha}\,\}, $$  
$$ V_{\alpha}^\prime \df \{ v \in V \, | \, \chi(\gamma)v=v \, \, \, \forall \gamma \in \Gamma_{\alpha}^\prime \,\},  $$
$$k_{\alpha} \df \dim V_{\alpha},$$
$$l_{\alpha} \df \dim V_{\alpha}^\prime, ~\text{and}$$
$$k \df k(\Gamma,\chi) \df \sum_{\alpha=1}^{\kappa} k_{\alpha}.$$
\end{defn}

The subspace $V_{\alpha}$ is called a \emph{singular} subspace while $V'_{\alpha}$ is called an \emph{almost singular}\footnote{ The notion of \emph{almost singularity} does not occur in the two dimensional case.} subspace,   $k_{\alpha}$ is called the \emph{degree of singularity} of the representation
$\chi$ at the cusp $\alpha$,  and $k=k(\Gamma,\chi)$  the degree
of singularity of $\Gamma$ relative to $\chi.$  
\bd
A representation $\chi$ is \emph{singular at} the cusp  $\zeta_{\alpha}$ if $k_{\alpha}>0,$  \emph{regular at $\zeta_{\alpha}$ } if $k_{\alpha}=0.$  It is    \emph{singular} if $k>0,$ and \emph{regular} if $k=0.$
\ed

\subsection{Unitary representions of the stabilizer subgroup}

\begin{lem}
\label{L:mylemma} 
There exist $E_\alpha,R_\alpha,S_\alpha \in \Gamma_{\alpha} $ with the following properties:

(1)  $\Gamma_{\alpha}=\{\, E_\alpha^{k}R_\alpha^{i}S_\alpha^{j}\,|\,0\leq k<m_\alpha,\, i,j\in\ZZ\,\}.$
Here $R_\alpha,S_\alpha $ are parabolic elements with $B_\alpha R_\alpha B_\alpha^{-1} (P) = P+1$ and $B_\alpha S_\alpha B_\alpha^{-1} (P) = P+\tau_\alpha$ (here $\Lambda_\alpha = \ZZ \oplus \ZZ \tau_\alpha $)  for all $P \in \HH,$  and  $E_\alpha$ is elliptic of order $m_\alpha.$ 

(2)  $\Gamma_{\alpha}^\prime =\{\, R_\alpha^{i}S_\alpha^{j}\,|\, i,j\in\ZZ\,\}.$

(3) The elements $R_\alpha$ and $S_\alpha$ commute but the group $\Gamma_{\alpha}$ is not abelian when $m_\alpha>1$. 

(4) If in addition, $m_\alpha > 1, $  
then  $\chi(E_\alpha)$ maps $V_\alpha^\prime $ onto itself. Furthermore, there exists
a basis of $V_\alpha^\prime$ so that $\chi(E_\alpha)|_{V_\alpha^\prime}$ is diagonal. 
\end{lem}
\begin{proof}
(1), (2), and (3) readily follow from \cite[Theorem 2.1.8]{Elstrodt}. We prove (4):
set $E = E_\alpha, R=R_\alpha, S=S_\alpha.$ Since $ERE^{-1}$ and $R^{-1}$ are both parabolic and in  $\Gamma_{\alpha}^\prime $ it follows that 
 the $A \df E^{-1}R^{-1}E R \in\Gamma_{\alpha}^\prime $ is parabolic.    Since $A$ is parabolic, $A=R^{i}S^{j}$
for some $i,j\in\ZZ$ (applying (2)). By definition, the restriction $\chi(A)|_{V_\alpha^\prime }=I_{V_\alpha^\prime }.$ 
We have  
\[ \chi(E)\chi(R)=\chi(R)\chi(E)\chi(A). \] Next applying an arbitrary $v\in V_\alpha^\prime $ we obtain 
\beq
\chi(E)v = \chi(E)\chi(R)v = \chi(R)\chi(E)\chi(A)v = \chi(R)\chi(E)v 
\eeq
Thus $\chi(R)$ fixes $\chi(E)v$ and similarly $\chi(S)$ fixes
$\chi(E)v,$ so by definition  $\chi(E)v\in V_\alpha^\prime.$  Since $\chi(E)$ is
 unitary, its restriction to $V_\alpha^\prime$ is also unitary, hence
$V_\alpha^\prime$ has a diagonalizing basis. 
\end{proof}

We remark that for a cusp $\alpha $ of $\Gamma, $ the group $\Gamma_\alpha^\prime $ is abelian.  Thus, the unitary representations $\chi $ restricted to $\Gamma_\alpha^\prime $ is diagonalizable and can be thought of as a direct sum of one-dimensional unitary representations of $\Gamma_\alpha^\prime. $

\section{Automorphic Functions }
The set of $\Gamma$-automorphic functions 
$A(\Gamma, \chi)$ is the set of all Borel-measurable functions $f:\HH \ra V$ that satisfy  $ f(\gamma P)= \chi(\gamma)f(P)$ for all  $\gamma \in \Gamma.$  Such functions are uniquely determined by their values on a 
fundamental domain. With this in mind we fix a fundamental domain 
$\F\subset\HH$ for $\Gamma$ and define 
$$\hs \df  L_{2}(\F,V,dv,\chi) \df \{ f \in A(\Gamma, \chi) ~|~  
\int_{\F}\langle f(P),f(P)\rangle_{V}\, dv(P)<\infty  \}. $$  
For $f,g\in\hs$ define the inner product $ \left< f,g \right > \df \int_{\F}\langle f(P),g(P)\rangle_{V}\, dv(P). $  With this inner product  $\hs$ is a Hilbert Space.  For $n \in \NN $ let  $C^{n}(\HH,V,\chi)$ be  space of $\Gamma$-automorphic functions which 
are $n$-times differentiable on $\HH$.

\subsection{The Automorphic Laplacian $\lp(\Gamma, \chi)$ }
Recall that the Laplace-Beltrami operator $\lp $ is defined on twice continuously differentiable functions $f:\HH \ra \CC. $  Thus $\lp $ can be defined on the set $A(\Gamma, \chi) \cap C^{2}(\HH,V,\chi),$ and for $f \in A(\Gamma, \chi) \cap C^{2}(\HH,V,\chi)$ it follows (since $\lp $ commutes with isometries) that $\lp f \in A(\Gamma, \chi)\cap C^{2}(\HH,V,\chi).$ 

 Consider the dense subspace  $D = \{ f \in \hs \cap C^{2}(\HH,V,\chi) \, | \, \lp f \in \hs \, \}.$  It follows that $\lp:D\ra \hs $ is essentially self-adjoint and has a unique positive self-adjoint extension to a space $ \widetilde{D}$ (see \cite[Section 4.1]{Elstrodt}).

\begin{defn}
Let  $\Gamma $ be a cofinite Kleinian group and $ \chi \in \rp. $
We define the automorphic Laplacian 
$$\lp(\Gamma, \chi) : \widetilde{D} \ra \hs $$ 
to be the self-adjoint extension of $\lp:D\ra \hs. $ 
\end{defn}
When there is little possibility for confusion, we will identify $\lp $ with $\lp(\Gamma, \chi).$  The dense subspace $ \widetilde{D} $ can be realized as the set of functions in $\hs $ whose \emph{distributional} Laplacian is in $\hs. $

\section{Fourier series expansion at a cusp}
\subsection{The scalar case} \label{secScalExpan}
Let $f \in A(\Gamma, 1) $ (that is $f:\HH \ra \CC, $ and for all $P \in \HH,~\gamma \in \Gamma,$  $f(\gamma P) = f(P)$), and let $\zeta_\alpha $ be a cusp of $\Gamma $ with $B_\alpha^{-1} \infty = \zeta_\alpha. $ Then it is easily seen that  $f \circ B_\alpha^{-1} \in A(B_\alpha \Gamma B_\alpha^{-1}, 1), $ in particular $ f \circ B_\alpha^{-1} \in A( B_{\alpha}\Gamma_{\alpha}^\prime B_{\alpha}^{-1},1).$   Since   $ B_{\alpha}\Gamma_{\alpha}^\prime  B_{\alpha}^{-1}$ is isomorphic to the lattice $\Lambda_\alpha, $  $~ f \circ  B_\alpha^{-1} (P + \omega) = f \circ B_\alpha^{-1}(P) $ for all $P \in \HH $ and $ \omega \in \Lambda_\alpha. $  That is, $ f \circ  B_\alpha^{-1} $ is invariant under that lattice $ \Lambda_\alpha. $  If we suppose further that $f$ is smooth then we can expand $f \circ B_\alpha^{-1} $ as
$$
f(B_\alpha^{-1}(z+rj)) = \sum _{\mu \in \Lambda_\alpha ^{0}}g_{\mu }(r)e^{2\pi i \left<\mu ,z \right>}.
$$ 
Here $ \left<~,~\right> $ is the (real) euclidean inner product in $\RR^2 = \CC, $  $ \Lambda_\alpha^0 $ is the lattice dual to $ \Lambda_\alpha.$  If we assume still further that $\lp f = \lambda  f,  $ and $ f(B_\alpha^{-1}(z+rj)) = O(r^N) $ for some $N \in \NN, $ then a simple  separation of variables argument (see \cite[Section 3.3]{Elstrodt}) shows that for $\lambda = 1-s^2,~s \neq 0,$  $g_\mu $ satisfies the Bessel equation  
\beq  \label{E:bessel}  \left( r^2 \frac{d^2}{dr^2} -r \frac{d}{dr} + \lambda - 4\pi^2 |\mu|^2 r^2 \right) g_\mu (r) = 0,
\eeq
 whose general solution is $g_\mu (r) = a_\mu r K_s(2 \pi |\mu| r) + b_\mu r I_s(2 \pi |\mu| r). $  Here 
$$  I_s(w) = \sum_{m = 0}^{\infty} \frac{w^{s+2m}}{m! 2^{s + 2m} \Gamma(s+m+1)}, \, \, \, \, \, K_s (w) = \frac{\pi}{2 \sin {s \pi} } \left(I_{-s}(w) - I_s (w) \right).  
$$
The function $ K_s(x) $ decreases exponentially, and  $ I_s(x) $ increases exponentially as $ x \ra \infty $.  Applying the growth bound of $f$ we obtain
\beq \label{E:scalarexpand}
f(B_\alpha^{-1}(z+rj)) = a_0 r^{1+s} + b_0 r^{1-s} + \sum_{ 0 \neq \mu \in \Lambda_\alpha^0} a_\mu r K_s(2 \pi |\mu| r) e^{2 \pi i \left<\mu,z \right >}.    \eeq

We will need the following fact later on: 
\bl \label{lemBesSumDec} 
Let $p$ be a non-negative integer.  Then 
$$ \frac{d^p}{dr^p}\sum_{ 0 \neq \mu \in \Lambda_\alpha^0} a_\mu r K_s(2 \pi |\mu| r) e^{2 \pi i \left<\mu,z \right >} = O(e^{-|c|r}) ~ \text{as}~ r \ra \infty,   $$ for some constant  $|c| > 0.$
\el
\subsection{Preliminary Lemmas} \label{secForExp}
Our goal in this section is to come as close as possible to diagonalizing the group of unitary transformations $ \{ \chi( \Gamma_\alpha ) \}. $
Recall that 
$$V_{\alpha}\df \{ v \in V \, | \, \chi(\gamma)v=v \, \, \, \forall \gamma \in \Gamma_{\alpha}\,\}, $$  
$$ V_{\alpha}^\prime \df \{ v \in V \, | \, \chi(\gamma)v=v \, \, \, \forall \gamma \in \Gamma_{\alpha}^\prime \,\},  $$
$$k_{\alpha} \df \dim V_{\alpha}, ~ \text{and} $$ 
$$l_{\alpha} \df \dim V_{\alpha}^\prime. $$

For each cusp $\zeta_\alpha $ of $\Gamma $ let $(I - P_\alpha) :V \ra V_\alpha^\perp $ denote the orthogonal projection onto $V_\alpha^\perp. $    By definition of singularity, for each $0 \neq v \in V_\alpha^\perp $ there exists $\gamma \in \Gamma_\alpha $ satisfying $\chi(\gamma) v \neq v.$

Fix a cusp $\zeta_\alpha $ of $\Gamma. $  We can partition $n = \dim V $ as follows:
$$1,\dots, k_\alpha, k_{\alpha + 1}, \dots, l_\alpha, l_{\alpha + 1}, \dots , n.  $$  We will next chose a particular basis of $V$ that respects the partions above and is a close as possible to diagonalizing $ \{ \chi( \Gamma_\alpha ) \}. $

Recall from Lemma \eqref{L:mylemma} that 
$\Gamma_{\alpha}=\{\, E_\alpha^{k}R_\alpha^{i}S_\alpha^{j}\,|\,0\leq k<m_\alpha,\, i,j\in\ZZ\,\}. $  Here $R_\alpha, S_\alpha $ are parabolic and $E_\alpha$ is elliptic.

Let $B_{s \alpha } \df  \{  v_{\alpha 1}, \dots, v_{ \alpha  k_\alpha } \} $ be an orthonormal basis for $V_\alpha. $   By lemma \eqref{L:mylemma} there exists a set  $B_{a \alpha} \df  \{ v_{\alpha (k_\alpha + 1)}, \dots, v_{\alpha l_\alpha} \} $ of orthonormal elements in $V,$ pairwise orthogonal to the elements of $B_{s \alpha}$  so that the following conditions are satisfied:  For each $v \in B_{a \alpha},  \chi(\gamma)v = v $ for all $\gamma \in \Gamma_\alpha^\prime, $ and $\chi(E_\alpha) v = \lambda v, $ where here $ E_\alpha $ is the primitive elliptic element in $\Gamma_\alpha $ chosen  in lemma \eqref{L:mylemma}, and $1 \neq \lambda \in \CC $ with $  \lambda^{m_\alpha} = 1. $  Finally, there exists a set $B_{r \alpha} \df  \{v_{\alpha (l_\alpha + 1) }, \dots, v_{\alpha n} \} $ of orthonormal elements of $V,$ pairwise orthogonal to the elements of $B_{s \alpha} \cup B_{a \alpha}$ so that the following conditions are satisfied:  For each $v \in B_{r \alpha}, \chi(\gamma) v = \lambda' v $ for each $\gamma \in \Lambda_\alpha^\prime $ and some $1 \neq \lambda' \in \CC $ with $ | \lambda' | = 1. $  We have proven the following.

\begin{lem} \label{lemChoBas}
For each cusp $\zeta_\alpha $ let $\Gamma_{\alpha}=\{\, E_\alpha^{k}R_\alpha^{i}S_\alpha^{j}\,|\,0\leq k<m_\alpha,\, i,j\in\ZZ\,\} $ (see lemma \eqref{L:mylemma}). Then  there exists  an orthonormal basis  
$$ \{  v_{\alpha 1}, \dots, v_{ \alpha  k_\alpha }, v_{\alpha (k_\alpha + 1)}, \dots, v_{\alpha l_\alpha}, v_{\alpha (l_\alpha + 1) }, \dots, v_{\alpha n} \} =  B_{s \alpha} \cup B_{a \alpha} \cup B_{r \alpha}  $$ of $V$ 
with the following properties:

(1) For each $v_{\alpha l} \in B_{s \alpha}~( 1 \leq l \leq k_\alpha), $ $\chi(\gamma) v_{\alpha l} = v_{\alpha l} $ for all $\gamma \in \Gamma_\alpha. $  

(2)  For each $v_{\alpha l} \in B_{a \alpha}~( k_\alpha + 1  \leq l \leq l_\alpha), $ $\chi(R_\alpha) v_{\alpha l} = \chi(S_\alpha) v_{\alpha l} = v_{\alpha l}, $ and $\chi(E_\alpha)v_{\alpha l} = \lambda_{\alpha l}v_{\alpha l}. $ Here $ 1 \neq \lambda_{\alpha l} \in \CC $ and $\lambda_{\alpha l}^{m_\alpha}  = 1.$

(3) For each $v_{\alpha l} \in B_{r \alpha}~( l_\alpha + 1  \leq l \leq n), $ $\chi(R_\alpha) v_{\alpha l} =  \lambda_{R_\alpha  l } v_{\alpha l}, $ and  $\chi(S_\alpha) v_{\alpha l} =  \lambda_{S_\alpha  l } v_{\alpha l}. $  Here $  \lambda_{R_\alpha  l }, \lambda_{S_\alpha  l } \in \CC, $  with $| \lambda_{R_\alpha  l }| = | \lambda_{S_\alpha  l }| = 1,$ and    $ \lambda_{R_\alpha  l }, \lambda_{S_\alpha  l } $ are not both equal to one. 
\end{lem}

For each cusp $\zeta_\alpha $ and $l_\alpha + 1  \leq l \leq n $ fix $\theta_{R_\alpha  l } $ and $\theta_{S_\alpha  l } $ so that 
\beq
e^{ 2 \pi i \theta_{R_\alpha  l } } = \lambda_{R_\alpha  l} ~ \text{and} ~  e^{ 2 \pi i \theta_{S_\alpha  l } } = \lambda_{S_\alpha  l}. 
\eeq

We remark that $ B_{s \alpha} $ and $ B_{r \alpha} $ are analogous to the two dimensional $\emph{singular} $ and regular  cases respectively (see \cite{Venkov} ).  We call $ B_{a \alpha} $ the \emph{almost singular} case.

Let $P_{\alpha l} $ be the orthogonal projection onto the subspace generated by $v_{\alpha l}. $
\begin{lem} \label{L:autolemma}
Let $f \in A(\Gamma, \chi ).  $ 

(1) If  $ 1 \leq l \leq l_\alpha, $ then  
$$P_{\alpha l}f( B_\alpha^{-1}(Q + \omega)) = P_{\alpha l}f( B_\alpha^{-1} Q ) $$ 
for all $Q \in \HH$ and $ \omega \in \Lambda_\alpha. $ 

(2) If $ l_\alpha + 1 \leq l \leq  n, $ then for $n,m \in \ZZ, $
$$P_{\alpha l}f( B_\alpha^{-1}(Q + m+ n \tau_\alpha )) =e^{2 \pi i  (m \theta_{R_\alpha  l }+ n \theta_{S_\alpha  l }) }     P_{\alpha l}f( B_\alpha^{-1} Q ), $$
where $\ZZ \oplus \ZZ \tau_\alpha = \Lambda_\alpha. $ 

\end{lem}
\pf
(1) Let $\omega \in \Lambda_\alpha.$  Since $ \Lambda_\alpha $ is isomorphic to   $ B_{\alpha}\Gamma_{\alpha}^\prime  B_{\alpha}^{-1}$ there exists $\gamma \in \Gamma_{\alpha}^\prime $ satisfying $ B_{\alpha}\gamma B_{\alpha}^{-1}Q = Q + \omega $ for all $Q \in \HH $ (See lemma \eqref{L:mylemma}).   Since $f \in A(\Gamma, \chi), $ it follows that   $P _{\alpha l}f(\gamma Q) =P _{\alpha l}\chi(\gamma) f(Q).$   However, $$   P _{\alpha l} \chi(\gamma) =  (\chi(\gamma)^* P _{\alpha l}^*)^* = (\chi(\gamma^{-1})P _{\alpha l})^* =   (P_{\alpha l})^* = P_{\alpha l}.   $$  Thus, 
$$P_{\alpha l}f (B_\alpha^{-1}(Q + \omega)) = P_{\alpha l}f (B_\alpha^{-1} B_{\alpha}\gamma  B_{\alpha}^{-1}(Q))  =  P_{\alpha l}f( \gamma  B_{\alpha}^{-1}Q) =P_{\alpha l}f( B_{\alpha}^{-1}Q).   $$

(2) A direct calculation similar to (1) shows that $P_{\alpha l}f( B_\alpha^{-1}(Q + m+ n \tau_\alpha )) = P_{\alpha l}f(R_\alpha^m S_\alpha^n B_\alpha^{-1}Q) = P_{\alpha l} \chi(R_\alpha)^m \chi(S_\alpha)^n f(B_\alpha^{-1}Q)= e^{2 \pi i  (m \theta_{R_\alpha  l }+ n \theta_{S_\alpha  l }) }P_{\alpha l} f(B_\alpha^{-1}Q).$ 
\epf
\subsection{Fourier Expansion: The vector case} \label{S:fourvec}
\bp \label{P:vectorexpan}
Let $\Gamma $ be cofinite,  $f \in A(\Gamma, \chi ) $ be  smooth satisfying \linebreak $ f(B_\alpha^{-1}(z+rj)) = O(r^N)$  for some $N \in \NN $ and  each cusp $\zeta_\alpha $ of $\Gamma,  $ and  $ \lp f = \lambda f  $ with $ \lambda = 1-s^2,~s \neq 0. $ 

(1) If  $ 1  \leq l \leq k_\alpha, $ then
$$P_{\alpha l}f( B_\alpha^{-1}(z + rj)) =  a_{\alpha l, 0} r^{1+s} + b_{\alpha l, 0} r^{1-s} + \sum_{ 0 \neq \mu \in \Lambda_\alpha^0} a_{\alpha l, \mu} r K_s(2 \pi |\mu| r) e^{2 \pi i \left<\mu,z \right >}.    $$

(2) If $ k_\alpha + 1  \leq l \leq l_\alpha, $ then
$$P_{\alpha l}f( B_\alpha^{-1}(z + rj)) =  \sum_{ 0 \neq \mu \in \Lambda_\alpha^0} a_{\alpha l, \mu} r K_s(2 \pi |\mu| r) e^{2 \pi i \left<\mu,z \right >}. $$

(3) If $ l_\alpha + 1  \leq l \leq n, $  $\theta_{R_\alpha  l },\theta_{S_\alpha  l } $ are chosen so that  $ \exp( 2 \pi i \theta_{R_\alpha  l }) = \lambda_{R_\alpha  l },$ and $ \exp( 2 \pi i \theta_{S_\alpha  l }) = \lambda_{S_\alpha  l },$
then 
$$P_{\alpha l}f( B_\alpha^{-1}(z + rj)) =  \sum_{  \mu \in \Lambda_\alpha^0} a_{\alpha l, \mu} r K_s(2 \pi |\mu'| r) e^{2 \pi i \left<\mu',z \right >}. $$
where $\mu '=\theta _{R}\mu _{1}+\theta _{S} \mu_{2}+\mu,$  $ \Lambda_\alpha^0 $ is the lattice dual to  $ \Lambda_\alpha$ generated by $\mu _{1},\mu _{2},$ where $$ ~\left< \mu_1, 1 \right > = 1,\left< \mu_1, \tau_\alpha  \right > = 0,  \left< \mu_2, 1 \right > = 0, \left< \mu_2, \tau_\alpha \right > = 1, $$ and the $ a_{ \alpha l, \mu }$ are scalar multiples of the vector $v_{\alpha l}.$
\ep

\pf
\noindent (1) By Lemma \ref{L:autolemma} $P_{\alpha l}f( B_\alpha^{-1}(z + rj)) $ is a scalar function multiplied by the vector $v_{\alpha l}, $ thus an application of equation \eqref{E:scalarexpand} concludes the proof.

\noindent(2) By (1) $$P_{\alpha l}f( B_\alpha^{-1}(z + rj)) =  a_{\alpha l, 0} r^{1+s} + b_{\alpha l, 0} r^{1-s} + \sum_{ 0 \neq \mu \in \Lambda_\alpha^0} a_{\alpha l, \mu} r K_s(2 \pi |\mu| r) e^{2 \pi i \left<\mu,z \right >}.    $$ We will show that $a_{\alpha l, 0}$ and  $b_{\alpha l, 0} $ are both zero.  To see this first, observe that (by \mbox{Lemma \ref{lemBesSumDec}}) 
$$ 
P_{\alpha l}f( B_\alpha^{-1}(z + rj)) =  a_{\alpha l, 0} r^{1+s} + b_{\alpha l, 0} r^{1-s} + \sum_{ 0 \neq \mu \in \Lambda_\alpha^0} a_{\alpha l, \mu} r K_s(2 \pi |u| r) e^{2 \pi i \left<u,z \right >} 
$$
$$
 =a_{\alpha l, 0} r^{1+s} + b_{\alpha l, 0} r^{1-s} +O(e^{-|c|r}) ~\text{as}~ r \ra \infty.
$$
Next, $$ P_{\alpha l}f( B_\alpha^{-1} (B_\alpha E_\alpha B_\alpha^{-1})  (z + rj ) ) = a_{\alpha l, 0} r^{1+s} + b_{\alpha l, 0} r^{1-s} +O(e^{-|c|r}) ~\text{as}~ r \ra \infty
$$
since $ B_\alpha E_\alpha B_\alpha^{-1} (z + rj) = z' + rj $ for some $z' \in \ZZ, $  that is $B_\alpha E_\alpha B_\alpha^{-1} $ fixes  the $r-$coordinate in $\HH $  (this follows from elementary computations and that fact that $B_\alpha E_\alpha B_\alpha^{-1} \infty = \infty $).  However, we also have 
$$ P_{\alpha l}f( B_\alpha^{-1} (B_\alpha E_\alpha B_\alpha^{-1})  (z + rj ) ) = P_{\alpha l} \chi(E_\alpha)f( B_\alpha^{-1}(z + rj) ) = $$ 
$$\lambda_{\alpha l} P_{\alpha l}f( B_\alpha^{-1}(z + rj ) ) = \lambda_{\alpha l}(  a_{\alpha l, 0} r^{1+s} + b_{\alpha l, 0} r^{1-s} +O(e^{-|c|r})), 
$$
a contradiction since $\lambda_{\alpha l} \neq 1. $

(3) Let $h(Q) \df h(z) \df  \exp( \left< \mu_1, z  \right> \theta_{ R \alpha l } + \left< \mu_2, z  \right> \theta_{ S \alpha l } ), $ where $Q = z + rj.$  Then $h(Q) $ satisfies $$ h(Q + m + \tau_\alpha n) = \exp(2 \pi i (m \theta_{ R \alpha l } +  n \theta_{ S \alpha l }) )h(Q),  $$ thus by lemma \eqref{L:autolemma} 
$$\frac{P_{\alpha l}f \circ B_\alpha^{-1}}{h}(Q + \omega)) = \frac{P_{\alpha l}f \circ B_\alpha^{-1} }{h}(Q)  $$ 
for all $Q \in \HH$ and $ \omega \in \Lambda_\alpha. $  Thus we can write  $$
\frac{P_{\alpha l}f \circ B_\alpha^{-1} }{h}(Q) = \sum _{\mu \in \Lambda_\alpha ^{0}}g_{\mu }(r)e^{2\pi i \left<\mu ,z \right>},
$$ which is equivalent to 
$$ P_{\alpha l}f \circ B_\alpha^{-1}(Q) =  \sum _{\mu \in \Lambda_\alpha ^{0}}g_{\mu }(r)h(z) e^{2\pi i \left<\mu ,z \right>}.   $$  To conclude the proof recall that $ \lp f = \lambda f, ~f $ grows by at most $O(r^N) $ near each cusp,   and apply the separation of variables technique (see \cite[page 105]{Elstrodt}). 
\epf

\section{The Eisenstein series}
As usual for this section let $\Gamma $ be cofinite with cusp representatives $ \{ \zeta_\alpha \}_{\alpha = 1}^{\kappa},  $  $B_{\alpha}^{-1}\infty=\zeta_{\alpha}, $ and $(\chi,V) \in \rp. $
\subsection{Definitions}
For  $P \in \HH,\, \R(s)>1, ~\text{and}~ v \in V_\alpha$ we define the  \emph{Eisenstein series}  by 
$$
E(P,s,\alpha ,v) \df  E(P,s,\alpha ,v,\Gamma,\chi) \df  \sum_{M \in \Gamma_\alpha \setminus \Gamma } \left(r(B_{\alpha}MP)
\right)^{1+s}\chi(M)^{*}v. $$
The series $E(P,s,\alpha ,v) $ converges uniformly and absolutely on compact subsets of  $\{ \R(s)> 1 \} \times \HH, $ and is a $\chi-$automorphic function   that  satisfies  $$ \lp E(~\cdot~,s,\alpha ,v) = \lambda E(~\cdot~,s,\alpha ,v). $$  In addition $$ E(B_\beta^{-1}(z+rj),s,\alpha ,v) = O(r^N)  $$ for any cusp $\zeta_\beta.$ 
\subsection{Fourier expansion of the Eisenstein series}
In this section we give the explicit Fourier expansion  of the Eisenstein series. 

For $\zeta_\alpha $ a cusp of $\Gamma $ and $1 \leq l \leq k_\alpha,  $ set  $E_{\alpha l} \df  E_{\alpha l}(s) \df   E(\,\cdot\,,s,\alpha,v_{\alpha l})$

\begin{lem}
Let $\zeta_\alpha, \zeta_\beta $ be cusps of $\Gamma, $ $1 \leq l \leq k_\alpha,  $ and $1 \leq k \leq k_\beta.  $  Then 
$$
P_{\beta} E_{\alpha l}(B_\beta^{-1}(z +rj)) = \delta_{ \alpha \beta } r^{1+s} v_{\beta l} + $$
$$ 
\frac{\pi }{|\Gamma_\alpha : \Gamma_\alpha^\prime  |  |\Lambda_\beta |^{s}}\left(\sum _{M \in \Gamma_{\alpha }^\prime \setminus \Gamma / \Gamma_{\beta }^\prime } P_\beta \chi ^{*}(M)v_{\alpha,l} |c|^{-2-2s}\right) r^{1-s} + 
 $$
$$
g(s)\sum _{0\neq |\mu |\in \Lambda_\beta^{0}}|\mu |^{s}\left(\sum _{M \in \Gamma_{\alpha}^\prime \setminus \Gamma / \Gamma_{\beta}^\prime } P_\beta \chi ^{*}(M)v_{\alpha,l} \frac{e^{2\pi i<\mu ,\frac{d}{c}>}}{|c|^{2+2s}}\right)rK_{s}(2\pi |\mu |r)e^{2\pi i<\mu ,z>} $$
where $$ g(s)\df \frac{2\pi ^{1+s}}{|\Gamma_\alpha : \Gamma_\alpha^\prime  | |\Lambda_\beta |\Gamma (1+s)},  $$

$c=c(B_\alpha M B_\beta^{-1}) ~ \text{and}~ d=d(B_\alpha M B_\beta^{-1})$  are the bottom left and bottom right indices respectively of the two-by-two matrix $B_\alpha M B_\beta^{-1}. $ 
\end{lem}
\pf
An application of Proposition~\ref{P:vectorexpan} part (1)  gives us the general form of the Fourier series expansion.  Then apply the  argument found in \cite[page 111]{Elstrodt}.
\epf

A \emph{Maa\ss~form} is a smooth function $f \in A(\Gamma, \chi ) $ satisfying $ f(B_\alpha^{-1}(z+rj)) = O(r^N)$  for some $N \in \NN $ and for each cusp $\zeta_\alpha $ of $\Gamma,  $ and  $ \lp f = \lambda f  $ with $ \lambda = 1-s^2,~s \neq 0. $   

We saw in proposition \eqref{P:vectorexpan} that for 
$ 1  \leq l \leq k_\alpha, $  a Maa\ss~form $f$ can be expanded as
$$P_{\alpha l}f( B_\alpha^{-1}(z + rj)) =  a_{\alpha l, 0} r^{1+s} + b_{\alpha l, 0} r^{1-s} + \sum_{ 0 \neq \mu \in \Lambda_\alpha^0} a_{\alpha l, \mu} r K_s(2 \pi |\mu| r) e^{2 \pi i \left<\mu,z \right >}.    $$  
The $r^{1+s}, r^{1-s} $ terms will play a promanent role in the what follows.  Define $ u_{f_{ \alpha l }} $ by 
$$ u_{f_{ \alpha l }}( B_\alpha^{-1}(z + rj)) =  a_{\alpha l, 0} r^{1+s} v_{\alpha l} + b_{\alpha l, 0} r^{1-s} v_{\alpha l} .    $$

\subsection{The scattering matrix} \label{secScaMat}
Since the Eisenstein series is a Maa\ss~form, we can apply the function $u$ to read off its  constant term in its Fourier expansion.  Define the \emph{scattering matrix}  $\smat_{\alpha ,l,\beta ,k}(s) $ by 
$$ u_{{\left( E_{\alpha l} \right) }_{\beta k}}(B_{\beta }^{-1}(z+rj)) = 
\delta_{\alpha \beta }\delta_{lk}r^{1+s}v_{\beta k} 
+ \smat_{\alpha ,l,\beta ,k}(s) r^{1-s}v_{\beta k}.   $$

We can put the scattering matrix   into a $k(\Gamma,\chi)
\times k(\Gamma,\chi)$ matrix. Let $\smat(s)$ be the $k(\Gamma,\chi)\times k(\Gamma,\chi)$ matrix given  block-wise,
$$(\smat_{\alpha \beta}(s)) \df  (\smat_{\alpha,l,\beta,k}(s)) $$  where $1 \leq l \leq k_\alpha$, $1 \leq k \leq k_\beta$.

\section{The Maa\ss-Selberg Relations}
Our main goal in this section is to prove the Maa\ss-Selberg relations in the case of cofinite Kleinian groups with finite dimensional unitary representations. The relations will be needed to prove the Selberg trace formula.   
For this section assume $\Gamma $ is a cofinite Kleinian group and $(\chi, V) \in \rp. $ 
\subsection{Greene's Theorem}
The basic tool for proving the Maa\ss-Selberg relations is Greene's Theorem.   Let $D \subset \HH $ be compact with a piece-wise smooth boundary, $f,g$ be complex valued smooth functions, $dv$ the volume form on $\HH, $ and $\widetilde{dv}$ the induced volume form on the boundary of $D,$ $\partial D. $  Then Greene's theorem is 
$$
\int _{D}(f\lp \bar{g}-\bar{g}\lp f)\, dv=\int _{\partial D}(f\frac{\partial }{\partial n}\bar{g}-\bar{g}\frac{\partial }{\partial n}f)\, \widetilde{dv}. $$

Our application of  Greene's theorem will be to the compact part of the fundamental domain $\F $ of $\Gamma. $  
Writing Greene's theorem in the global coordinates of $\HH $ 
using $dv=\frac{dx\wedge dy\wedge dr}{r^{3}},$ $\frac{\partial }{\partial n}=r\frac{\partial }{\partial r},$
$\widetilde{dv}=\frac{\partial }{\partial n}\perp dv=\frac{dx\wedge dy}{r^{2}}$
and Greene's Theorem becomes, \[
\int _{D}(f\lp \bar{g}-\bar{g}\lp f)\, \frac{dxdydr}{r^{3}}=\int _{\partial D}(f\frac{\partial }{\partial r}\bar{g}-\bar{g}\frac{\partial }{\partial r}f)\, \frac{dxdy}{r}\]
when $\partial D$ is a \td region parallel to $\CC $ on the
boundary of $\HH .$

\subsection{The scalar Maa\ss-Selberg relations}
For this section, assume $\Gamma$ has only one class of cusps and that $\chi = 1.$
Let $f \in A(\Gamma, 1) $ and assume that $f$ has polynomial growth as $z+rj $ approaches the cusp at $\infty, $ is smooth,  and satisfies $\lp f = \lambda f $ (in other words $f$ is a Maa\ss~form).  Let $\Lambda $ be the lattice corresponding to $\Gamma_\infty^\prime.$ Then we can expand   
\[
f(z+rj)=\sum _{\mu \in \Lambda _{0}}g_{\mu }(r)e^{2\pi i<\mu ,z>}\]
where  $$
(r^{2}\frac{d^{2}}{dr^{2}}-r\frac{d}{dr}+\lambda -4\pi ^{2}|\mu |^{2}r^{2})g_{\mu }(r)=0.  $$

\begin{lem}
\label{lem:ms1}Let $\lp f=\lambda f$, $\lp g=\nu g$ with $\lambda =1-s^{2},\nu =1-t^{2},s\neq 0,t\neq 0.$
Assume $\Gamma $ has only one cusp at $\infty $ and \[
f(P)=ar^{1+s}+br^{1-s}+\sum _{\mu \in \Lambda _{0}}f_{\mu }(r)e^{2\pi i<\mu ,z>}\]
\[
g(P)=cr^{1+t}+dr^{1-t}+\sum _{\mu \in \Lambda _{0}}g_{\mu }(r)e^{2\pi i<\mu ,z>}\]
 Let $u_{f}(P)=ar^{1+s}+br^{1-s}, $ $u_{g}(P)=cr^{1+t}+dr^{1-t}$
\[
h^{Y}(P)=\left\{ \begin{array}{cc}
 h(P)-u_{h}(P), & P\in \F (Y)\\
 h(P), & P\in \F _{0}\end{array}
\right.\]
for $h=f,g.$ Then \[
(s^{2}-\bar{t}^{2})\int _{\F }f^{Y}\bar{g}^{Y}\, dv=\frac{|\Lambda |}{|\Gamma _{\infty }:\Gamma '_{\infty }|}\left((s-\bar{t})a\bar{c}Y^{s+\bar{t}}+(s+\bar{t})a\bar{d}Y^{s-\bar{t}}\right)\]

\end{lem}
\[
-\frac{|\Lambda |}{|\Gamma _{\infty }:\Gamma '_{\infty }|}\left((s+\bar{t})b\bar{c}Y^{-s+\bar{t}}+(s-\bar{t})b\bar{d}Y^{-s-\bar{t}}\right)\]
 
\begin{proof}

Let $D = \F_0 = \F \setminus \F(Y) $ (here $ \F(Y)$ is the cusp sector that goes to infinity). Then by Greene's theorem 
$$
\int _{D}(f\lp \bar{g}-\bar{g}\lp f)\, dv=\int _{\partial D}(f\frac{\partial }{\partial n}\bar{g}-\bar{g}\frac{\partial }{\partial n}f)\, \widetilde{dv}. $$
Since the (vertical) sides of $D$ are pair-wise identified by isometries (with opposite orientation) the only surviving boundary term is the horizontal cross section (the boundary where $\F $ was separated into $\F_0 $ and $\F(Y)$).  Greene's theorem becomes   
\beq  \label{E:ms1}
\int _{D}(f\lp \bar{g}-\bar{g}\lp f)\, \frac{dxdydr}{r^{3}}=\int _{P}(f\frac{\partial }{\partial r}\bar{g}-\bar{g}\frac{\partial }{\partial r}f)\, \frac{dxdy}{r}, \eeq
where $dv=\frac{dx\wedge dy\wedge dr}{r^{3}},$ $\frac{\partial }{\partial n}=r\frac{\partial }{\partial r},$
$\widetilde{dv}=\frac{\partial }{\partial n}\perp dv=\frac{dx\wedge dy}{r^{2}},$ and
$P \df \partial D$ is the horizontal boundary region described above.  

The set $P$ is a vertical translation of a fundamental domain of  the action of $\Gamma_\infty $ on $\CC. $ That is the $r-$coordinate projection is constant and equal to $Y.$  Set    $$P' \df \bigcup _{ \gamma \in \Gamma _{\infty }^\prime  \setminus  \Gamma_{\infty}} \gamma P. $$ Then $P'$ is a fundamental domain for $\Lambda = \Lambda_\infty.$

With our notation the right side of equation \eqref{E:ms1}  becomes 
$$
\frac{1}{|\Gamma _{\infty }:\Gamma '_{\infty }|}\int _{P'}((\left(\sum _{\mu \neq 0}f_{\mu }(Y)e^{2\pi i<\mu ,z>}\right)\left(\sum _{\mu \neq 0}\bar{g}'_{\mu }(Y)e^{-2\pi i<\mu ,z>}\right) $$
$$
-\left(\sum _{\mu \neq 0}\bar{g}_{\mu }(Y)e^{-2\pi i<\mu ,z>}\right)(\left(\sum _{\mu \neq 0}f'_{\mu }(Y)e^{2\pi i<\mu ,z>}\right))\, \frac{dxdy}{Y} $$
 $$
+\frac{|\Lambda |}{|\Gamma _{\infty }:\Gamma '_{\infty }|}\left((s-\bar{t})a\bar{c}Y^{s+\bar{t}} +(s+\bar{t})a\bar{d}Y^{s-\bar{t}}\right) $$
 $$
-\frac{|\Lambda |}{|\Gamma _{\infty }:\Gamma '_{\infty }|}\left((s+\bar{t})b\bar{c}Y^{-s+\bar{t}}+ (s-\bar{t})b\bar{d}Y^{-s-\bar{t}}\right). $$
In the terms above involving multiplications of  lattice sums, after
integration over $P'$ the only surviving terms are (by the orthogonality of the family $\{ e^{2\pi i<\mu ,z>} \}$) \[
\frac{|\Lambda |}{|\Gamma _{\infty }:\Gamma '_{\infty }|}\sum _{\mu \neq 0}\left(f_{\mu }(Y)\bar{g}'_{\mu }(Y)-\bar{g}_{\mu }(Y)f'_{\mu }(Y)\right)\frac{1}{Y}.\]
 A direct calculation using the fact that $g_\mu $ and $f_\mu$ satisfy the Bessel equation  \eqref{E:bessel} shows that $$
\frac{d}{dr}(\frac{\left(f_{\mu }(r)\bar{g}'_{\mu }(r)-\bar{g}_{\mu }(r)f'_{\mu }(r)\right)}{r})=\frac{1}{r^{3}}(s^{2}-\bar{t}^{2})f_{\mu }(r)\bar{g}_{\mu }(r) $$
 thus $$
\frac{|\Lambda |}{|\Gamma _{\infty }:\Gamma '_{\infty }|}\sum _{\mu \neq 0}\left(f_{\mu }(Y)\bar{g}'_{\mu }(Y)-\bar{g}_{\mu }(Y)f'_{\mu }(Y)\right)\frac{1}{Y} $$
 \[
=(s^{2}-\bar{t}^{2})\frac{|\Lambda |}{|\Gamma _{\infty }:\Gamma '_{\infty }|}\sum _{\mu \neq 0} \int _{Y}^{\infty }f_{\mu }(r)\bar{g}_{\mu }(r)\, \frac{dr}{r^{3}}=(s^{2}-\bar{t}^{2})  \int _{\F (Y)}f^{Y}\bar{g}^{Y}\, dv. \]
By definition  $ \int _{D}f^{Y}\bar{g}^{Y} \, dv  = \int _{D}f\bar{g} \, dv     $ and thus combing the integral of $D$ and $F(Y) $ we obtain the lemma.
 \end{proof}

\subsection{The vector form of the Maa\ss-Selberg relations}
Let $\Gamma $ be cofinite with cusps $\zeta_{\alpha} = B_{\alpha}^{-1} \infty $ for $\alpha = 1\dots \kappa,$ with fundamental domain  $\F =\F _{0}\cup \F _{1}(Y)\cup \dots \F _{\kappa }(Y).$  For $(\chi, V) \in \rp. $

Let  $f \in A(\Gamma, \chi ) $ be a Maa\ss~ form, that is $f$ is smooth,  satisfies $ f(B_\alpha^{-1}(z+rj)) = O(r^N)$  for some $N \in \NN $ and for each cusp $\zeta_\alpha $ of $\Gamma,  $ and  $ \lp f = \lambda f  $ with $ \lambda = 1-s^2,~s \neq 0. $  We define
\[
f^{Y}(P)=\left\{ \begin{array}{cc}
 f(P)- \sum_{l=1}^{k_\alpha} u_{f_{\alpha l}}(P) & \text{if} \, P\in \F _{\alpha}(Y), ~\alpha=1\dots \kappa \\
 f(P) & \text{otherwise}  \end{array}
\right. \]

\begin{thm} \label{thmVecMasSel}
(Maa\ss-Selberg Relations) Let $\Gamma $ be cofinite, $(\chi, V) \in \rp,$ $f$ and $g$ be Maa\ss ~forms satisfying  $\lp f=\lambda f$, $\lp g=\nu g$
with $\lambda =1-s^{2},\nu =1-t^{2},s\neq 0,t\neq 0.$ 
Suppose that 
 $u_{f,\alpha  l}(B_{\alpha}^{-1}P)=a_{\alpha  l}r^{1+s} +b_{\alpha l}r^{1-s}$
and $u_{g,\alpha l}(B_{\alpha}^{-1}P)= c_{\alpha l} r^{1+s}+d_{\alpha l}r^{1-s}.$ Then 

\begin{multline*}
(s^{2}-\bar{t}^{2})\int _{\F }\left<f^{Y},g^{Y}\right>_V\, dv \\=\sum _{\alpha=1}^{ \kappa} \frac{|\Lambda_\alpha |}{|\Gamma_{\alpha}:\Gamma_\alpha^\prime |} \sum _{ l=1}^{k_{\alpha}} 
\left[ ~(s-\bar{t})\left<a_{\alpha l},c_{\alpha l}\right>_VY^{s+ \bar{t}}+ (s+\bar{t})\left<a_{\alpha l},d_{\alpha l}\right>_VY^{s-\bar{t }} \right. \\-
  \left.
(s+\bar{t})\left<b_{\alpha l}, c_{\alpha l}\right>_VY^{-s+\bar{t}}+(s-\bar{t}) \left<b_{\alpha l },d_{\alpha l }\right>_VY^{-s-\bar{t}}~ \right] 
\end{multline*}
\end{thm}

\begin{proof}
For notational simplicity we assume $\Gamma $ has only one class of cusps at $\infty .$ We can write $P_{\infty }=P_{1}+\dots +P_{k_\infty }$ (see \S\ref{secForExp}).
\begin{multline}
(s^{2}-\bar{t}^{2})\int _{\F }\left<f^{Y},g^{Y}\right>_V\, dv \\ =(s^{2}-\bar{t}^{2})\int _{\F }\left<(P_{\infty }+(1-P_{\infty })  )f^{Y},(P_{\infty }+(1-P_{\infty }) )g^{Y}\right>_V\, dv    
\\=(s^{2}-\bar{t}^{2})\int _{\F }\left(\left<P_{\infty }f^{Y},P_{\infty }g^{Y}\right>_V+\left< (1-P_{\infty }) f^{Y},(1-P_{\infty })g^{Y}\right>_V\right)\, dv.  \label{eq:ms1}
\end{multline}
The equality follows from the fact that $P_\infty $ is an orthogonal projection.   We will first show 
\begin{equation}
(s^{2}-\bar{t}^{2})\int _{\F }\left<(1-P_{\infty })f^{Y},(1-P_{\infty })g^{Y}\right>_V\, dv=0\label{eq:ms2}
\end{equation}

Since $(1 - P_\infty )f $ has no \emph{constant} Fourier expansion coefficients (see Proposition \ref{P:vectorexpan})   $(1 - P_\infty )f^{Y}= (1 - P_\infty) f, $ and by \mbox{Lemma \ref{lemBesSumDec}}
\beq \label{eqMsDec}
\frac{d^{p}}{dr^{p}}(1-P_\infty) f(B_\alpha^{-1}(z+rj)) = O(e^{-|c|r}) \eeq
 for $p \geq 0$ and as $r \ra \infty. $

Next to prove equation \eqref{eq:ms2} we apply the 
 vector form of Greene's theorem, $$
\int _{D}(\left<f,\lp g\right>_V-\left<g,\lp f) \right> \, dv=\int _{\partial D}(\left<f,\frac{\partial }{\partial n}g\right>_V-\left<g,\frac{\partial }{\partial n}f\right>_V)\, \widetilde{dv} $$ and we obtain (by \eqref{eqMsDec})

\begin{multline*}
\lim_{R \ra \infty} (s^{2}-\bar{t}^{2})\int _{\F _{R}}\left< (1-P_\infty)  f, (1-P_\infty) g\right>_V\, dv \\ 
=\lim_{R \ra \infty}  \int _{P(R)} \left( \left< (1-P_\infty) f,\frac{\partial }{\partial r}  (1-P_\infty) g\right>_V \right. \\
\left. -\left< (1-P_\infty) g,\frac{\partial }{\partial r} (1-P_\infty) f\right>_V \right) \, \frac{dxdy}{R}  = 0. \end{multline*}

To conclude the proof write \[
(s^{2}-\bar{t}^{2})\int _{\F }\left<P_{\infty }f^{Y},P_{\infty }g^{Y}\right>_V\, dv=(s^{2}-\bar{t}^{2})\sum _{j=1}^{k}\int _{\F }\left<P_{j}f^{Y},P_{j}g^{Y}\right>_V\, dv\]
 and apply Lemma (\ref{lem:ms1}) to $(s^{2}-\bar{t}^{2})\int _{\F }\left<P_{j}f^{Y},P_{j}g^{Y}\right>_V\, dv.$
\end{proof}
\section{The resolvent kernel} \label{secResKerEst}
\subsection{The definition of the resolvent}
One of our main goals in this paper is to understand the spectrum of $\lp.$  The spectrum is best understood in terms of the \emph{resolvent} operator.  The resolvent set of $\lp,~$  $ \rho(\lp), $ is the set of all $z \in \CC $ so that $ \left( \lp -zI \right) $ has a bounded inverse, defined on the entire Hilbert space $\hs$.  The spectrum of $\lp, $ $\sigma(\lp), $ is the complement of the resolvent set $\CC \setminus \rho(\lp). $  For each $z \in  \rho(\lp)$ we have the resolvent operator $ R_z = \left( \lp -zI \right)^{-1}. $ 

\subsection{The resolvent kernel for $\lp(\Gamma, \chi) $}
We will understand the resolvent operator through its kernel.
A bounded linear operator $K:\hs \mapsto \hs$ is said to have a \emph{kernel}
if there exists a function $k:\HH \times \HH \times V \mapsto V$ so that for all
$f\in\hs\,,P\in\HH,~$ \mbox{$
Kf(P)=\int_{\F}k(P,Q)f(Q)\, dv(Q).$}  Recall (from \S\ref{secFreCas}) that for $t > 1, s \in \CC,$ 
$$ \py_{s}(t) \df \frac{\left(t+\sqrt{t^{2}-1}\right)^{-s}}{\sqrt{t^{2}-1}}, $$ 
and  for $P,Q \in \HH, ~ \delta(P,Q)  = \cosh(d(P,Q)), $ where $d$ is the hyperbolic distance function.
\bd
Let $\Gamma $ be a cofinite Kleinian group, $\chi \in \rp,$ and  $s \in \CC $ with $\R(s)> 1.$  Define the Maa\ss-Selberg Series 
\beq \label{S:rk}
\rv=\frac{1}{4\pi}\sum_{M\in\Gamma} \chi(M)\varphi_{s}(\delta(P,MQ)). 
\eeq
\ed

Since $ \sum_{M\in\Gamma}  \left| \varphi_{s}(\delta(P,MQ))  \right| $  converges uniformly when $P,Q, ~ \text{and}~ s$ are restricted to compact subsets of \newline $\{ \R(s)> 1 \} \times \left(\HH \times 
\HH \setminus  \, \{\,(P,Q)\in \HH \times \HH  \, | \, \Gamma P = 
\Gamma Q
 \} \right) $ (see \cite[page 96]{Elstrodt}) and $\chi $ is a unitary representation  ($\chi(\gamma) $ has norm one for all $\gamma \in \Gamma $) the  sum in equation \eqref{S:rk} converges absolutely and uniformly when $P,Q, ~ \text{and}~ s$ are restricted to compact subsets of \newline $\{ \R(s)> 1 \} \times \left(\HH \times 
\HH \setminus  \, \{\,(P,Q)\in \HH \times \HH  \, | \, \Gamma P = 
\Gamma Q
 \} \right). $

The function $\rv $ is the kernel for the resolvent operator of $\lp. $  
\begin{thm} 
Let $\Gamma $ be cofinite with $\chi \in \rp $, and  $\lambda = 1-s^2 \in \rho(\lp) $ with  $\R(s)> 1.$ Then 
 $R_\lambda:\hs \ra \widetilde{D}$ has kernel $\rv$.
\end{thm} 
To prove the theorem one shows that  for $\lambda = 1-s^2 $ with $\R(s)> 1, $ any element 
$u \in \widetilde{D}$ can be represented by a continuous function 
which satisfies  \beq  \label{vres}
u(Q)=\frac{1}{4\pi}\int_{\F}\rv(\lp-\lambda)u(P)\, dv(P). 
\eeq
Equation \eqref{vres} can be seen by unfolding the integral in \eqref{vres} and applying Lemma \ref{lemScalRes}.  See \cite[theorem 4.2.6, page 150]{Elstrodt} for more details.

\subsection{Fourier expansion of the resolvent kernel}
In this section we will concern ourselves with a sub-sum of the resolvent kernel that causes the resolvent kernel to not be of Hilbert-Schmidt type. For notational simplicity we assume that $\Gamma $ is cofinite with only one cusp at $\zeta_\alpha = \infty. $ Fix a fundamental domain $\F = \F_Y \cup \F(Y),$ and let $\chi \in \rp. $ 

Recall that for $\R(s) > 1 $ the resolvent kernel can be expressed as a sum 
$$ \rv=\frac{1}{4\pi}\sum_{M\in\Gamma} \chi(M)\varphi_{s}(\delta(P,MQ)).  $$  Let 
$$F_\infty (P,Q) \df \frac{1}{4\pi}\sum_{M \in \Gamma_\infty} \chi(M)\varphi_{s}(\delta(P,MQ)),   $$ 
and $$ F_0 (P,Q) =\rv - F_\infty (P,Q).    $$  

Let $ \{ v_{ \infty l} \} $ be the basis constructed in \S\ref{secForExp}.  For convenience we drop the $\infty$ subscript and consider $ \{ v_{l} \} $   and corresponding orthogonal projections $ \{P_l \}. $

The following lemmas are all proved using a standard separation of variables technique.  See \cite[Pages 20-23]{Venkov} and \S\ref{secScalExpan} 
\begin{lem}
Let $P,Q \in \HH$ satisfy $ P \neq Q. $  Then
$$ P_\infty F_\infty (P,Q) =  P_\infty R_0(r_P, r_Q, s) + 
 P_\infty  \sum_{ 0 \neq \mu \in \Lambda_\alpha^0} a_{\mu} R_\mu (r_P,r_Q,s) e^{2 \pi i \left<\mu,z_P - z_Q \right >}, 
$$ 
where 
$$ R_0 = \frac{c}{s}
\left\{ \begin{array}{cc}  r_P^{1+s}r_Q^{1-s}   &  r_P \leq r_Q  \\
  r_P^{1- s}r_Q^{1+s} &  r_P \geq r_Q \\
\end{array}  \right.
$$ for some constant $c,$ and 
$$ R_\mu (r_P,r_Q,s) = \left\{ \begin{array}{cc}  
  r_P K_s(2 \pi |\mu| r_P) ~ r_Q I_s(2 \pi |\mu| r_Q) &  r_P \geq r_Q  \\
 r_P I_s(2 \pi |\mu| r_P) ~r_Q K_s(2 \pi |\mu| r_Q)  &  r_P \leq r_Q \\
\end{array}  \right..   $$

\end{lem}

\begin{lem}
Let $P,Q \in \HH$ satisfy $ P \neq Q, k_\infty + 1 \leq l \leq l_\infty. $  Then
$$ P_l F_\infty (P,Q) = $$
 $$ P_l  \sum_{ 0 \neq \mu \in \Lambda_\alpha^0} a_{l, \mu} R_\mu (r_P,r_Q,s) e^{2 \pi i \left<\mu,z_P - z_Q \right >}, 
$$ 
where 
$$ R_\mu (r_P,r_Q,s) = \left\{ \begin{array}{cc}  
  r_P K_s(2 \pi |\mu| r_P) ~ r_Q I_s(2 \pi |\mu| r_Q) &  r_P \geq r_Q  \\
 r_P I_s(2 \pi |\mu| r_P) ~r_Q K_s(2 \pi |\mu| r_Q)  &  r_P \leq r_Q \\
\end{array}  \right..   $$
\end{lem}

\begin{lem}
Let $P,Q \in \HH$ satisfy $ P \neq Q, l_\infty + 1 \leq l  \leq n. $  Then
$$ P_l F_\infty (P,Q) = $$
 $$ P_l  \sum_{ 0 \neq \mu \in \Lambda_\alpha^0} a_{l, \mu} R_\mu' (r_P,r_Q,s) e^{2 \pi i \left<\mu',z_P - z_Q \right >}, 
$$ 
where 
$$ R_\mu' (r_P,r_Q,s) = \left\{ \begin{array}{cc}  
  r_P K_s(2 \pi |\mu'| r_P) ~ r_Q I_s(2 \pi |\mu'| r_Q) &  r_P \geq r_Q  \\
 r_P I_s(2 \pi |\mu'| r_P) ~r_Q K_s(2 \pi |\mu'| r_Q)  &  r_P \leq r_Q \\
\end{array}  \right.,   $$
and $\mu '=\theta _{R}\mu _{1}+\theta _{S} \mu_{2}+\mu.$  (See Proposition \ref{P:vectorexpan} for an explanation of the notation used).
\end{lem}

\subsection{Estimates of the resolvent kernel}
It follows from \cite[Theorem 4.5.2]{Elstrodt} that 
$$ \int_\F |F_0 (P,Q)|_V^2~dv(P)~dv(Q) < \infty.   $$
In other words, $F_0 (P,Q) $ is a Hilbert-Schmidt kernel.  On the other hand,  $F_\infty (P,Q) $ is not.  However, from the Fourier series expansion of  $F_\infty (P,Q), $ and the decay properties of the Bessel functions $r I_s(2 \pi |\mu| r), r K_s(2 \pi |\mu| r) $ it follows that the only component of the resolvent kernel that destroys square-integrability is, 
$$ 
R_\infty(P,Q) \df \frac{c}{s}P_\infty \left\{ \begin{array}{cc}  r_P^{1+s}r_Q^{1-s}   &  r_P \leq r_Q  \\
  r_P^{1- s}r_Q^{1+s} &  r_P \geq r_Q \\
\end{array}  \right.  $$

\bd  Let $\Gamma $ be cofinite with one cusp at $\zeta_\alpha = \infty. $ For $B \geq Y$ decompose  $\F = \F_0 \cup \F(B),  $ and let  $\chi \in \rp, $  define \beq   
F^B(P,Q) \df  \left\{ \begin{array}{cc}  R_\infty(P,Q)  &  P,Q \in \F(B)  \\
  0 &  \text{else} \\
\end{array}, \right. \eeq 
and $ F_B(P,Q) \df \rv - F^B(P,Q). $
\ed
We have shown the following:
\begin{lem}
$$\int_\F |F_B(P,Q)|_V^2 ~dv(P)~dv(Q) < \infty.    $$ 
\end{lem}

\section{Analytic continuation of the Eisenstein Series}  
The meromorphic continuation of $E_{\alpha l}(P,s) $ is necessary for the proof of the  spectral decomposition theorem, and is highly non-trivial.  Fortunately, there are several well-known methods available,  \cite{Faddeev}, \cite{Selberg2}, \cite{Selberg3}, \cite{Verdi}, and \cite{Langlands}.  In \cite{Elstrodt} an adaptation of the methods in \cite{Verdi} is used to prove the three-dimensional scalar case, and a similar adaptation works\footnote{Alternatively, Faddeev's method can also be used  to express the meromorphic continuation of $E_{\alpha l}(P,s) $ by adapting \cite[Chapters 2 and 3]{Venkov}.} for the vector case which we show now.

\subsection{The idea of the proof}
As in Faddeev's proof the key idea is the resolvent kernel of $\lp. $ We know from  \S\ref{secResKerEst} that it is not a Hilbert-Schmidt operator. However the key step in the proof is  to construct a self adjoint operator $\lp_a $ defined on a closed subspace (the set of functions whose constant Fourier term vanishes for $r > a$ ) so that  $ (\lp_a - \lambda)^{-1}$ is Hilbert-Schmidt and is closely related to $\lp.$   

\subsection{Statement of the theorem}

\begin{thm} \label{T:Selberg} Let $\Gamma $ be cofinite, $(\chi,V) \in \rp, $ and $ 1 \leq l \leq k_\alpha.$  Then 

(1) $E_{\alpha l}(P,s)$ admits a meromorphic continuation in the 
following sense. 
For each fixed $P \in \HH$,  $E_{\alpha l}(s,P)$ is a meromorphic function in $s \in \CC$. The poles of  $E_{\alpha l}(P,s)$ 
depend only on $s$ and not on $P$.  

(2) If $U \subset \CC$ is open with $E_{\alpha l}(P,s)$ regular on $U$
then $E_{\alpha l}$ is real analytic on $U \times \HH$.

(3) For each regular $s \in \CC$,  $\lp E_{\alpha l}(P,s) = (1-s^2)
E_{\alpha l}(P,s)$.

(4) The scattering matrix $\smat(s) $ and $E_{\alpha l}(P,s)$ are both finite ordered meromorphic functions with order $ \leq 4. $

(5) $\mathcal{E}(P, -s) = \smat(-s) \mathcal{E}(P, s) $

(6) $ \smat(s) \smat(-s) = 1 .$

(7) $ \smat(s) $ is a unitary matrix on the critical axis $\R(s) = 0.$
\end{thm}

\subsection{ Proof of theorem \eqref{T:Selberg} part I: Preliminaries}
In order to avoid complicated notation (quadruple subscripts) we will assume the following:
\begin{as} \label{asOne}
The Kleinian group $\Gamma$ has only one class of cusps at \mbox{$\zeta = \infty \in \PP,$} and $\chi \in \rep.$
\end{as}
Then we can choose  fundamental domain of the form $\F = \F_0 \cup \F(Y) $ (here $\F_0 = \F_Y,$ see \S\ref{secCofKleGro} ). We can write the orthogonal projection onto $V_\infty $ by    $$P_\infty = \sum_{k = 1}^{k_\infty}P_k, $$ where $P_k $ is the orthogonal projection onto $v_k = v_{\infty k} $  (see \S\ref{secForExp}).  

Since we assume that $\infty $ is the only cusp we have (by definition) $k_\infty $ distinct Eisenstein series ( $ E_{\infty l}(P,s),~l = 1\dots k_\infty    $ ) which have an expansion of the form (see \S\ref{S:fourvec}, \S\ref{secScaMat}). 
\beq 
P_k E_{\infty l}(P,s) = \delta_{k l}r^{1+s}v_k + \phi_{k l}(s) r^{1-s} v_k  + O(e^{-|c|r}). 
\eeq

We will need various spaces of functions.
Let $\pi_\Gamma : \HH \ra \Gamma \setminus \HH$ be the standard (quotient map) projection, and set  
$$ \D^\infty \df  \{ \Phi  \in  A(\Gamma, \chi) \cap C^\infty(\HH,V)  ~ | ~ \pi_\Gamma( \text{supp}(\Phi)) ~\text{is compact in }~ \Gamma\setminus\HH \}. 
$$ 
We will restrict functions in $\D^\infty $ to $\F.$  The space of distributions $\D' $ is the set of all continuous (complex valued) linear functionals on $\D^\infty. $  A function $f \in L_\text{loc}^1 (\Gamma, \chi)  $ can be identified with the distribution $U$ where 
\beq \label{E:373721}
U(\Phi) \df \int_{\F}  \left< \Phi(P) ,f(P) \right>_V~dv(P). 
\eeq  
Noting that $ \hs \subset L_\text{loc}^1 (\Gamma, \chi),$ we  say that a distribution $U \in \hs $ if and only if  there exists $f \in \hs $ so that equation \eqref{E:373721} holds for all $\Phi \in \D^\infty. $

For $Y$ chosen in the decomposition $\F = \F_0 \cup \F(Y)$ define  $$\D(Y,\infty) \df C_c^\infty ( (Y,\infty), \CC), $$ and $\D'(Y,\infty) $ the corresponding space of distributions (the space of continuous, complex valued, linear functionals of $\D(Y,\infty)$ ).  For each $h \in L_\text{loc}^1(Y,\infty) $ we associate the distribution $$ h(\phi) \df \int_Y^\infty \phi(r)  \overline{h(r)}~\frac{dr}{r^3}~(\phi \in \D(Y,\infty).  $$

 For $U \in \D', \Phi \in \D^\infty, i \in \NN \cup \{0 \}, $  we define the following differential operators at each point $(x,y,r) \in \HH: $
$$
 D_x^{[i]}U  (\Phi) \df (-1)^i U(r^i \frac{\partial^i}{\partial x^i} \Phi ), ~ D_y^{[i]}U (\Phi) \df (-1)^i U(r^i \frac{\partial^i}{\partial y^i} \Phi ),
$$
$$  D_r^{[i]}U (\Phi) \df (-1)^i U(r^{3}  \frac{\partial^i}{\partial x^i}(r^{i-3} \Phi) ),  $$
$$\lp U \df  \left( D_x^{[2]} +D_y^{[2]} + D_r^{[2]} - D_r^{[1]} \right) U. $$ 
An unraveling of the definitions shows that $ \lp U(\Phi) = U(\lp \Phi ). $
For \mbox{$f \in \hs $} we abuse notation and define $$ \textbf{grad}(f) \in \hs ~ \text{iff} ~ D_x^{[1]}f, D_y^{[1]}f,  D_r^{[1]}f \in \hs. $$ 

Let $f,g \in \hs $ with $\textbf{grad}(f), \textbf{grad}(g) \in \hs. $  Then we define $$   
Q(f,g) \df \int_\F \left(  \left< D_x^{[1]}f, D_x^{[1]}g \right>_V + \left< D_y^{[1]}f, D_y^{[1]}g \right>_V  + \left< D_r^{[1]}f, D_r^{[1]}g \right>_V     \right)~dv.
$$

An application of Greene's theorem shows that for $\Phi \in \D^\infty, $
$$  Q(f,\Phi) = \int_\F \left< \lp \Phi,f\right>_V ~dv = ( \lp f)(\Phi), $$ where  $ \lp f  $  is taken  in the distributional sense. 

We next define useful maps between the function spaces $\D^\infty, \D(Y,\infty) $ and the distribution spaces $\D',\D'(Y,\infty). $  In order for the notation to be simple, we will use uppercase Greek letters for functions in $\D^\infty ,$ uppercase Roman letters for distributions in  $\D',$ lower case Greek letters for functions in  $\D(Y,\infty), $ and lower case Roman letters for distributions in $\D'(Y,\infty). $  It is also useful to think of $\D^\infty, \D' $ as \emph{BIG} spaces and $\D(Y,\infty), \D'(Y,\infty)$ as \emph{small} spaces.

For $F \in \hs, 1 \leq l \leq k_\infty $ set 
$$ L_l[F](r) \df  \frac{1}{ |\Lambda| } \int_{ \Lambda \setminus \CC} \left< P_l F(x,y,r), v_l \right>_V ~dxdy.   $$  
By Fubini's theorem $ L_l[F]  \in L^1_\text{loc} ( (Y, \infty), r^{-3}dr ).$
For any function $$\phi : (Y, \infty) \ra \CC, $$ and $1 \leq l \leq k_\infty$  set  
$$ B_l [\phi](P) \df  
\left\{ \begin{array}{cc}  \phi(r)v_l   &  \text{if $P = z+rj \in \F_\infty(Y) $ } \\
0  & \text{ if  $P  \in \F_0$ } \\
\end{array}  \right.  . $$  We will identify  $L_l [\phi](P)$ with its automorphic extension in $A(\Gamma, \chi ).$ 

For each $l \in 1\dots k_\infty $ define $ B_{l *}:  D' \ra D'(Y, \infty )  $ by  
$$  B_{l *}[U](\phi) \df U ( B_l[\phi] ), $$ where $U \in \D', \phi \in \D(Y, \infty ),  $ and define $ L_{l}^*  \D'(Y, \infty ) \ra \D' $ by 
$$ L_{l}^*[u](\Phi) = u(L_l[\Phi]),   $$ where $u \in  \D'(Y, \infty ), \Phi \in \D^\infty. $

Define $\lp_0 \D'(Y,\infty) \ra \D'(Y,\infty) $ by $$ (\lp_0 u)(\phi) = u \left( \left( -r^2 \frac{d^2}{dr^2} + r \frac{d}{dr} \right) \phi \right) \quad (u \in \D'(Y,\infty)), \phi \in \D(Y,\infty).$$  Then it follows that for $u \in \D'(Y,\infty),U \in \D',  l \in 1\dots k_\infty, $
$$ L_{l}^*(\lp_0 u) = \lp L_{l}^*(u)  ~\text{and}~ B_{l *}(\lp U) = \lp_0 B_{l *}(U). $$

For $a > Y, s\in \CC $ set 
$$  \eta_{a,s} \df  \left\{ \begin{array}{cc}  r^{1+s} - a^{2s} r^{1-s}  &  r > a  \\
0  &  r < a \\
\end{array}  \right.,   $$ and 
$ T_l^a = L_l^* [\delta_a],  $ where $\delta_a \in \D'(Y,\infty) $ is the Dirac delta-distribution at the point $a.$

A simple calculation using Greene's function techniques for second order differential equations (see \cite[Pages 194-195]{Iwaniec} for details on the technique) shows that 
$$ ( \lp_0 - (1-s^2)) \eta_{a,s} = -2sa^{s-1} \delta_a,  $$  and an application of $ L_{l}^* $ ($l \in 1\dots k_\infty $) gives us 
\beq
 ( \lp -(1-s^2) ) B_l[  \eta_{a,s} ] = -2sa^{s-1} T_l^a.
\eeq

\subsection{ Proof of theorem \eqref{T:Selberg} part II: Lemmas}

\begin{lem} \label{lemTruLap}
Let  $$ \mathfrak{H}_a \df \{ F \in \hs ~|~ L_l[F](r) = 0, r \geq a, l \in 1\dots  k_\infty    \},   $$ 
$$\mathfrak{D}_a  \df \{ F \in \mathfrak{H}_a ~|~ \emph{\textbf{grad}}(F) \in \hs \}, $$ and for $F,G \in \mathfrak{D}_a $    
$$ Q_a (F,G) \df  Q(F,G).   $$

Then 

(1) $ \mathfrak{H}_a $ is a closed subspace of $\hs.$  

(2) There exists a self-adjoint operator $\lp_a$ with domain $\mathfrak{D}_{ \lp_a} \subset \mathfrak{D}_a \subset  \mathfrak{H}_a $ with
 $$\lp_a : \mathfrak{D}_{ \lp_a} \ra \mathfrak{H}_a $$ 
that satisfies 
$$  \int_\F \left< \lp_a F,G \right>_V ~dv   = Q_a(F,G)     $$ 
for all $F \in  \mathfrak{D}_{ \lp_a}, G \in  \mathfrak{D}_a,$

(3) $ \mathfrak{D}_{ \lp_a} $ consists of those $F \in \mathfrak{D}_a $ which satisfy 
\beq  \lp F = G + \sum_{l = 1}^{k_\infty} c_l T_l^a, \label{E:doma}  \eeq
for some $G \in \mathfrak{H}_a. $  If equation  \eqref{E:doma} is satisfied, then $\lp_a F = G .$
\end{lem}

The proof (of the above lemma) is similar to \cite[Pages 237-239]{Elstrodt}.

\begin{lem} \label{lemCutRes}

(1) For $\R(s) > 1, $  $( \lp_a - (1-s^2) )^{-1} $ is a Hilbert-Schmidt operator. 

(2) The self-adjoint operator $ \lp_a$ has a purely discrete spectrum.

(3) The resolvent $( \lp_a - (1-s^2) )^{-1} $ can be continued to an operator valued meromorphic function on all of $\CC$ of order $ \leq 4. $
\end{lem}
\pf
(1) and (2) follow from a straight forward adaptation \cite[Prop. 6.1.8, page 240-241]{Elstrodt}.  (2) follows from standard functional analysis while (3) is proved in \cite[Lemma 6.1.9, page 240]{Elstrodt}.
\epf

Let $h_s(r) = h(r) r^{1+s} $ where for $Y < Y' < Y'' $ $h:(Y, \infty) \ra \RR $ smooth and $h(y) = 1 $ for  $y \geq Y'' $ and $h(y) = 0 $ for $y < Y'. $

Let $H_l (P,s) = (( \lp -(1-s^2) )B_l [h_s](P). $ Then it is in $\hs $ and smooth. In fact it is also in $\mathfrak{H}_a. $

Let  $$\Omega \df \{ s \in \CC~|~ \R(s)>0, 1-s^2 \notin \sigma(\lp) \}.$$
\begin{lem} \label{lemAdHocEis}
Let $Y < Y' < Y''$ and $h:(Y,\infty) \ra \RR$ be a smooth function satisfying $h(r) = 1$ for $ r \geq Y''$ and  $h(r) = 0 $ for $r < Y'.$  For each $l \in 1 \dots k_\infty$ set 
$$H_l (P,s) \df  (( \lp -(1-s^2) )B_l [h_s](P), $$ and for $s \in \Omega $ set 
$$ E^l(P,s) \df B_l [h_s](P) - (( \lp -(1-s^2) )^{-1}H_l (P,s).   $$
 
Then 

(1) $B_l [h_s] \in A(\Gamma,\chi) $ is a smooth. 

(2) $H_l (z + jr ,s) = 0  $ for all $r > Y'' $ and hence is in $\hs.$ Moreover, if   $ Y'' < a$ then  $ H_l(\cdot,s) \in \mathfrak{H}_a. $  

(3) The function $H_l(P,s) $ is holomorphic for $s \in \CC.$

(4) The function $E^l(P,s) $ is holomorphic for $s \in \Omega. $

(5)  The function $E^l(P,s) $ is the unique function in $A(\Gamma,\chi ) $ that satisfies  
$$ ( \lp -(1-s^2))E^l(P,s)=0,~ E^l(P,s) - B_l [h_s](P) \in \hs.  $$

(6)  For all $P \in  \HH,~ s \in \Omega, $ $E^l(P,s) = E_{\infty l}(P,s).   $  
 \end{lem}
See \cite[Pages 234-235]{Elstrodt}.

\subsection{Proof of theorem \eqref{T:Selberg} part III: The meromorphic continuation}
We exhibit the meromorphic continuation of the Eisenstein series to all of $\CC.$

For $s \in \CC $ with $1-s^2 \notin \Omega_a \df \sigma(\lp_a) $ set  
$$A_l(P,s) \df  ( \lp_a - (1-s^2) )^{-1} H_l (P,s). $$
By   Lemma \ref{lemCutRes} and Lemma \ref{lemAdHocEis},  $A_l(P,s)$ varies  meromorphically.  
Next for each $s\in \Omega_a$ $A_l(P,s)  \in \mathfrak{D}_{\lp_a} $ (since it is in the range of the resolvent of $\lp_a $), and an appeal to Lemma \ref{lemCutRes} give us
\beq  \label{E:tmp383821} ( \lp -(1-s^2) ) A_l(P,s) = H_l (P,s) + \sum_{j = 1}^{k_\infty} c_{l,j} (s) T_j^a. \eeq
Plugging in (to the distributional equation above) arbitrary test functions $\Phi \in \D^\infty, $ which do not vanish at the support of each $T_j^a,$ into  equation \eqref{E:tmp383821} shows that the functions $c_{l,j} (s) $ are meromorphic on all of $\CC $ (the remaining terms in equation \eqref{E:tmp383821} are  meromorphic).

Next define the following meromorphic functions, $$ e_{l,j} (s) \df  \frac{1}{2} s^{-1} a^{1-s} c_{l,j}(s). $$
Finally, let  
$$N_l(P,s) \df A_l(P,s) + \sum_{j = 1}^{k_\infty} e_{l,j} (s) B_l [\eta_{a,s}](P) + B_l[h_s](P).    $$
Then  (in the  distributional sense)  $( \lp -(1-s^2) ) N_l (P,s) = 0, $ and moreover, by the elliptic regularity theorem $N_l(P,s) $ is smooth as a function of $P\in \HH.$

Since $e_{l,j} (s),  B_l[h_s](P)$ are meromorphic, and of order $\leq 4,$    $N_l(P,s) $ is meromorphic, and of order $ \leq 4. $ 

To relate   $N_l(P,s) $ to $E_{\infty l}(P,s), $ express $E_{\infty l}(P,s) $ as a linear combination of the $N_l(P,s) $ and apply the uniqueness part of Lemma \ref{lemAdHocEis}. 

We have proved parts (1)-(4) of Theorem  \ref{T:Selberg}.  Part (5) can be seen by plugging in $-s$ into  $E_{\infty l}(P,s) $ and applying the uniqueness part of Lemma \ref{lemAdHocEis}. Part(6) follows immediately from (5). Part (7): from the Dirichlet series representation of the components of $\smat(s) $ we see that for $\R(s) > 1, $
$$ \smat(s) = \smat^*(\overline{s}). $$ Observing that for $s = it, t\in \RR,$ $\overline{it} = -it$ and appealing to (6) proves (7). 

\section{Proof of the spectral decomposition theorem}
There are several proofs of the spectral decomposition theorem.  All of them use the analytic continuation of the Eisenstein series to the critical line $\R(s) = 0 $ in a crucial manner.  

One particularly nice approach is to use the general theory of eigenpackets.  This approach makes  the deep connection between the Eisenstein series and the spectral decomposition theorem quite clear.  We shall see that integral of the  Eisenstein series on segments of the critical line ($\R(s)=0$) is an eigenpacket.

\subsection{Statement of the theorem}
Let $\Gamma $ be a cofinite group with $(\chi,V) \in \rp. $

Let $\mathfrak{D}$ denote a countable indexing set of all  $u \in D$ (the domain of $\lp$)  such that 
$\lp u = \lambda u$ for some $\lambda \in \CC $.  

We remark that in general, we do not know if $\mathfrak{D} $ is an infinite set.

\begin{thm}\label{T:spectral}
Let $f \in D$.  Then $f$  has an expansion of the form
\begin{multline}\label{E:spectral}
 f(P)  = \\ \sum_{m \in \mathfrak{D}} <f,e_m>e_m(P) + 
\frac{1}{4\pi}\sum_{\alpha = 1}^h \sum_{l = 1}^{k_\alpha}
\frac{\left[ \Gamma_\alpha:\Gamma_\alpha^\prime \right]}{\left| \Lambda_\alpha \right|}
\int_{\RR} <f,E_{\alpha l}(\,\cdot \,,it)>E_{\alpha l}(P,it) \, dt,
\end{multline}
The sum and integrals converge pointwise absolutely and uniformly on compact subsets of $\HH.$
\end{thm}

\subsection{Eigenpackets}
Throughout this section, let $A:\D_A \ra \bh $ be a self-adjoint operator with domain $\D_A $ in the Hilbert space $\bh,$ with inner product $ <\cdot,\cdot>.$  

We will outline the spectral decomposition theorem for eigenpackets.  See \cite[Section 6.2]{Elstrodt} for more details.

An eigenpacket of $A$ is a map $v:\RR \ra \bh, \lambda \mapsto v_\lambda  $ having the following properties:

(1) $v_0 = 0 $ and $v_\lambda \in \D_A $ for all $\lambda \in \RR.$

(2) The map $v:\RR \ra \bh$ is continuous (in the norm sense).

(3) $Av_\lambda = \int_0^\lambda \mu~ dv(\mu),$ where the integral is the limit in the norm sense of the corresponding Stieltjes sums.  By an integration by parts, we can rewrite (3) as $ Av_\lambda = \lambda v_\lambda - \int_0^\lambda v_\mu ~d\mu, $ where the integral is the limit of Riemann sums in the norm sense.

Eigenpackets have many similarities to eigenfunctions.  If $f$ is an eigenvector of $A$ then $ <f, v_\lambda> = 0 $ for all $\lambda \in \RR. $  In addition, if $w$ is also an eigenpacket of $A$ and $[\alpha, \beta], ~[\gamma,\delta]$ are two intervals having at most one point in common, then
\beq
<v_\beta - v_\alpha, w_\delta-w_\gamma> = 0.
\eeq

Fix an eigenpacket $v$ of $A.$  Then there exists a non-decreasing right continuous function $F:\RR \ra \RR$ which is unique up to an additive constant satisfying the equation 
\beq
F(\beta) -F(\alpha) = \|v_\beta - v_\alpha \|^2
\eeq
for all $\alpha < \beta.$    

For  $x \in \bh $ define a function $\varphi :\RR \ra \CC, $ by 
\beq
\varphi(\lambda) = <x,v_\lambda>.
\eeq 
Next let $T$ be a partion of $[\alpha,\beta],$ of the form
$$
T:\alpha = \lambda_0 < \lambda_1 < \dots < \lambda_n = \beta, 
$$
and define
$$
x_T \df  \sum_{j=1}^n \frac{\varphi(\lambda_j) - \varphi(\lambda_{j -1})}{F(\lambda_j) - F(\lambda_{j -1})}(v_{\lambda_j}-v_{\lambda_{j-1}}).
$$
For a sequence of partitions (sub partitions) $T_n$ with maximal width of its subdivsion converging to zero set
$$
  \int_\alpha^\beta \frac{d\varphi(\lambda) ~dv(\lambda)}{dF(\lambda)} = \lim_{n \ra \infty} x_{T_n}.
$$
Further define 
$$P_v x \df  \int_{-\infty}^\infty \frac{d\varphi(\lambda) ~dv(\lambda)}{dF(\lambda)} \df  \lim_{\alpha \ra -\infty} \lim_{\beta \ra \infty}  \int_\alpha^\beta \frac{d\varphi(\lambda) ~dv(\lambda)}{dF(\lambda)}.    $$

We can now write down the general spectral decomposition theorem. Let $\mathfrak{D} $ be a countable indexing set for a maximal set of orthonormal eigenvectors  $A,$ and $\mathfrak{N} $ a countable indexing set for a maximal set of orthogonal eigenpackets of $A$.  

\begin{lem} Let $ \{ e_m \}_{m \in \mathfrak{D}} $ be a maximal set of orthonormal eigenvectors of $A,$ and $ \{ v_n \}_{n \in \mathfrak{N}} $ be a maximal set of orthogonal eigenpackets of $A.$  Then every $x \in \bh $  can be decomposed as 
$$
x = \sum_{m \in \mathfrak{D}} <x,e_m> e_m + \sum_{n \in \mathfrak{N}}  \int_{-\infty}^\infty \frac{d\varphi_n(\lambda) ~dv_n (\lambda)}{dF_n(\lambda)}.
$$ 
The sum and integrals converge in the norm sense in $\bh.$
\end{lem}

\subsection{Proof of the theorem}
\begin{lem} Let $\Gamma $ be cofinite with $\chi \in \rp. $
For $P \in \HH,$ $\alpha \in 1 \dots \kappa,$ 
\newline $l \in 1\dots k_\alpha, $ we define 
$$
V_\lambda^{\alpha l}(P) \df \left\{ \begin{array}{cc} 0    & \text{for} ~\lambda < 1     \\
  \int_0^{\sqrt{\lambda -1}} E_{\alpha l}(P, it)~dt   & \text{for} ~\lambda \geq 1  \\
\end{array}  \right. 
$$
with positive choice of square root.  Then $V_\lambda^{[\alpha l]}(\cdot) \in \hs $ and the family of $V^{[\alpha l]}$ are a maximal system of orthogonal eigenpackets for $\lp.$  Further, 
\beq \label{E:ocon}
\| V_{\lambda_1}^{[\alpha l]} - V_{\lambda_2}^{[\alpha l]}   \|^2 = 2 \pi \frac{\left| \Lambda_\alpha \right|}{\left[ \Gamma_\alpha:\Gamma_\alpha^\prime \right]}(T_1 - T_2)
\eeq

for $\lambda_1 \geq \lambda_2 \geq 1$ and $T_1 = \sqrt{\lambda_1 - 1}, T_2 = \sqrt{\lambda_2 - 1}.$ 
\end{lem}
\pf
For notational simplicity we work under Assumption \ref{asOne}. From the Fourier expansion of the Eisenstein series we see that  $E_{\infty l}(P,s) = A_l(P,s) + O(e^{-|c|r}) $ as $r \ra \infty, $ where $A_l(P,s) =  r^{1+s}v_l + \phi_l(s)r^{1-s}  $, and $ \phi_l(s)$ is defined by $<\phi_l(s), v_k>_V =\smat(s)_{l, k} \df \smat(s)_{\infty l,\infty k}.$  Since $\smat(t) $ is unitary on the critical axis, $|\phi_l(it)|_V = 1 $ for $t \in \RR.$ 

We next show that  $V_\lambda^{[l]} \df V_\lambda^{[\alpha l]} \in \hs$ for each $\lambda \in \RR.$  Put $T = \sqrt{\lambda -1},$ and set  
$$ \beta_l (P,T) \df \int_0^T A_l(P,it) ~dt.$$  
An integration by parts shows that\footnote{The implicit constant depends on $T.$} $$\beta_l(P,T) = O(\frac{r}{\log r}) \quad \text{as $r \ra \infty.$}  $$ 
Hence $\beta_l(P,T) \in L^2(\F).$
The upper bound on $\beta(P,T) $ and the dominated convergence theorem imply that $V_\lambda^{[l]}$ is continuous in the norm sense.

Next, for $\lambda \geq 1,$
$$ \lp V_\lambda^{[l]}(P) =  \int_0^{\sqrt{\lambda -1}} (1+t^2)E_{\infty l}(P, it)~dt = \int_0^\lambda \mu ~dV_\mu^{[l]}. $$ The last equality was obtained by subsituting $\mu = 1+t^2.$ We have now shown that $V_\mu^{[l]} $ is an eigenpacket.

We next prove the orthogonality relation
$$\left< V_\lambda^{[l]}, V_\lambda^{[k]}\right> = 0  ~\text{for $ l \neq k$}, $$ using the Maa\ss-Selberg relations (Theorem \ref{thmVecMasSel}).  We have 
 $$\left< V_\lambda^{[l]}, V_\lambda^{[k]} \right>  = \lim_{Y \ra \infty} \left< \int_0^{T} E_{l}^Y (P, it)~dt   , \int_0^{T} E_{k}^Y(P, it')~dt' \right>  $$
$$ = \lim_{Y \ra \infty}\frac{[\Gamma_\infty: \Gamma_\infty^\prime]}{|\Lambda_\infty|}  \int_0^{T} \int_0^{T} \left< E_{l}^Y (P, it)   ,  E_{k}^Y(P, it') \right>~dt~dt'  $$ $$= \int_0^{T} \int_0^{T} \lim_{Y \ra \infty} \left[~ -Y^{it+it'}i \frac{ \overline{\smat_{k,l}(it')}}{t+t'}+ Y^{-it-it'}i \frac{ \overline{\smat_{k,l}(it)}}{t+t'}  \right. $$ $$ \left. -Y^{-it+it'}i \frac{ \sum_{j=1}^{k_\infty} \smat_{j,l}(it)\overline{\smat_{j,k}(it')} }{t'-t} ~ \right] ~dt~dt'  = 0,  $$  since all three terms above are bounded for $t,t' \in [0,T] $ ($\smat(it)$ is a unitary matrix so its terms are bounded and rows are orthonormal), and by the Riemann-Lebesgue lemma.

Equation \eqref{E:ocon} follows in a similar manner.
Letting $\lambda_1 = 1+T_1^2, \lambda_2 = 1+T_2^2$, we have
$$\| V_{\lambda_1}^{[l]} - V_{\lambda_2}^{[l]}   \|^2 = \lim_{Y \ra \infty}   \int_{T_2}^{T_1} \int_{T_2}^{T_1} \left< E_{l}^Y (P, it)   ,  E_{l}^Y(P, it') \right>~dt~dt'   $$
$$= \frac{\left| \Lambda_\alpha \right|}{\left[ \Gamma_\alpha:\Gamma_\alpha^\prime \right]} \lim_{Y \ra \infty} \int_{T_2}^{T_1} \int_{T_2}^{T_1} \left[ -Y^{it+it'}i \frac{ \overline{\smat_{l,l}(it')}}{t+t'}+ Y^{-it-it'}i \frac{ \overline{\smat_{l,l}(it)}}{t+t'}  \right. $$ $$  \left. +Y^{it-it'}i\frac{1}{t'-t}  - Y^{-it+it'}i \frac{ \sum_{j=1}^{k_\infty} \smat_{j,l}(it)\overline{\smat_{j,l}(it')} }{t'-t}  \right] ~dt~dt' = 
 2 \pi \frac{\left| \Lambda_\alpha \right|}{\left[ \Gamma_\alpha:\Gamma_\alpha^\prime \right]}(T_1 - T_2). $$   See \cite[Pages 272-275]{Elstrodt} and \cite[Page 97]{Iwaniec} for the details. The trick is to add and subtract the conjugate to the second to last term (on the left side of the equal sign) and combine one of the extra terms with the last term. 

Finally, it remains to prove completeness.  Let $A > Y,$ and extend the function    
$$  f_A^{[l]}(P) \df  \left\{ \begin{array}{cc} 
v_l    & \text{for} ~P \in \F(A)     \\
0   & \text{else}   \\
\end{array}  \right. 
$$
to $\hs. $ Then one can verify that $f_A^{[l]} $ is orthogonal to the set of \emph{cusp forms}\footnote{The space of cusp forms is the space spanned by the eigenvalues of $\lp$ that are in $\hs,$ and whose constant term in their Fourier expansion is identically zero. Only finitely many eigenvalues of $\lp$ are not cusp forms. }, and to all eigenpackets $V^{[k]} $ for $k\neq l.$  The function  $f_A^{[l]} $ is not orthogonal to the exceptional eigenfunctions (the residues of the Eisenstein   series $E_l(P,s)$  at the poles $\sigma_j \in (0,1]$ ) and the eigenpacket $V^{[l]}.$ Let  $T_l$ denote the orthogonal projection onto the space spanned by the exceptional eigenfunctions and the the span of $\{ V_\alpha^{[l]} - V_\beta^{[l]} \}_{\alpha,\beta \in \RR}. $  Then a computation shows that 
$\|  f_A^{[l]}  \|^2 = \| T_l f_A^{[l]} \|^2, $ that is our system eigenfunctions-eigenpacket system is complete with respect to expanding $f_A^{[l]}. $  Now suppose that there exists an eigenpacket $W$ orthogonal to all of the $ \{ V^{[l]} \}_{l =1\dots k_\infty}, $ and cusp forms, and let $U$ be the span of $\{ W_\alpha - W_\beta \}_{\alpha,\beta \in \RR}. $  Then for $u \in U$ by the completeness with respect to $f_A^{[l]}, $ $ \left< f_A^{[l]},u \right> = 0$ for $l =1\dots k_\infty,$ that is $u$ is a cusp form, a contradiction unless $u = 0.$  See \cite[Pages 274-276]{Elstrodt} for more details.
\epf

\chapter{The Selberg Trace Formula}
In this chapter we derive the Selberg Trace Formula for cofinite Kleinian groups with finite dimensional unitary representations.
\section{The Selberg Trace Formula}
\label{sectionSelberg}
\begin{thm}{(Selberg Trace Formula)} \label{T:SelbergTrace}
Let $\Gamma $ be a cofinite Kleinian group,  $\chi \in \rep,~h $ be  a holomorphic function on  $ \{ s \in \CC \, | \, |\I(s)| < 2+ \delta \}$ for some $\delta > 0,$ satisfying $ h(1+z^2) = O( 1+|z|^2)^{3/2 - \epsilon}) $ as  $|z| \ra \infty,$ and let
$$ g(x) = \frac{1}{2\pi} \int_{\RR} h(1+t^2)e^{-itx}\,dt. $$  Then
\begin{multline}
\sum _{m \in \D }h(\lambda _{m}) 
= \frac{\vol \left( \Gamma \setminus \HH \right)}{4\pi ^{2}}\dim_\CC V \int _{\RR }h(1+t^{2})t^{2}\, dt  \\ + 
 \sum_{ \{R \} \text{\emph{nce}}}\frac{\tr _{V}\chi (R) g(0)\log N(T_{0})}{4|\ren |\sin ^{2}(\frac{\pi k}{m(R)})}  + 
\sum_{\{ T \} \text{\emph{lox}} } \frac{\tr _{V}\chi (T) g(\log N(T))}{\oen |a(T)-a(T)^{-1}|^{2}}\log N(T_{0}) \\ - 
\frac{\tr (\smat (0))h(1)}{4}+\frac{1}{4\pi }\int _{\RR }h(1+t^{2})\frac{\phi'}{\phi }(it)\, dt  \\ +
\sum_{\alpha=1}^{\kappa}\sum _{k=1}^{e_{\alpha}}\frac{\tr _{V}\chi (g_{\alpha k})}{|C(g_{\alpha k})|}\left(g(0)c_{\alpha k}+d_{\alpha k}\int _{0}^{\infty }g(x)\frac{\sinh x}{\cosh x-1+\alpha _{\alpha k}}\, dx\right)  \\ + 
\sum _{\alpha =1}^{\kappa}\left(\frac{l_{\alpha}}{|\Gamma _{\alpha}:\Gamma '_{\alpha}|}\left(\frac{h(1)}{4}+g(0) (\frac{1}{2}\eta_{\Lambda_{\alpha}}-\gamma )- \frac{1}{2\pi }\int _{\RR }h(1+t^{2})\frac{\Gamma '}{\Gamma }(1+it)\, dt\right)\right)  \\ +
\sum _{\alpha=1}^{\kappa}\frac{g(0)}{|\Gamma _{\alpha}:\Gamma '_{\alpha}|}\sum _{k=l_{\alpha}+1}^{\dim_\CC V}L(\Lambda_{\alpha},\psi_{k \alpha}).
\end{multline}
\end{thm}
Here, $ \{ \lambda _{m} \}_ {m \in \D } $ are the eigenvalues of $\lp$ counted with multiplicity.  Following \cite{Elstrodt} section 5.2,  the summation with respect to $\{R \}_\text{nce} $ extends over the finitely many $\Gamma-$conjugacy classes of the non cuspidal elliptic elements (elliptic elements that do not fix a cusp) $R \in \Gamma,$ and for such a class $N(T_0)$ is the minimal norm of a hyperbolic or loxodromic element of the centralizer $\mathcal{C}(R).$  The element $R$ is understood to be a $k-$th power of a primitive non cuspidal elliptic element $R_0 \in \mathcal{C}(R)$ describing a hyperbolic rotation around the fixed axis of $R$ with minimal rotation angle $\frac{2 \pi}{m(R)}.$  Further, $\ren$ is the maximal finite subgroup contained in $  \mathcal{C}(R).$  The summation with respect to $\{ T \}_\text{lox}$ extends over the  $\Gamma-$conjugacy classes of hyperbolic or loxodromic elements of $\Gamma,$ $T_0$ denotes a primitive hyperbolic or loxodromic element associated with $T.$  The element $T$ is conjugate in $\pc$ to the transformation  described by the diagonal matrix with diagonal entries $a(T), a(T)^{-1}$ with $|a(T)| > 1, $ and $N(T) =   |a(T)|^2. $ For $s \in \CC,$ $\smat(s)$ is  a $k(\Gamma,\chi) \times  k(\Gamma,\chi)$ matrix-valued meromorphic function, called the \emph{scattering matrix}  of $\lp,$ and $\phi(s) = \det \smat(s).$  Furthermore $c_{\alpha k}, g_{\alpha k},$ and $d_{\alpha k}$ are constants depending on $\Gamma$ which will be determined in the case of $\Gamma $ having only one cusp at $\infty. $  The remaining notation will be defined in this chapter\footnote{Please note that there is a typographical error in the loxodromic and non cuspidal elliptic terms in \cite{Elstrodt} Theorem 6.5.1; both terms are missing a factor of $\frac{1}{4 \pi}. $ \\ }.

For $P = z+rj,~P'=z'+r'j \in \HH$ set 
$$ \delta(P,P') \df  \frac{|z-z'|^{2}+r^{2}+ r'^{2}}{2rr'}. $$ It follows that $\delta(P,P') = \cosh(d(P,P'))$, where $d$ denotes the hyperbolic distance  in $\HH.$  Next, for $k \in \scz([1,\infty))$ a Schwartz-class function, define 
$$K(P,Q) = k(\delta(P,Q)),~\text{and}~ K_\Gamma(P,Q) \df \sum_{\gamma\in\Gamma} \chi(\gamma)K(P,\gamma Q). $$ 

The series above converges absolutely and uniformly on compact subsets of
$\HH \times \HH$, and is the kernel of a bounded operator $\mathcal{K} : \hs \mapsto \hs.$  The Selberg trace formula is essentially\footnote{We say ``essentially'' because $\K$ is not of trace class.  Selberg's  procedure is used to define and compute the \emph{regularized trace.}\\} the trace of $\mathcal{K}$ evaluated in two different ways: the  first using spectral theory, and the second as an explicit integral.

The function $h$ that appears in the Selberg trace formula is the  Selberg--Harish-Chandra transform\footnote{If $f:\HH \mapsto V$ is a smooth function satisfying $\lp f = \lambda f$, then $\mathcal{K}f= h(\lambda)f.$  That is $f$ is an eigenfunction of $\mathcal{K}$ with an eigenvalue that depends only on $\lambda.$ \\} of $k,$ defined as follows:  
\beq \label{eqSHC}
h(\lambda) = h(1-s^2) \df \frac{\pi}{s}
\int_{1}^{\infty}k\left(\frac{1}{2}\left(t+\frac{1}{t}\right)\right)
(t^{s}-t^{-s})\left(t-\frac{1}{t}\right)\,\frac{dt}{t},~~\lambda = 1-s^2.
\eeq

\subsection{Expansion of $K_\Gamma$}
We now \emph{diagonalize} $\mathcal{K}.$
For $v,w \in V $ define a map $v \otimes \overline{w} : V \ra V$ by \beq 
v \otimes \overline{w}(x) = <x,w>v.
 \eeq
An immediate application of the spectral decomposition theorem and the Selberg transform give us, 
\begin{lem} \label{lemKerExp}
Let  $k \in \scz $ and  $h:\CC\ra\CC$ be the Selberg Transform
of $k.$ Then
\begin{multline}\label{E:kernel expansion}
K_\Gamma(P,Q) =   \sum_{m \in \D} h(\lambda_m)e_m(P) \otimes 
\overline{e_m(Q)} \\ +
\frac{1}{4\pi}\sum_{\alpha = 1}^h \sum_{l = 1}^{k_\alpha}
\frac{\left[ \Gamma_\alpha:\Gamma_\alpha^\prime \right]}
{\left| \Lambda_\alpha \right|}
\int_{\RR} h \left( 1+t^2 \right) E_{\alpha l}(P,it) \otimes \overline{E_{\alpha l}(Q,it)} \, dt.
\end{multline}
The sum and integrals converge absolutely and uniformly on compact subsets 
of $\HH \times \HH$.
\end{lem}

Set 
$$
H_\Gamma(P,Q) =  \frac{1}{4\pi}\sum_{\alpha = 1}^h 
\sum_{l = 1}^{k_\alpha}
\frac{\left[ \Gamma_\alpha:\Gamma_\alpha^\prime \right]}
{\left| \Lambda_\alpha \right|}
\int_{\RR} h \left( 1+t^2 \right) E_{\alpha l}(P,it) 
\otimes \overline{E_{\alpha l}(Q,it)} \, dt,
$$
$$
L_\Gamma(P,Q) =  \sum_{m \in \D} h(\lambda_m)e_m(P) \otimes 
\overline{e_m(Q)}.
$$  
Then clearly $ K_\Gamma = L_\Gamma + H_\Gamma. $
We have 
\begin{lem}
$$ \int_\F \|L_\Gamma(P,Q)\|_{V}^{2} \, dv(P) \,dv(Q) < \infty, $$
$$ \int_\F \tr_V (L_{\Gamma}(P,P))\,dv(P) = 
\sum_{m \in \D} h(\lambda_m) ~ \text{and},$$
$$ \sum_{m \in \D} |h(\lambda_m)| < \infty.
$$
\end{lem}
\begin{proof}
The first equation follows from the decay properties of 
$h$, the orthonormality of the $e_m$ and that fact that $$ \int_\F \|L_\Gamma(P,Q)\|_{V}^{2} \, dv(P) \,dv(Q) \leq 
C \sum_{m \in D}|h(\lambda_m)|^2 < \infty. $$
For the second, note that
$$ \tr_V(e_m(P) \otimes 
\overline{e_m(P)}) = \< e_m(P),e_m(P)\>_V $$
 thus by the orthonormality of the $e_m$, 
$$ \int_\F \tr_V (L_{\Gamma}(P,P))\,dv(P) = \sum_{m \in \D} 
h(\lambda_m) \int_\F \< e_m(P),e_m(P)\>_V \, dv(P) =
\sum_{m \in \D} h(\lambda_m).  $$

By the decay properties of $h$, 
\beq \sum_{m \in \D} |h(\lambda_m)| < \infty \eeq
\end{proof}
If we naively try to take the trace of $\K$ by $
 \int_\F \tr_V (K_{\Gamma}(P,P))  \,dv(P) 
$  we would see that the integral does not converge.  However we can \emph{regularize} the trace by subtracting off the term
$\int_\F \tr_V (H_{\Gamma}(P,P))  \,dv(P)$ in the following careful
manner. 

\begin{lem} \label{lemSinCanOut}
Let   $ f(P) = \tr_V (K_{\Gamma}(P,P)) - \tr_V (H_{\Gamma}(P,P)).$ 
Then $f$ is in $L^1(\F).$
Further,
\begin{multline} 
\lim_{A \ra \infty} \int_{\F_A} \tr_V (K_{\Gamma}(P,P))  \,dv(P)
- \lim_{A \ra \infty} \int_{\F_A} \tr_V (H_{\Gamma}(P,P))  \,dv(P)
\\=\int_\F f(P) \,dv(P).
\end{multline}
(The set $\F_A$ is defined in Lemma~\ref{lemFunDom}.) In particular, the singularities arising from each integral must
cancel out.
\end{lem}
The fact that the singularities arising from each integral must
cancel out will be crucial in our derivation of the Selberg trace formula, and imply some interesting group relations of $\Gamma. $
\subsection{Proof of the Spectral Truncated Trace}
For this section Assumption \ref{asOne} is in effect. So far we know little about the $t-$dependence of the Eisenstein series on the critical line.  The following consequence of Lemma~\ref{lemKerExp} will remedy the situation.

\begin{lem} \label{lemEisUpBou} 
Let $l \in 1\dots k_\infty.$  Then
$$ \int_{-T}^T |E_l(z+jr,it)|_V^2~dt = O(r^2 T + T^3).   $$
\end{lem}
\pf
For $0 < \epsilon < 1/16 $ let $k_\epsilon:[1,\infty) \mapsto \RR $ be a smooth function satisfying  $$
\text{supp}(k_\epsilon) \subset [1,1+\epsilon ],~\text{max} (|k_\epsilon|) \leq c \epsilon^{-\frac{3}{2}},~ \int_{\HH} k_\epsilon(\delta(P,Q))~dv(Q) = 1,
$$
for some $c > 0,$ and for all $P\in \HH$ (See \cite[page 292]{Elstrodt} for the construction).  Then it follows (\cite[page 292]{Elstrodt}) that  $$
\int_\F \| K_\Gamma (P,Q) \|_V^2~dv(Q) \leq c_1 \left( r_P^2 \epsilon^{-\frac{1}{2}} + \epsilon^{-\frac{3}{2}} \right)
$$
for all $P = (z_P,r_P) \in \F.$  On the other hand we have (Lemma~\ref{lemKerExp})
\begin{multline*}
K_\Gamma(P,Q) = \\ \sum_{m \in \D} h(\lambda_m)e_m(P) \otimes 
\overline{e_m(Q)}+
\frac{1}{4\pi} \sum_{l = 1}^{k_\infty}
C_\infty \int_{\RR} h \left( 1+t^2 \right) E_{l}(P,it) \otimes \overline{E_{l}(Q,it)} \, dt,
\end{multline*}
and by Parseval's equality (for the Hilbert space $\hs$), setting $\epsilon = (16T^2)^{-1},$ we obtain $$  
c_1 (r_P^2 T + T^3) \geq \| K_\Gamma (P, \cdot) \|_{\hs}^2 \geq \frac{C_\infty}{4\pi}\int_{-T}^{T} |h(1+t^2)|^2 |E_l(P,it)|^2~dt. 
$$
By \cite[page 122]{Elstrodt} we get $$|h(1+t^2) -1| < \frac{7}{2} (1+|t|) \frac{1}{4T},  ~ \text{for}~ |t| \leq T, $$ which gives a positive lower bound on $|h(1+t^2)|^2 $ for $ t \in [-T,T].  $ 
The  lemma now follows. 
\epf
The following lemma will allow us to compute the spectral trace.  More specifically the Maa\ss-Selberg relations will give the leading terms of the truncated trace while the error term (below) will vanish as $A \ra \infty. $   We adapt a nice argument found in \cite[Page 142]{Iwaniec} which shortens the exposition found in \cite[Page 295]{Elstrodt}.
\begin{lem} \label{lemEisErrDie}
Let $h(1+z^2) = O((1+|z|^2)^{-3/2 - \epsilon}).$ Then  
$$ \int_\RR \int_{\F(A)} h(1+t^2)|E_l^A(P,it)|_V^2~dv(P)~dt =O(A^{-1})~\text{as}~A \ra \infty.  $$
\end{lem}
\pf
Let $C_\infty = \frac{[\gi:\gip]}{|\Lambda_\infty|}$ and let $\mathcal{P}$ be a fundamental domain for the lattice $\Lambda_\infty.$
By Parseval's equality\footnote{$\phi_{l,\mu}(s)$ is the Fourier coefficient of $E_l(P,s).$ \\} (for eigenfunctions of the euclidean Laplacian on $\CC=\RR^2),$) 
\beq \label{eqEucPar} 
\frac{1}{|\Lambda_\infty|} \int_{\mathcal{P}} |E_l^A(z+rj,it)|_V^2~dz =  \sum_{\mu \neq 0} |\phi_{l,\mu}(it)|_V^2 |rK_{it}(2 \pi |\mu| r)|^2. 
\eeq  
We have the well known formula 
$$ K_{it}(r) = \pi^{1/2}\Gamma(it + \frac{1}{2})^{-1} \left( \frac{r}{2}\right)^{it} \int_1^\infty (x^2 -1)^{it - 1/2} e^{-xr}~dx, $$ which implies  (by an elementary argument) that for any $ 0< n < m,  $ constant  $A,$ sufficiently large, and implied constant \emph{independent} of $t,$
\beq \label{eqBesBou}
\int_A^\infty |K_{it}(r)|^2 r^{-n}~dr << \int_{A/2}^\infty |K_{it}(r)|^2 r^{-m}~dr.
\eeq
An application of  Parseval's equality, and \eqref{eqEucPar} gives   

\begin{multline*}
\int_{\F(A)} |E_l^A(P,it)|_V^2~dv(P)  \\ =\frac{1}{[\gi:\gip]} \int_A^\infty \int_{\mathcal{P}} |E_l^A(z+rj,it)|_V^2 dz \frac{dr}{r^3}  \\ = \frac{1}{C_\infty} \int_A^\infty \sum_{\mu \neq 0} |\phi_{l,\mu}(it)|_V^2 |rK_{it}(2 \pi |\mu| r)|^2~\frac{dr}{r^3} \\=C_\infty \int_A^\infty \sum_{\mu \neq 0} |\phi_{l,\mu}(it)|_V^2 |K_{it}(2 \pi |\mu| r)|^2~\frac{dr}{r}  \\
<<\int_{A/2}^\infty  \sum_{\mu \neq 0} |\phi_{l,\mu}(it)|_V^2 |K_{it}(2 \pi |\mu| r)|^2~\frac{dr}{r^5} 
\\=\int_{A/2}^\infty  \sum_{\mu \neq 0} |\phi_{l,\mu}(it)|_V^2 |rK_{it}(2 \pi |\mu| r)|^2~\frac{dr}{r^7} \\= \int_{A/2}^\infty  \frac{1}{|\Lambda_\infty|} \int_{\mathcal{P}} |E_l^{A/2}(z+rj,it)|_V^2~dz~\frac{dr}{r^7} 
 \\ \leq \int_{A/2}^\infty  \frac{1}{|\Lambda_\infty|} \int_{\mathcal{P}} |E_l(z+rj,it)|_V^2~dz~\frac{dr}{r^7}.
\end{multline*}
Hence, we infer from Lemma~\ref{lemEisUpBou}, and  another application of Parsavel's equality (for eigenfunctions of the Euclidean Laplacian)  that 
\begin{multline*}
\int_{-T}^{T} \int_{\F(A)} |E_l^A(P,it)|_V^2~dv(P) ~ dt \\ << \int_{-T}^{T} \int_{A/2}^\infty  \int_{\mathcal{P}} |E_l^{A/2}(z+rj,it)|_V^2~dz~\frac{dr}{r^7}
 \\ \leq  C \int_{A/2}^\infty \int_{-T}^{T} |E_l (P,it)|_V^2 dt \frac{dr}{r^7} \\ << \int_{A/2}^\infty  c_1(r^2 T + T^3) ~\frac{dr}{r^7} << \frac{T^3}{A^4}.
\end{multline*}
The lemma follows from the bound on $h(1+z^2)$ using an elementary integration by parts argument.
\epf

Define  $\phi(s) =  \det \smat(s). $
By Theorem \ref{T:Selberg} $\phi(s)$ is a meromorphic, finite ordered function on $\CC$,  $|\phi(s)| = 1$ on  $ \{ \R(s) = 0 \},$ and 
$\lds(s)$ is regular on $ \{ \R(s) = 0 \}.$ Let $h(1+z^2) = O((1+|z|^2)^{-3/2 - \epsilon}), $ and $$g(x) = \frac{1}{2 \pi} \int_\RR h(1+t^2) e^{-itx}~dt.$$ Then we have 
\begin{lem} \label{lemSpeTra}
\begin{multline}
\int_{\F_A} \tr_V (H_{\Gamma}(P,P))  \,dv(P) = \\ 
g(0)k(\Gamma,\chi) \ln(A) - \frac{1}{4\pi} \int_{\RR}\lds(it)
h(1+t^2)\,dt + \frac{h(1)\tr \smat(0)}{4} +  
\lto_{A \ra \infty}.
\end{multline}
The integral converges absolutely.
\end{lem}
\pf
\begin{multline*}
4 \pi  \int_{\F_A} \tr_V (H_{\Gamma}(P,P))  \,dv(P)  \\ = \sum_{l=1}^{k_\infty} \int_{\F_A} \tr C_\infty \int_{\RR} h \left( 1+t^2 \right) E_{l}(P,it) \otimes \overline{E_{l}(P,it)} \, dt  ~dv(P) \\
=\int_{\RR} h \left( 1+t^2 \right)  \sum_{l=1}^{k_\infty} \int_{\F_A} \tr C_\infty |E_{l}(P,it)|_V^2 ~dv(P) ~dt
\\
=\int_{\RR} h \left( 1+t^2 \right)  \sum_{l=1}^{k_\infty} \int_{\F_A} \tr C_\infty |E_{l}^A (P,it)|_V^2 ~dv(P) ~dt
\\
=\int_{\RR} h \left( 1+t^2 \right)  \sum_{l=1}^{k_\infty} \int_{\F} \tr C_\infty |E_{l}^A (P,it)|_V^2 ~dv(P) ~dt
\\ 
- \int_{\RR} h \left( 1+t^2 \right)  \sum_{l=1}^{k_\infty} \int_{\F(A)} \tr C_\infty |E_{l}^A (P,it)|_V^2 ~dv(P) ~dt.
\end{multline*}
By Lemma~\ref{lemEisErrDie} the term on the last line tends to zero as $A \ra \infty.$  We evaluate the first term above with the Maa\ss-Selberg relations on the critical line $\R(s)=0.$  For $t\in \RR, t' = t+r, r > 0, C_\infty = \frac{[\Gamma_\infty: \Gamma_\infty^\prime]}{|\Lambda_\infty|}$ and $\mathcal{E}^Y (P, it) $ the column vector of Eisenstein series ($l=1\dots k_\infty$).  Consider the  $k_\infty \times k_\infty $ matrix\footnote{Each entry is the inner product of two truncated Eisenstein series.  The inner product takes place in the Hilbert space $\hs.$\\},   $  \left< C_\infty \mathcal{E}^Y (P, it), \mathcal{E}^Y(P, it')^t \right>.$  Then  
$$  
\tr  \left< C_\infty \mathcal{E}^Y (P, it)   ,  \mathcal{E}^Y(P, it')^t \right> = \sum_{l=1}^{k_\infty} \int_{\F} \tr C_\infty |E_{l}^A (P,it)|_V^2 ~dv(P).
$$ 
Applying the Maa\ss-Selberg relations  $k_\infty \times k_\infty $ many times and placing the result in matrix form gives us\footnote{We are only interested in the trace of the matrix $ \left< C_\infty \mathcal{E}^Y (P, it)   ,  \mathcal{E}^Y(P, it')^t \right>.$ However by working with the entire matrix we can use the identity \eqref{eqTraLog} and derive a nice expression for the trace. \\}
$$ \left< C_\infty \mathcal{E}^Y (P, it)   ,  \mathcal{E}^Y(P, it')^t \right> =  $$
$$  -Y^{it+it'}i \frac{ \smat(it')}{t+t'}+ Y^{-it-it'}i \frac{ \smat(it)^*}{t+t'} +    Y^{it-it'}i\frac{I}{t'-t}  - Y^{-it+it'}i \frac{  \smat(it)\smat^*(it') }{t'-t}.  $$
Next applying the approximations, 
$$Y^{ir} = 1+ir \ln{Y}+ \dots, $$ 
$$ Y^{-ir} = 1-ir \ln{Y}+ \dots,  $$
$$ \smat( it + ir) =\smat(it) + ir \smat( it )^\prime + \dots,  $$
identities, 
$$\smat(-it) = \smat^* (it) = \smat^{-1}(it), $$  
$$\frac{d}{dt} \smat(it) \smat(-it) = -i \smat(it) \smat^\prime (-it) + i \smat^\prime (it)\smat(-it) = \frac{d}{dt} I  = 0,   $$
\beq  \label{eqTraLog} 
\tr{\smat^\prime (it) \smat^{-1} (it)} = \frac{d}{dt} \log \det{\smat(it)} = \frac{\phi^\prime}{\phi}(it),  
\eeq
and letting $r \ra 0$ we obtain the lemma.  See  \cite[Page 305]{Elstrodt} and \cite[Page 139-142]{Iwaniec} for more details.
\epf

\section{Lattice Characters and Sums} \label{secDoubleSum}
As mentioned earlier the parabolic elements of a cofinite Kleinian group have an associated lattice. 
\subsection{$ Z(x, \Lambda , \psi)$}
Let $\Lambda =  \ZZ \oplus \ZZ \tau \subset \CC$ be a lattice with $\I(\tau) > 0. $  A (lattice) character  $ \psi $ of $\Lambda $ is a one-dimensional unitary representation of $\Lambda. $
\bd
For  $x > 0 $ set 
$$ Z(x,\Lambda , \psi) \df \sum_{ \substack{ \mu \in \Lambda \\  | \mu |^2 \leq x \\ \mu \neq 0 } } \frac{ \psi(\mu)}{ |\mu |^2},    $$   
and when the limit exists 
$$  L(\Lambda , \psi ) \df  \lim_{x\ra \infty} Z(x,\Lambda , \psi). $$
\ed

\bp \label{P:latticesum} Let $ \psi $ be a  character of $\Lambda $. 

(1) If $\psi = \text{id}$, the trivial character, then 
\beq Z(x,\Lambda , \psi) = \frac{\pi}{|\Lambda |}  (\ln x + \eta_\Lambda ) + O \left( x^{ - \frac{1}{2}} \right) \space  \text{ \space \space \space as $x \ra \infty $.}  \eeq

(2) If $\psi \neq \text{id}$  then $ \lim_{x\ra \infty} Z(x,\Lambda , \psi) $ exists and 
$$ Z(x,\Lambda , \psi) = L(\Lambda , \psi ) + O \left( x^{- \frac{1}{2}} \right) \space  \text{ \space \space \space as $x \ra \infty $.}  $$
\ep
Here $ \eta_\Lambda $ can be thought of as an analogue of the Euler constant $\gamma$ for the lattice $\Lambda \subset \RR^2.$ 
\pf
(1) This is proven in  \cite[page 298]{Elstrodt}.  

(2) By Lemma~\ref{L:latticetail} 
$$  \sum_{ \substack{ \mu \in \Lambda \\  x < | \mu |^2 < p  } } \frac{ \psi(\mu)}{ |\mu |^2} = O(x^{- \frac{1}{2}}).   $$ Here the implied constant does not depend on $p.$  Hence, convergence follows from the Cauchy criterion.  Letting $ p \ra \infty $ shows that 
$$ Z(x, \Lambda, \psi) = L(\Lambda, \psi) + O \left(x^{-\frac{1}{2}} \right).
$$
\epf
 
\subsection{Kronecker's Second Limit Formula}
Let $u,v$ be real numbers which are both not integers and  write  $ \tau = x + i y $ $(y > 0).$   For $ \R(s) > 1 $ set 
$$  E_{u,v}(\tau, s) \df  \lim_{x \ra \infty} \psum_{|m \tau + n|^{2} < x}  e^{ 2\pi i (mu+nv)} \frac{y^s}{|m \tau + n|^{2s}}. 
$$
Here the prime in the sum means to leave out zero.
The series converges  uniformly and absolutely on compact subsets of $ \R(s) > 1. $  Thus we can also define $$  E_{u,v}(\tau, s)  = \sum_{ (n,m) \neq (0,0) }  e^{ 2\pi i (mu+nv)} \frac{y^s}{|m \tau + n|^{2s}} $$ since the order of summation is not important when a sum converges absolutely.
In the sums, $n$ and $m$ are understood to vary over  the integers. 

We have  (\cite{Siegel} or  \cite[page 276]{Lang}\footnote{There appears to be a typographical error in the definition of the Siegel function on page 276 of the second edition of \cite{Lang}.  The correct definition appears on page 262.\\})
\begin{lem} \label{P:kron}
The function $ E_{u,v} (\tau, s) $ can be continued to an entire function of $s \in \CC,$ and one has   \beq  E_{u,v} (\tau, 1) = - 2 \pi \log | g_{-v,u} (\tau) |, \eeq 
where $g_{a_1,a_2} $ is the Siegel function, \beq
g_{a_1,a_2}(\tau) = -q_{\tau}^{(1/2)\textbf{B}_{2} (a_1)} e^{ 2 \pi i a_2 (a_1-1)/2 }(1-q_z) \prod_{n=1}^{\infty} (1 - q_{\tau}^n q_z)(1 - q_{\tau}^n /q_z),
\eeq
$ \textbf{B}_{2} (X) = X^2 -X + 1/6, $ $ q_{\tau} = e^{ 2 \pi i \tau},$  $ q_z =  e^{ 2 \pi i z}, $ and $z = a_1\tau+a_2. $
\end{lem}

We now explain the relationship between $L(\Lambda, \psi)$ and  $ E_{u,v} (\tau, 1). $ The character $\psi $ is determined by $u$ and $v.$  That is  $u$ and $v$ can be chosen so that  
$  \psi(1) =  e^{2 \pi i u} ~\text{and}~ \psi(\tau) = e^{2 \pi i v}. $ 
We now can rewrite 

\beq \label{E:l30034}
Z( x ,\Lambda , \psi) =   \psum_{|m \tau + n|^{2} < x} \frac{e^{ 2\pi i (mu+nv)}}{|m \tau + n|^{2}}. 
\eeq
Here the prime means we leave off zero from the sum.

Formally, ignoring convergence and order of summation, 
$$ \lim_{x \ra \infty}  Z( x ,\Lambda , \psi) = \sum_{(m,n) \neq (0,0)} \frac{e^{ 2\pi i (mu+nv)}}{|m \tau + n|^{2}} =  \frac{1}{y}E_{u,v} (\tau, 1)   $$

In the next section we will make the above argument rigorous.
\subsection{Evaluation of $L(\Lambda, \psi)$}
\bp Let $\Lambda =  \ZZ \oplus \ZZ \tau \subset \CC$ be a lattice with $\I(\tau) > 0, ~ \psi $ a character of  $\Lambda,$  and $u,v \in \RR $ are both not integers satisfying  
$  \psi(1) =  e^{2 \pi i u} ~~\text{and}  $
$ \psi(\tau) = e^{2 \pi i v}. $
Then 
\beq L( \Lambda, \psi ) = \frac{-2 \pi}{y} \log \left| g_{-v, u} \left( \tau \right) \right|.      \eeq

\ep

\pf
For $\R(s) > 1 $ recall that 
$$ y^{-s} E_{u,v}(\tau,s)  =  \lim_{x \ra \infty} \psum_{|m \tau + n|^{2} < x} \frac{e^{ 2\pi i (mu+nv)}}{|m \tau + n|^{2s}}. $$  
For $t \in [1,2] $ set  
$$ f(t) \df   \lim_{x \ra \infty} \psum_{|m \tau + n|^{2} < x} \frac{e^{ 2\pi i (mu+nv)}}{|m \tau + n|^{2t}}.   $$
By proposition \eqref{P:latticesum} $ \lim_{x \ra \infty} Z(x, \Lambda, \psi) = L( \Lambda, \psi ) $ converges, and by
equation \eqref{E:l30034} $L( \Lambda, \psi ) = f(1). $   By  Lemma~\ref{P:kron} $y^{-s} E_{u,v}(\tau,s)$ can be continued to an entire function (also denoted by $y^{-s} E_{u,v}(\tau,s)$).  In particular $y^{-s} E_{u,v}(\tau,s)$ is continuous and 
$$
\lim_{t \ra 1^{+}} y^{-t} E_{u,v}(\tau,t) =  y^{-1} E_{u,v}(\tau,1)= \frac{-2 \pi}{y} \log \left| g_{-v, u} \left( \tau \right) \right|.
$$
By definition, $f(t) =y^{-t} E_{u,v}(\tau,t)  $ for $t \in (1,2].$ If we can show that $f$ is continuous on  $[1,2] $  then 
\begin{multline*} L( \Lambda, \psi ) = f(1) = \lim_{t \ra 1^{+}} f(t) = \lim_{t \ra 1^{+}}  y^{-t} E_{u,v}(\tau,t)=  y^{-1} E_{u,v}(\tau,1) \\= \frac{-2 \pi}{y} \log \left| g_{-v, u} \left( \tau \right) \right|. \end{multline*}

The function $f$  is a continuous function on $[1,2] $
by Lemma  \ref{L:latticetail}.  To see this observe that  by the Cauchy criterion the sum $f$ converges uniformly for $t \in [1,2]. $  Since a uniformly convergent sum of continuous functions is continuous, the proposition is proved. 
\epf

\subsection{Proof of lemma \eqref{L:latticetail}   }
\begin{figure}[htb]
\begin{center}
\includegraphics[height=3in,width=4in]{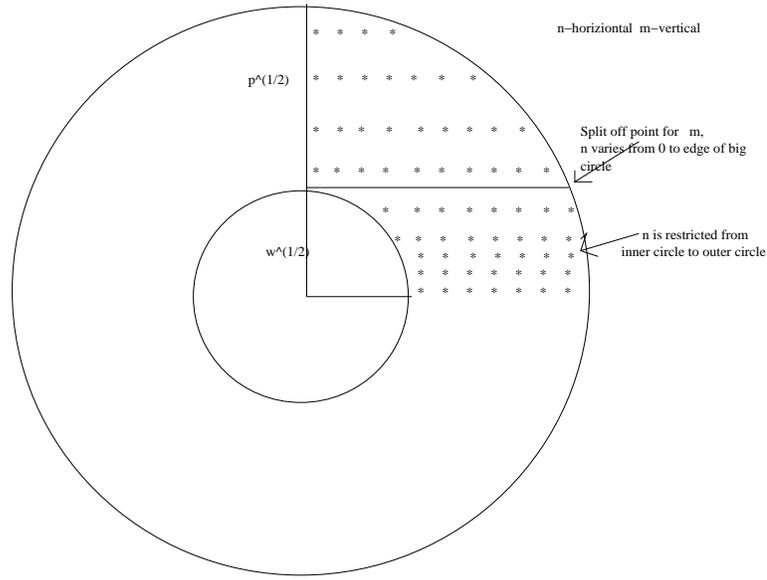}
\caption{Lattice points for the sum.}
\end{center}
\label{fig:latticesum}
\end{figure}

\begin{lem} \label{L:latticetail}
For $ s \in [1,2], p > w > 0, $ and $v,u$ both not integers, $\tau = x+ iy$, $\I(\tau) > 0$  let 
$$ h(w,s) \df  \sum_{ w < |m \tau + n |^2 < p}  \frac{e^{ 2\pi i (mu+nv)}}{|m \tau + n|^{2s}}.  $$  
Then $$ h(w,s) = O \left( w^{\frac{1}{2} - s} \right) $$ uniformly for  $s \in [1,2]. $
\end{lem}
The implied constant does not depend on $w$ or $s.$  Our  proof will work for $s \in ( \frac{1}{2}, \infty). $  However the interval $[1,2]$ is sufficient for our application.
\pf
Without loss of generality we assume $v \notin \ZZ $ and $\R(\tau) > 0.$
We estimate the sum 
\beq \label{S:0101} \sum_{ w < |m \tau + n |^2 < p}  \frac{e^{ 2\pi i (mu+nv)}}{|m \tau + n|^{2s}}. 
\eeq
Sum \eqref{S:0101} is a sum over lattice points of $\Lambda $ that are outside of the circle of radius $ \sqrt{w}$ but inside the circle of radius $\sqrt{p}.$  It suffices to restrict our sum to the first quadrant of the plane since the sum over the other three quadrants can be estimated similarly.  With this assumption  $ n \geq 0 $ and $ m \geq 0 $ and our sum can be written as an explicit double sum with $n$ and $m$ separated as,
\beq \label{S:8838}
\sum_{m = 0}^{  [\sqrt{w}/| \tau | ] }e^{ 2\pi i (mu)}~ \sum_{n=q(w) }^{q(p) }  \frac{e^{ 2\pi i nv}}{|m \tau + n|^{2s}} + \sum_{m =  [\sqrt{w}/| \tau | ] + 1}^{  [\sqrt{p}/| \tau | ] }e^{ 2\pi i (mu)}~ \sum_{n=0}^{ q(p) }  \frac{e^{ 2\pi i nv}}{|m \tau + n|^{2s}},
\eeq where  
$q(\alpha) =  [\sqrt{\alpha - m^2 y^2} -mx ]. $  The brackets, $ [~] $ represent the greatest integer function.  The sum above is split over $m$. For in the first case horizontal translations of lattice points will intersect the inner circle $ |z| = \sqrt{w} $ and hence $n$ must be restricted while the second sum $n$ is free to start at zero.  See figure \eqref{fig:latticesum}.  
We encourage the reader to write down the explicit double sum in the simple case of  $\ZZ \oplus \ZZ i. $ 

We estimate the inner sum using a Stieltjes integral.  For $t \geq 0 $ let $\pi(t) = [t]. $ Then for $0 \leq a < b$
$$\sum_{n=a }^{b}  \frac{e^{ 2\pi i nv}}{|m \tau + n|^{2s}} = \int_{a'}^{b'}\frac{e^{ 2\pi i tv}}{|m \tau + t|^{2s}}~d\pi(t) =  
 \int_{a'}^{b'}\frac{e^{ 2\pi i tv}}{((t + mx)^2 + (my)^2 ) ^{s}}~d\pi(t), $$
where $a' = a - \delta,~b' = b + \delta $ for any  $\delta > 0$ sufficiently small  (the Stieltjes integral must have a slightly larger integration domain to include the endpoints of the sum).
Next using integration by parts with 
$$U(t) = \frac{1}{((t + mx)^2 + (my)^2 ) ^{s}}, $$ 
$$ dV(t) =  e^{ 2\pi i tv}~d\pi(t), $$
$$  V(t) = \int_{0}^{t} e^{ 2\pi i x v}~d\pi(x),   $$
$$ dU(t) = \frac{d}{dt} U(t)~dt, $$
and noting that 
$$ V(t) = \int_{0'}^{t'} e^{ 2\pi i x v}~d\pi(x) = \sum_{k=0}^{[t]} e^{ 2\pi i k v}  $$ 
satisfies $ -C \leq  V(t) \leq  C  $ since $v$ is not an integer, we can estimate the   integral  
 $$  \left| \int_{a'}^{b'}\frac{e^{ 2\pi i tv}}{((t + mx)^2 + (my)^2 ) ^{s}}~d\pi(t)  \right| \leq 
\left|  \left. U(t)V(t) \right]_{a'}^{b'}   \right| + \left|  \int_{a'}^{b'} |V(t)|~|dU(t)| \right| 
$$
$$  
\leq 
C \left|  \left. U(t) \right]_{a'}^{b'}   \right| + C \left|  \int_{a'}^{b'} |dU(t)| \right| \leq 2C \left|  \left. U(t) \right]_{a'}^{b'} \right| \leq 2C |U(a)| + 2C|U(b)|.
$$
Note that $ |dU(t)| $ is easily estimated since $U(t) $ is monotone decreasing and that the primes can be dropped since $\delta $ can be made  arbitrarily small.

What remains is to apply the above estimate to each of the inner sums and crudely put absolute values around all the other terms.  

The first sum can be estimated as follows,
\begin{multline*}
\left| \sum_{m = 0}^{  [\sqrt{w}/| \tau | ] }e^{ 2\pi i (mu)}~ \sum_{n=q(w) }^{q(p) }  \frac{e^{ 2\pi i nv}}{|m \tau + n|^{2s}}  \right|  \\ \leq 
\sum_{m = 0}^{  [\sqrt{w}/| \tau | ] } \left| e^{ 2\pi i mu} \right|~ \left| \sum_{n=q(w) }^{q(p) }  \frac{e^{ 2\pi i nv}}{|m \tau + n|^{2s}} \right| \\  = 
\sum_{m = 0}^{  [\sqrt{w}/| \tau | ] } ~ \left| \sum_{n=q(w) }^{q(p) }  \frac{e^{ 2\pi i nv}}{|m \tau + n|^{2s}} \right| \\= \sum_{m = 0}^{  [\sqrt{w}/| \tau | ] } ~ \left| \int_{q(w)' }^{q(p)' }  \frac{e^{ 2\pi i t v}}{|m \tau + n|^{2s}}~d\pi(t)  \right|  \\ \leq 2C \sum_{m = 0}^{  [\sqrt{w}/| \tau | ] } \left( |U(q(w)| + |U(q(p))| \right).
\end{multline*}
Simplifying $ U(q(w))$ and $U(q(p))$ using 
$$U(t) = \frac{1}{((t + mx)^2 + (my)^2 ) ^{s}} $$ and  $$q(\alpha) =  [\sqrt{\alpha - m^2 y^2} -mx ] $$
yields 
\begin{multline*}
2C \sum_{m = 0}^{  [\sqrt{w}/| \tau | ] } \left( |U(q(w))| + |U(q(p))| \right) \\ \leq  2C \sum_{m = 0}^{  [\sqrt{w}/| \tau | ] }  \left(  \frac{1}{[w]^s}  +  \frac{1}{[p]^s}  \right) \\ \leq 2C[\sqrt{w} / | \tau | +1]( \left(    \frac{1}{[w]^s} +  \frac{1}{[p]^s} \right) \\ \leq
2C[\sqrt{w} / | \tau | +1]( \left(    \frac{1}{[w]^s} +  \frac{1}{[w]^s} \right)  = O \left( w^{\frac{1}{2} - s} \right).
\end{multline*}
The last inequality follows from $w \leq p. $

We apply the same method to the second sum and obtain,
\begin{multline*}
 \left| \sum_{m =  [\sqrt{w}/| \tau | ] + 1}^{  [\sqrt{p}/| \tau | ] }e^{ 2\pi i (mu)}~ \sum_{n=0}^{ q(p) }  \frac{e^{ 2\pi i nv}}{|m \tau + n|^{2s}} \right| \\ \leq 2C \sum_{m =  [\sqrt{w}/| \tau | ] + 1}^{  [\sqrt{p}/| \tau | ] } \left( |U(0)| + |U(q(p))| \right)  \\ \leq
 2C \sum_{m =  [\sqrt{w}/| \tau | ] + 1}^{  [\sqrt{p}/| \tau | ] }  \left( \frac{1}{((mx)^2 + (my)^2 ) ^{s}} + \frac{1}{[p]^s} \right) \\ = 2C  \sum_{m =  [\sqrt{w}/| \tau | ] + 1}^{  [\sqrt{p}/| \tau | ] }  \frac{1}{((mx)^2 + (my)^2 ) ^{s}} + 2C \sum_{m =  [\sqrt{w}/| \tau | ] + 1}^{  [\sqrt{p}/| \tau | ] }   \frac{1}{[p]^s}. 
\end{multline*}

It remains to show that both sums above are bounded by a term of growth $ O \left( w^{\frac{1}{2} - s} \right). $  

For the right hand sum,
$$  2C \sum_{m =  [\sqrt{w}/| \tau | ] + 1}^{  [\sqrt{p}/| \tau | ] }   \frac{1}{[p]^s}  \leq  2C  [\sqrt{p}/| \tau | ] \frac{1}{[p]^s} = O \left( w^{\frac{1}{2} - s} \right) $$ since $p > w. $ 

For the left hand sum
$$ 2C  \sum_{m =  [\sqrt{w}/| \tau | ] + 1}^{  [\sqrt{p}/| \tau | ] }  \frac{1}{((mx)^2 + (my)^2 ) ^{s}} \leq C_1 \sum_{m= [\sqrt{w}]}^\infty \frac{1}{m^{2s}}= O \left( w^{\frac{1}{2} - s} \right),$$ for some $C_1 >0.$
The last equality holds uniformly for $s \in [1,2] $ by the integral test.  
\epf
\section{The Explicit Trace}
In this section our main goal is to give an explicit formula for $$\int_{\F_A} \tr_V (K_{\Gamma}(P,P))  \,dv(P).  $$  
Following Selberg's original method we  decompose 
$$K_\Gamma(P,Q) = \sum_{\gamma \in \Gamma} \chi(\gamma)
K(P,\gamma Q)$$  into various sub-sums.  Depending on their \emph{type.} The types are as follows,   ``id'' is the identity,  ``par''  are the parabolic elements, ``ce'' are the cuspidal elliptic elements\footnote{These elliptical elements share a common fixed point in $\PP$ with some parabolic element in $\Gamma.$\\ }, ``nce'' are the non-cuspidal elliptical elements, ``lox'' are the hyperbolic and loxodromic elements, and ``cusp'' $=$ ``par'' $\union$ ``ce''.
For each $ S \in \{\ID,\PAR,\CE,\NCE,\LOX,\CUSP\}$ set 
$$
K_{\Gamma}^{S}(P,Q) = \sum_{\gamma \in \Gamma^S} \chi(\gamma)
K(P,\gamma Q).
$$
Here $\Gamma^S$ denotes the subset of $\Gamma $ consisting of elements of type $S.$  

Following \cite[Section 5.2, Theorem 6.5.1]{Elstrodt}  we have,
\begin{lem} \label{lemNonCuspidal}
\begin{multline}
\int_{\F_A} \tr_V \left( K_{\Gamma}^{\ID}+ K_{\Gamma}^{\NCE}
+  K_{\Gamma}^{\LOX}  \right)(P,P)  \,dv(P)  \\
= \frac{ \vol (\Gamma \setminus \HH) \dim_\CC (V) } { 4 \pi^2 } 
\int_{\RR} h(1+t^2)t^2 \, dt 
+ \sum_{ \{ R \}_{\NCE} } \frac{ \tr(\chi(R)) g(0) \log N(T_0) }{ 4 |\E(R)| \sin^2 \left(\frac{ \pi k }{ m(R) } \right) } 
\\ + \sum_{ \{ T \}_{\LOX}} \frac{\tr(\chi(T)) g( \log N(T))
\log N(T_0)}{|\E(T)| |a(T) - a(T)^{-1}|^2}  + o(1) ~\text{as $A \ra \infty.$}
\end{multline}
\end{lem}
The notations above were defined in \S\ref{sectionSelberg}.

\section*{The Cuspidal Elliptic Elements}
Our next immediate goal is to evaluate  
$$ \int_{\F_A} \tr_V  K_{\Gamma}^{\CE}(P,P)  \,dv(P).  $$

For notational simplicity, we will adopt Assumption \ref{asOne}  from this point on until the end of this thesis.

Denote by $ \cuspi $  set of elements of $\Gamma$ 
which are $\Gamma$-conjugate to an element of $ \Gamma_\infty
\setminus \Gamma_\infty^\prime = \{ \gamma \in \gi~|~\gamma~\text{is not parabolic nor the identity element}~\}. $  
We fix representatives  of conjugacy classes of  $\cuspi,~
g_{1}, \dots , g_{d}\footnote{There are only finitely many distinct conjugacy classes of elliptic elements in a cofinite Kleinian group.} $ that have the form 
\beq g_{i} =  \left(\begin{array}{cc}
 \epsilon_{i} & \epsilon_{i} \omega_{i} \\
 0 & \left( \epsilon_{i}  \right)^{-1}\end{array}
\right).
\eeq
For $g \in \cuspi,$ let $\mC (g) $ denote the centralizer in $\Gamma $ of $g.$   In addition, let   $\{ p_i,\infty \} $ be the set of fixed points in $\cinf$ of the element $g_i.$  Since $g_i $ is a cuspidal elliptic element it follows that  $p_i $ is a cusp of $\Gamma $ (see \cite{Elstrodt} page 52).  Hence by Assumption~\ref{asOne} there is an element $\gamma_i \in \Gamma $ with $\gamma_i \infty = p_i. $  Suppose that $c_i $ is the lower left hand (matrix) entry of $\gamma_i.$ Then we have (see \cite[Pages 302-304]{Elstrodt}),
\bl \label{lemCuspidalElliptic}
\begin{multline}
\int_{\F_A} \tr K_\Gamma^{\CE}(P,P)~ dv(P) \\ = \sum_{i=1}^d \frac{\tr \chi(g_i)}{|\mC(g_i)|}\left[  \frac{2 g(0) (\log|c_i| + \log A)}{|1-\epsilon_i^2|^2} \right. \\+ \left. \frac{1}{|1-\epsilon_i^2|^2} \int_0^\infty g(x) \frac{\sinh x}{\cosh x - 1 +\frac{|1-\epsilon_i^2|^2}{2} }~dx \right]  \\ + o(1) ~\text{as $A \ra \infty.$}
\end{multline}
\el

\section{The Parabolic Elements} \label{sec.par}
This section contains the new features of the Selberg trace formula that are not present in the \tdv and \thds cases.  We remind the reader that Assumption~\ref{asOne} is in effect.

Our main goal for this section is to evaluate 
\beq \int_{\F_A} \tr_V  K_{\Gamma}^{\PAR}(P,P)  \,dv(P). \label{eqParInt} \eeq

\section*{Evaluation of Integral \ref{eqParInt} }

Let  $\pg$ be a fundamental domain for the action\footnote{See \S\ref{secStaCus} for more details on the action.}  of  $\Gamma_{\infty} $ on $\CC $,  
$$\widetilde{\pg} \df  \{ (z,r) \in \HH \, | \, z \in   \pg \, \},$$ and $$ \pg_{A} \df    \{ (z,r) \in \HH \, | \, z \in   \pg, \, r \leq A \,  \}. $$
It follows that 
$\widetilde{\pg} $ is a fundamental domain for the action of  $ \Gamma_{\infty} $  on $\HH $. 

Recall that $ \gip $ is canonically isomorphic to a lattice $  \Lambda_{\infty}. $  For $\mu \in \Lambda_{\infty}$ let $\widehat{\mu}$ denote the corresponding parabolic element in $\gip. $ We will need the following (see \cite[Pages 300-301]{Elstrodt})
\bl \label{lemParIntToSum}
\beq   \int_{\F_A} \tr K_{\Gamma}^{\PAR}(P,P) \, dv(P) =  \psum_{\mu \in \Lambda_{\infty}} \tr( \chi ( \widehat{\mu} ) ) \int_{ \pg_{A}} K(P, \widehat{\mu} P) \,dv(P) + \lto_{A \ra \infty}. \eeq
\el
  
Since $\gip $ is an abelian group, $\chi $ restricted to $\gip $ can be diagonalized.  In other words, there exist lattice characters $ \{ \psi_l \}_{l = 1\dots n} $ so that 
\beq \label{eqTraceChi}
\tr \chi |_{\gip} = \sum_{l=1}^n \psi_l.  
\eeq  Thus it suffices to consider lattice characters instead of unitary representations.
\bl \label{lemParMain} Let $\psi $ be a lattice character of $\Lambda_\infty. $  Then 

(1) For $ \psi = \text{id}, $   
\begin{multline*}
 \psum_{ \mu \in \Lambda_\infty } \psi (\mu)  \int_{\pg_{A} }K(P,  \widehat{\mu} P) =  \\
\frac{1}{[ \Gamma_{\infty}:\Gamma_{\infty}^{\prime} ]} 
\left(  g(0) \log A + \frac{h(1)}{4} + 
g(0) \left( \frac{ \eta_{\infty}}{2} 
- \gamma \right) - 
\frac{1}{2\pi} 
\int_{\RR} h(1+t^2) \frac{\Gamma'}{\Gamma}(1+it) \,dt   \right) \\ + o(1) ~\text{as $A \ra \infty $}.
\end{multline*}

(2) For  $ \psi \neq \text{id}, $ 
$$
 \psum_{ \mu \in \Lambda_\infty} \psi (\mu)  \int_{\pg_{A} }K(P,  \widehat{\mu} P) = 
\frac{g(0)}{[\Gamma_{\infty}:\Gamma_{\infty}^{\prime}]} L(\Lambda, \psi ) +o(1) ~\text{as $A \ra \infty $}.
$$
\el
\pf
(1) is proved in \cite[pages 300-302]{Elstrodt}.

The proof of (2) is  a modification of (1).   Let $C_\infty = \frac{ | \Lambda_\infty  |}{ [ \Gamma_{\infty}:\Gamma_{\infty}^{\prime} ]}, $ then  by the definition of the action of $\widehat{\mu}$ on $\HH, $
\begin{multline} \label{eq393}
 \psum_{ \mu \in \Lambda_\infty } \psi (\mu)  \int_{\pg_{A} }K(P,  \widehat{\mu} P) = 
C_\infty \psum_{ \mu \in \Lambda_\infty } \psi( \mu ) \int_0^A k \left( \frac{ | \mu |^2 }{2 r^2} + 1 \right) \, \frac{dr}{r^3}.
\\ = C_\infty \psum_{ \mu \in \Lambda_\infty } \frac{ \psi (\mu) }{ | \mu |^2 } \int_{ \frac{| \mu |^2 }{2A^2}}^\infty k(u+1) \,du. 
\end{multline}

Since $$ Z(x,\Lambda , \psi) = \sum_{ \substack{ \mu \in \Lambda_\infty  \\  | \mu |^2 \leq x \\ \mu \neq 0 } } \frac{ \psi(\mu)}{ |\mu |^2}, $$  using summation by parts, we can rewrite
\eqref{eq393} as 
\beq
C_\infty \int_0^\infty k(u+1) Z(2A^2u,\Lambda, \psi )\,du.
\eeq
Next we apply Proposition \ref{P:latticesum} to  obtain 
\begin{multline}
C_\infty \int_0^\infty k(u+1) Z(2A^2u,\Lambda, \psi )\,du \\= C_\infty \int_0^\infty k(u+1) \left( L(\Lambda, \psi) + O \left( (2A^{2}u)^{-1/2} \right) \right) \,du. 
\end{multline}
Now we  show that the resulting error term is $ o(1). $  Since $ k(u+1) = O \left( (1+u)^{-4}    \right)  $  ($k$ is a rapid decay function)
$$ \left|  \int_0^\infty k(u+1) O \left( (2A^{2}u)^{-1/2} \right) \,du \right| \leq 
\frac{D}{A} \int_0^\infty (1+u)^{-4} u^{-1/2} \,du = o(1) \,\, \text{as} \,\, A \ra \infty 
 $$ for some $D > 0.$
To complete the proof note that $ \int_0^\infty k(u+1) \,du = g(0). $
\epf
Finally we can evaluate \eqref{eqParInt}:
\bl \label{lemParabolic}
\begin{multline}
\int_{\F_A} \tr K_{\Gamma}^{\PAR}(P,P) \, dv(P) =  \\
\frac{l_\infty }
{| \Gamma_{\infty}: \Gamma_{\infty}^{\prime} |} 
\left(  g(0) \log A + \frac{h(1)}{4} + 
g(0) \left( \frac{ \eta_{\infty}}{2} 
- \gamma \right) - 
\frac{1}{2\pi} 
\int_{\RR} h(1+t^2) \frac{\Gamma'}{\Gamma}(1+it) \,dt   \right)  \\ + 
 \frac{g(0)}{| \Gamma_{\infty}: \Gamma_{\infty}^{\prime} |} \sum_{ l = l_\infty + 1 }^{n} L(\Lambda_{\infty}, \psi_{ l}  ). 
\end{multline}
Here $n=\dim_\CC V, $ $\psi_{ l}$ are the lattice characters associated to the lattice $\Lambda_\infty, $ $l_\infty = \dim_\CC V_\infty^\prime, $ and  $\eta_\infty $ is the analogue of the Euler constant for the lattice $\Lambda_\infty. $ 
\el
\pf
The proof follows immediately from Lemma~\ref{lemParIntToSum}, Equation~\ref{eqTraceChi}, and Lemma~\ref{lemParMain}.
\epf

We have evaluated the truncated trace of $\mathcal{K}$ explicitly as an integral, and by using spectral theory.  Notice that as $A \ra \infty $ the integral over the parabolic sum (Lemma~\ref{lemParabolic}) has a divergent term.  So does the corresponding cuspidal elliptic integral (Lemma~\ref{lemCuspidalElliptic}). By Lemma~\ref{lemSinCanOut} the divergent terms \emph{must} equal the divergent term of the spectral (truncated) trace (Lemma~\ref{lemSpeTra}).  It follows that
$$2 g(0) \log A \sum_{i=1}^d \frac{\tr \chi(g_i)}{|\mC(g_i)||1-\epsilon_i^2|^2} +  
 g(0) \log A \frac{l_\infty }
{| \Gamma_{\infty}: \Gamma_{\infty}^{\prime} |}  - g(0) k_\infty \log A = 0. $$  
By choosing a suitable $k$ so that $g(0) \neq 0 $ we obtain
\bl \label{lemCuspElip}
$$
2  \sum_{i=1}^d \frac{\tr \chi(g_i)}{|\mC(g_i)||1-\epsilon_i^2|^2} +  
  \frac{l_\infty }
{| \Gamma_{\infty}: \Gamma_{\infty}^{\prime} |}  =  k_\infty. 
$$
\el
The formula\footnote{A similar formula is valid for the general case of $\kappa$-many cusps. \\} above is an  application of spectral theory to the group relations of a cofinite hyperbolic three-orbifold\footnote{Notice that all of the terms above are defined simply in terms of group relations. }.   We will use the above lemma  to give a meromorphic continuation of the Selberg zeta-function.

\section*{Completion of the proof of the Selberg Trace Formula}
The Selberg trace formula now follows: combine    Lemma~\ref{lemParabolic}, Lemma~\ref{lemCuspidalElliptic}, Lemma~\ref{lemCuspElip}, Lemma~\ref{lemNonCuspidal}, Lemma~\ref{lemSpeTra}, and Lemma~\ref{lemSinCanOut}.  Note that the divergent terms all cancel by Lemma~\ref{lemSinCanOut} (or we can use  Lemma~\ref{lemCuspElip}).  Finally take the limit as $A \ra \infty.$  See \cite[Section 6.5]{Elstrodt} for more details on combining the lemmas above.

\chapter{The Selberg Zeta Function} 
In this section we define the Selberg zeta-function $Z(s,\Gamma,\chi)$    for cofinite Kleinian groups with finite-dimensional unitary representations, in the right half-plane $\Re(s) > 1.$  We then  evaluate the logarithmic derivative of  $Z(s,\Gamma,\chi)$ and show that $Z(s,\Gamma,\chi)$ admits a meromorphic continuation, subject to some technical assumptions concerning the stabilizer subgroup $\gi.$

\section{The Definition and Motivation} \label{secSZF}

In the celebrated paper \cite{Selberg1} Selberg first defined what is now called ``The Selberg zeta-function\footnote{More precisely, the Selberg zeta-function of a cocompact Fuchsian group.\\}'' as an infinite product over lengths of \emph{primitive closed geodesics}\footnote{Geodesics that do not trace over themselves multiple times.\\}, bearing a strong resemblance to the Riemann zeta-function. Surprisingly, the Selberg zeta-function satisfies a Riemann hypothosis, and encodes both geometric and spectral data of the quotient orbifold\footnote{A Riemann surface if $\Gamma $ is torsion-free.\\} $\Gamma \setminus \hh^2.$  The spectral and geometric connection is made clear when one understands the Selberg zeta-function as a by-product of the Selberg trace formula applied to the resolvent kernel of $\lp.$ 

In Defintion~\ref{defSZ} we will define the Selberg zeta-function for our case of interest. A natural question  arises: \emph{what does our zeta-function have in common with  the original Selberg zeta-function?}  The answer\footnote{An alternative answer is that in the cocompact case, both zeta functions are factors of the regularized (functional) determenant $\det \left( \lp - (1-s^2)\right).$  See \cite{Sarnak} for more details. \\ }: the logarithmic derivatives of both zeta-functions are directly related to the loxodromic (or hyperbolic) contribution of the Selberg trace formula applied to the resolvent kernel of $\lp.$ The term (from the trace formula) in question for our case has the form 
$$
 \sum_{ \{ T \}\LOX}  \frac{   \tr (\chi(T)) \log N(T_{0})}{m(T)|a(T)-a(T)^{-1}|^{2}}N(T)^{-s}.
$$
We show in Lemma~\ref{lemLogDerZ} that the term above is the logarithmic derivatives of a meromorphic function $Z(s,\Gamma,\chi),$ and that it has a product expansion in the right half-plane $R(s)>1.$ 

In order to define $Z(s,\Gamma,\chi)$ we will need some notions concerning centralizer subgroups of loxodromic elements.  For more details see \cite[Sections  5.2,5.4]{Elstrodt}. 

Let  $\Gamma $ be a cofinite  Kleinian group and let $\chi \in \rep.$  Suppose $T\in\Gamma$ is loxodromic (we consider hyperbolic elements as loxodromic elements). Then $T$ is conjugate in $\pc$ to a unique element of the form 
$$
D(T)=
\left(\begin{array}{cc}
a(T) & 0\\
0 & a(T)^{-1}
\end{array}\right) $$
such that $a(T)\in\CC$ has $|a(T)|>1$.  Let $N(T)$  denote  the \emph{norm} of $T,$   defined by  $$N(T) \df |a(T)|^{2},$$ and  let   by  $\mC(T) $ denote  the centralizer of $T$ in $\Gamma.$  There exists a (primitive)  loxodromic element $T_0,$ and a finite cyclic elliptic subgroup  $\en$ of order $m(T), $ generated by an element $E_{T_0} $   such that 
$$\mC(T) = \langle T_{0} \rangle \times \en. $$
Here $\langle  T_{0} \rangle = \{\, T_{0}^{n} ~ | ~ n \in\ZZ ~ \}.$\footnote{Note that by definition $\mathcal{E}(T)=\mathcal{E}(T_0), $  and that $T_0$ is unique up to multiplication by an element of $ \mathcal{E}(T).$ \\}  Next,
Let $\mathfrak{t}_1,\dots, \mathfrak{t}_n, $ and $\mathtt{t'_1},\dots,  \mathtt{t'_n}$ denote the eigenvalues of $\chi(T_0)$ and $\chi(E_{T_0})$ respectively.   The elliptic element $ E_{T_0}$ is conjugate in $\pc$ to an element of the form 
$$\left(\begin{array}{cc}
\zeta(T_0) & 0 \\
0 & \zeta(T_0)^{-1}
\end{array}\right), $$ 
where here $\zeta(T_0)$ is a primitive $2m(T)$-th root of unity.

\bd \label{defSZ}
For $\R(s)>1 $ the Selberg zeta-function $Z(s,\Gamma,\chi)$ is defined by
$$
Z(s,\Gamma,\chi) \df \prod_{ \{T_0 \} \in \mathcal{R}} ~ \prod_{j=1}^{ \dim V} \prod_{  \substack{ l,k \geq 0 \\  c(T,j,l,k)=1   } } \left( 1-\mathfrak{t}_{j} a(T_0)^{-2k} \overline{ a(T_0) ^{-2l}} N(T_0)^{-s - 1}    \right).  
$$
Here the product with respect to $T_0$ extends over a maximal reduced system $\mathcal{R} $ of $\Gamma$-conjugacy classes of primitive loxodromic elements of $\Gamma.$ The system  $\mathcal{R} $ is called reduced if no two of its elements have representatives with the same centralizer\footnote{See \cite{Elstrodt} section 5.4 for more details \\}.  The function  $c(T,j,l,k)$ is defined by 
$$c(T,j,l,k)= \mathtt{t'_j} \zeta(T_0)^{2l}  \zeta(T_0)^{-2k}.$$
\ed

\bl \label{lemLogDerZ}For $\R(s)>1,$ 
$$\frac{d}{ds} \log Z(s,\Gamma,\chi)  = \sum_{ \{ T \}\LOX}  \frac{\tr (\chi(T)) \log N(T_{0})}{m(T)|a(T)-a(T)^{-1}|^{2}}N(T)^{-s}. $$
\el
\pf
It follows from the proof of \cite[Lemma 5.4.2]{Elstrodt} that 
\begin{multline}
\sum_{ \{ T \}\LOX}  \frac{\tr (\chi(T)) \log N(T_{0})}{m(T)|a(T)-a(T)^{-1}|^{2}}N(T)^{-s} \\
 = \sum_{ \substack{ \{T_0 \} \in \mathcal{R} \\ n \geq 0 \\ 1 \leq v \leq m(T_0)}}  \frac{ \tr \chi(T_0^{n+1}E_0^v) \log{N(T_0)}  }{m(T_0) \left| \zeta(T_0)^v a(T_0)^{n+1} - \zeta(T_0)^{-v} a(T_0)^{-n-1} \right|^2} N(T_0)^{-s(n+1)}. 
\end{multline}
Next since $T_0$ commutes with $E_{T_0}$ we can diagonalize the restriction of $\chi$  to $\mathcal{C}(T) $ and continue the equality to

\begin{multline*}
  = \sum_{ \substack{ \{T_0 \} \in \mathcal{R} \\ n \geq 0 \\ 1 \leq v \leq m(T_0) \\ 1 \leq j \leq \dim V}} \frac{ \mathfrak{t}_j^{n+1} \mathtt{t'}_j^v \log{N(T_0)}  }{m(T_0) \left| \zeta(T_0)^v a(T_0)^{n+1} - \zeta(T_0)^{-v} a(T_0)^{-n-1} \right|^2} N(T_0)^{-s(n+1)} \\ 
  = \sum_{ \substack{ \{T_0 \} \in \mathcal{R} \\ n \geq 0 \\ 1 \leq v \leq m(T_0) \\ 1 \leq j \leq \dim V  }} \frac{ \mathfrak{t}_j^{n+1} \mathtt{t'}_j^v \log{N(T_0)}  }{m(T_0) \left( 1-d(T_0)  \right) \left( 1-\overline{d(T_0)}  \right) } N(T_0)^{-s(n+1)} \\
  = \sum_{ \substack{ \{T_0 \} \in \mathcal{R} \\ n \geq 0 
  \\ 1 \leq v \leq m(T_0) \\ l,k \geq 0 \\ 1 \leq j \leq \dim V    }   }  \frac{ \mathfrak{t}_j^{n+1} \mathtt{t'}_j^v N(T_0)^{-s(n+1)} \log{N(T_0)}  \left( d(T_0)  \right)^k \left( \overline{d(T_0)}  \right)^l}{m(T_0)}.
\end{multline*}

where $$d(T_0) =  \zeta(T_0)^{-2v}a(T_0)^{-2(n+1)}  $$
Next we sum over the $v-$index (note that it is a geometric sum of an $m(T_0)-$th root of unity) observe that the sum is non-zero only when 
$$ \mathtt{t'}_j \zeta(T_0)^{2l}\zeta(T_0)^{-2k}=1 $$ or using our notation $c(T,j,l,k)=1.$   The  equality continues as 
\begin{multline*}
= \sum_{ \substack{ \{T_0 \} \in \mathcal{R} \\ n \geq 0  \\ l,k \geq 0 \\ 1 \leq j \leq \dim V \\ c(T,j,l,k)=1 }}   \mathfrak{t}_j^{n+1}  N(T_0)^{-s(n+1)}\log{N(T_0)}  \left( a(T_0)^{-2(n+1)}  \right)^k \left(  \overline{a(T_0)^{-2(n+1)}}  \right)^l  \\
= \sum_{ \substack{ \{T_0 \} \in \mathcal{R} \\ l,k \geq 0 \\ 1 \leq j \leq \dim V \\ c(T,j,l,k)=1}}  \frac{\mathfrak{t}_j a(T_0)^{-2k} \overline{a(T_0)^{-2l}} N(T_0)^{-(s+1)} \log{N(T_0)}}{1-\mathfrak{t}_j a(T_0)^{-2k} \overline{a(T_0)^{-2l}} N(T_0)^{-(s+1)}}
= \frac{Z'(s,\Gamma,\chi)}{Z(s,\Gamma,\chi).}  
\end{multline*}

\epf

\section{The Logarithmic Derivative of the Selberg Zeta-Function}
The first step in obtaining the meromorphic continuation of the zeta-function is to relate its logarithmic derivative to the trace formula.
From this point on Assumption~\ref{asOne} is in effect. 

We apply the Selberg trace formula to the pair of functions,  $$h(w)=\frac{1}{s^2+w-1} - \frac{1}{B^2+w-1} ~~\text{and}~ $$ 
$$ g(x) = \frac{1}{2s}e^{-s|x|} - \frac{1}{2B}e^{-B|x|}, $$ where $1 < \R(s) < \R(B)$  and obtain

\begin{lem}
\begin{multline} \label{eqLogDer}
\frac{1}{2s} \frac{Z^\prime}{Z}(s) - \frac{1}{2B} \frac{Z^\prime}{Z}(B)  
=\frac{1}{2s} \sum_{ \{ T \}\LOX}  \frac{   \tr (\chi(T)) \log N(T_{0})}{m(T)|a(T)-a(T)^{-1}|^{2}}N(T)^{-s} \\
 -\frac{1}{2B}\sum_{ \{ T \}\LOX} \frac{ \tr (\chi(T))  \log N(T_{0})}{m(T)|a(T)-a(T)^{-1}|^{2}}N(T)^{-B} \\
= \sum_{n \in D} \left(\frac{1}{s^2 - s_n^2} - \frac{1}{B^2 - s_n^2}  \right) 
- \frac{1}{4 \pi} \int_\RR  \left(\frac{1}{s^2 +w^2} - \frac{1}{B^2 + w^2}  \right) \frac{\phi^\prime}{\phi}(i w)~dw  \\
+ \frac{l_\infty}{2 \pi [\gi:\gip]} \int_\RR  \left(\frac{1}{s^2 +w^2} - \frac{1}{B^2 + w^2}  \right) \frac{\Gamma^\prime}{\Gamma}(1 + iw)~dw 
+ \frac{\tr \smat(0)}{4s^2} - \frac{\tr \smat(0)}{4B^2} \\
- \frac{l_\infty}{4[\gi:\gip]s^2} +\frac{l_\infty}{4[\gi:\gip]B^2}  \\
- \sum_{i=1}^l \frac{\tr \chi(g_i) }{|C(g_i)||1-\epsilon_i^2|^2} \int_0^\infty  \left( \frac{e^{-sx}}{2s} - \frac{e^{-Bx}}{2B} \right)   \frac{\sinh x}{\cosh x -1 +\frac{|1-\epsilon_i^2|^2}{2} }~dx 
   \\ - \left( \frac{1}{2s}-\frac{1}{2B}\right) \sum_{ \{R \} \text{\emph{nce}}}\frac{\tr _{V}\chi (R) \log N(T_{0})}{4|\ren |\sin ^{2}(\frac{\pi k}{m(R)})}
+ \frac{\vol \left( \Gamma \setminus \HH \right)\dim V  }{4\pi}(s-B)
 \\ - \left( \frac{1}{2s}-\frac{1}{2B}\right) \sum_{i=1}^l \frac{2 \tr \chi(g_i) \log|c_i| }{|C(g_i)||1-\epsilon_i^2|^2} 
 \\ - \left( \frac{1}{2s}-\frac{1}{2B}\right) \frac{1}{[\gi:\gip]} \left(  l_\infty \left( \frac{ \eta_{\infty}}{2} 
- \gamma \right)  +  \sum_{ l = l_\infty + 1 }^{n} L(\Lambda_{\infty}, \psi_{ l}  )   \right).
\end{multline}
\end{lem}
\begin{proof}
The first equality follows from Lemma~\ref{lemLogDerZ}. The second equality follows directly from the Selberg trace formula.
\end{proof}

Equation \eqref{eqLogDer} is  used to exhibit the meromorphic continuation of $Z(s,\Gamma,\chi). $  If we fix $B$ and multiply through   by $2s,$  it is not hard to see that each term on the right of  \eqref{eqLogDer} is meromorphic.  In order to  see that $\sz$ is meromorphic, we must compute the residues of each term on the right of \eqref{eqLogDer}.  We will show that the residues are fractional and that for some $N \in \NN,~\sz^N$ is a meromorphic function.

\begin{thm}
Let  $\Gamma$ be cofinite with one class of cusps at $\zeta = \infty,$ and let $\chi \in \rep.$  
\begin{enumerate}
\item \label{itMer01} If  $[\gi:\gip] = 1$ or $[\gi:\gip] = 2,$  then  $Z(s,\Gamma,\chi)$ is a meromorphic function. 
\item \label{itMer02}If $[\gi:\gip] = 3,$  then there exists a natural number $N,$ $1 \leq N \leq 6,$ so that $\left(Z(s,\Gamma,\chi)\right)^N$ is a meromorphic function.
\end{enumerate}
\end{thm}
\pf
The proof follows from a careful study of \eqref{eqLogDer}.   We must show that after multiplying by $2s,$ each term on the right (of the second equal sign) has at most simple poles with integral or rational residues\footnote{For case (1) the residues must be integer while for case (2) it suffices to show that the residues are rational with bounded denominator. \\}.   This is demonstrated in  Lemma~\ref{lemSpecDiv}, Lemma~\ref{lemTopCaseOne}, Lemma~\ref{lemTopCaseTwo}, and Lemma~\ref{lemTopCaseThree}.  
\epf
We remark that the divisor of the Selberg zeta-function is readily read off from Lemma~\ref{lemSpecDiv}, Lemma~\ref{lemTopCaseOne}, Lemma~\ref{lemTopCaseTwo}, and Lemma~\ref{lemTopCaseThree}. 

Our zeta function satisfies a functional equation.  A standard argument (\cite[Theorem 5.1.5, page 85]{Venkov}) using \eqref{eqLogDer},  Lemma~\ref{lemSpecDiv}, Lemma~\ref{lemTopCaseOne}, Lemma~\ref{lemTopCaseTwo}, and Lemma~\ref{lemTopCaseThree} yields:
\begin{thm} \label{thmFuncEq}Suppose that $[\gi:\gip]=1$ or $[\gi:\gip]=2.$   Then $\sz$ satisfies:
$$
Z(-s,\Gamma,\chi) = Z(s,\Gamma,\chi)\phi(s)\Psi(s,\Gamma,\chi). 
$$
For $[\gi:\gip]=1,$ 
\begin{equation*}
\Psi(s) \df   \left(\frac{\Gamma(1-s)}{\Gamma(1+s)}   \right)^{k_\infty}   \exp \left(-\frac{\vol \left( \Gamma \setminus \HH \right)\dim V  }{3\pi}s^3+Es + C \right) \end{equation*}
and 
\begin{equation*}
 E \df \sum_{ \{R \} \text{\emph{nce}}}\frac{\tr _{V}\chi (R) \log N(T_{0})}{4|\ren |\sin ^{2}(\frac{\pi k}{m(R)})}
 +  \left(  k_\infty \left( \frac{ \eta_{\infty}}{2} 
- \gamma \right) +  \sum_{ l = k_\infty + 1 }^{n} L(\Lambda_{\infty}, \psi_{ l}  )   \right).
\end{equation*} 
For $[\gi:\gip]=2,$  \begin{multline*}
\Psi(s) \df \\  \left(\frac{\Gamma(1-s)}{\Gamma(1+s)}   \right)^{l_\infty} \left(\prod_{k=1}^\infty \exp(-k(-1)^{k} \frac{(k-1)^2-s^2}{(k+1)^2-s^2}  \right)^{k_\infty/2-l_\infty/2}  \\ \cdot \exp \left(-\frac{\vol \left( \Gamma \setminus \HH \right)\dim V  }{3\pi}s^3+Es + C \right) \end{multline*}
and 
\begin{multline*}
 E \df \sum_{ \{R \} \text{\emph{nce}}}\frac{\tr _{V}\chi (R) \log N(T_{0})}{4|\ren |\sin ^{2}(\frac{\pi k}{m(R)})}
+  \sum_{i=1}^l \frac{2 \tr \chi(g_i) \log|c_i| }{|C(g_i)||1-\epsilon_i^2|^2} 
 \\ + \frac{1}{[\gi:\gip]} \left(  l_\infty \left( \frac{ \eta_{\infty}}{2} 
- \gamma \right) +  \sum_{ l = l_\infty + 1 }^{n} L(\Lambda_{\infty}, \psi_{ l}  )   \right).
\end{multline*}

 The constant\footnote{The value of $C$ can be read off by letting $s \ra 0$ in the functional equations.  Its value  depends on whether  $\phi(0)$ is $1$ or $-1$ and the multiplicity of $\sz$  at $s=0.$ \\} $C$ satisfies the equation: $\exp(C)=\pm 1.$  
\end{thm}

\bl \label{lemSpecDiv}
The expression 
$$\sum_{n \in D} 2s \left(\frac{1}{s^2 - s_n^2} - \frac{1}{B^2 - s_n^2}  \right) 
- \frac{1}{4 \pi} \int_\RR 2s \left(\frac{1}{s^2 +w^2} - \frac{1}{B^2 + w^2}  \right) \frac{\phi^\prime}{\phi}(i w)~dw   $$ has only simple poles and integral residues:

(a) at the points $\pm s_j$ on the line $\R(s)=0$ and on the interval $[-1,1].$  Each point $s_j$ is related to an eigenvalue $ \lambda_j $ of the discrete spectrum of $\lp $ by $1-s_j^2 = \lambda_j.$  The residue of each $s_j$  is equal to the multiplicity of the corresponding eigenvalue. If $\lambda = 1,$ is an eigenvalue of $\lp, $ then the residue of the point $s_j = 0,$ is twice the multiplicity of $\lambda;$ 

(b) at the points $\rho_j $ that are poles of  $\smat(s),$ which lie in the half-plane $\R(s) < 0. $ The residue of each $\rho_j $ is non-negative\footnote{The point $\rho_j $ is a zero of $Z(s,\Gamma,\chi)$ and a pole of $\smat(s).$ We understand the multiplicity of a pole  as non-negative number (not as a negative number). \\} and equal to its multiplicity as a pole of $\smat(s).$  
\el
\pf
The computation involves elementary complex analysis. See \cite[Section 5.1]{Venkov}.
\epf

The residues above come  from terms that are related to the spectral and scattering theory of $\lp.$  The remaining residues are computed using group theoretic data involving $\Gamma$ and $\chi.$  The poles and zeros of $\Z(s,\Gamma,\chi)$ that correspond to these residues are commonly called \emph{topological} or \emph{trivial}\footnote{We refer to them as topological.\\}.

\section{The Topological Zeros and Poles} \label{secTopZeros}
The computation of the topological residues   is considerably more complicated than the corresponding spectral computation.  Poles can only  arise from the following terms (excluding the spectral terms previously dealt with) of \eqref{eqLogDer} (note that we multiplied all terms through by $2s$):
 \begin{multline} \label{eqLogDerTop}
 \frac{l_\infty}{2 \pi [\gi:\gip]} \int_\RR  \frac{2s}{s^2 +w^2}  \frac{\Gamma^\prime}{\Gamma}(1 + iw)~dw  - \frac{l_\infty}{2[\gi:\gip]s}  +  \frac{\tr \smat(0)}{2s}
\\ - \sum_{i=1}^l \frac{\tr \chi(g_i) }{|C(g_i)||1-\epsilon_i^2|^2} \int_0^\infty   e^{-sx}  \frac{\sinh x}{\cosh x -1 +\frac{|1-\epsilon_i^2|^2}{2} }~dx. 
\end{multline}

The first two terms come  from the parabolic elements of $\Gamma$, the third from the spectral trace, and the last from the cuspidal elliptic elements of $\Gamma.$  It is remarkable that the last three terms need to be taken together in order to compute the residue at $s=0,$ while the first and last are needed to compute the residues on the negative real axis.
 
It is well known that 
\begin{multline} \label{eqGammaInt}
 \frac{l_\infty}{2 \pi [\gi:\gip]} \int_\RR  \frac{2s}{s^2 +w^2}  \frac{\Gamma^\prime}{\Gamma}(1 + iw)~dw \\ = 
\frac{l_\infty}{[\Gamma_\infty:\Gamma_\infty^\prime]} \left( \frac{\Gamma^\prime}{\Gamma}(1-s) + \sum_{k=1}^\infty \left( \frac{1}{s+k} + \frac{1}{s-k} \right)  \right).
\end{multline}
In order to obtain a  similarly explicit formula for 
\beq \label{eqCE}
 \sum_{i=1}^l \frac{\tr \chi(g_i) }{|C(g_i)||1-\epsilon_i^2|^2} \int_0^\infty   e^{-sx}  \frac{\sinh x}{\cosh x -1 +\frac{|1-\epsilon_i^2|^2}{2} }~dx
 \eeq
 we must make some technical assumptions.
 \subsection{Case One: $[\gi:\gip] = 1$} In this case, \eqref{eqCE} is not applicable and  $l_\infty = k_\infty = k(\Gamma,\chi)$ (the last equality follows from our assumption that $\infty$ is the only cusp).  Since $\smat(0)$ is a unitary self-adjoint matrix of dimension $k \times k$, its trace consists of a sum of $k$ terms of the form $ \pm 1.$  It follows that  $\frac{1}{2}(\tr \smat(0) - k)$ is an integer. After applying equations \eqref{eqLogDerTop} and \eqref{eqGammaInt} we have: 
 \bl \label{lemTopCaseOne}
 Suppose $[\gi:\gip] = 1.$ Then the poles of \eqref{eqLogDerTop} are simple and are located at the points $s=-1,-2,\dots,$ with residue $k_\infty$ and at the point $s=0,$ with residue $\frac{1}{2}(\tr \smat(0) - k_\infty).$
 \el
 \subsection{Case Two: $[\gi:\gip] = 2$ }
 In this case for all $i,~ \epsilon_i = \epsilon \df  \sqrt{-1}  $ and \eqref{eqCE} becomes 
 \beq \label{eqIntTwo}
 \sum_{i=1}^l \frac{\tr \chi(g_i) }{|C(g_i)||1-\epsilon|^2} \int_0^\infty   e^{-sx} \frac{\sinh x}{\cosh x +1  }~dx.
 \eeq
 An application of Lemma~\ref{lemCuspElip} gives us the coefficient of the integral above 
 $$
 \sum_{i=1}^l \frac{\tr \chi(g_i) }{|C(g_i)||1-\epsilon|^2}  =  \frac{1}{2} \left(
  k_\infty  - \frac{l_\infty }{2}  \right). 
 $$ 
 Next, to evaluate the integral in \eqref{eqIntTwo} we appeal to the following formula:
 \beq
 \int_0^\infty e^{-\mu x} (\cosh x - \cos t)^{-1}~dx = \frac{2}{\sin t} \sum_{k=1}^\infty \frac{\sin kt}{\mu + k},
 \eeq
valid for $\R(\mu) > -1$ and $t \neq 2n\pi, $ (see \cite{Grad} formula  3.543.2).  Averaging for the two values  $\mu -1, \mu+1,$ we obtain 
\beq \label{eqIntegral}
\int_0^\infty e^{-sx}\frac{\sinh x}{\cosh x - \cos t} ~dx  = \frac{1}{\sin t} \sum_{k=1}^\infty \sin{kt} \left( \frac{1}{s-1+k} - \frac{1}{s+1+k} \right).
\eeq
Finally, we can evaluate the integral in \eqref{eqIntTwo} by taking the limit as $t \ra \pi $ 
\beq \label{eqIntTwoEx}
  \int_0^\infty   e^{-sx} \frac{\sinh x}{\cosh x +1  }~dx = 
 \sum_{k=1}^\infty k(-1)^{k+1} \left( \frac{1}{s-1+k} - \frac{1}{s+1+k} \right). 
\eeq
Combining the residues above with those from \eqref{eqGammaInt} give us: 
\bl \label{lemTopCaseTwo}
Suppose $[\gi:\gip]=2.$   Then the poles of \eqref{eqLogDerTop} are simple and are located at the points $s =n,$ with $n = -1,-2,\dots,$ and residue\footnote{Note that by defintion $l_\infty \geq k_\infty.$ \\}:  
$$ m_n = \begin{cases}
k_\infty, 	&  \text{if $n$ is odd}, \\
l_\infty - k_\infty,	& \text{else}
\end{cases}
$$ and at the point $s=0$ with residue $\frac{1}{2}(\tr \smat(0) - k_\infty).$
\el
Since the Picard group $\Gamma = \PSL(2,\ZZ[\sqrt{-1}])$ satisfies $[\gi:\gip]=2$ we have: 
\begin{cor}
Let $\Gamma = \PSL(2,\ZZ[\sqrt{-1}]),$ and let $\chi \in \rep.$ Then $Z(s,\Gamma,\chi)$ is a meromorphic function.
\end{cor}

\subsection{Case Three: $[\gi:\gip] = 3$ } 
In this case, the cuspidal elliptic elements $g_i$ that are in \eqref{eqCE} must all be of order three.  Hence   $\epsilon_i \in \{ \sqrt[3]{1}, \sqrt[3]{-1}~\},$ but when  such $\epsilon_i$ is plugged into $|1-\epsilon_i^2|^2$ one obtains $3.$  Thus we can rewrite \eqref{eqCE} as 
\beq \label{eqCaseThreeA}
\sum_{i=1}^l \frac{\tr \chi(g_i) }{|C(g_i)||1-\epsilon_i^2|^2} \int_0^\infty   e^{-sx}  \frac{\sinh x}{\cosh x + \frac{1}{2} }~dx.
\eeq
Next, applying \eqref{eqIntegral} with $t = \frac{2}{3} \pi, $ and Lemma~\ref{lemCuspElip}, we can rewrite \eqref{eqCaseThreeA} as
\begin{multline}
\frac{1}{2} \left(k_\infty - \frac{l_\infty}{3} \right) \left( \left( \frac{1}{s-1+1} - \frac{1}{s+1+1} \right) - \left( \frac{1}{s-1+2} - \frac{1}{s+1+2} \right) \right.  \\  \left.
+ \left( \frac{1}{s-1+4} - \frac{1}{s+1+4} \right) - \left( \frac{1}{s-1+5} - \frac{1}{s+1+5} \right) + \dots  \right). 
\end{multline}
Finally  combining the residues above with the  residues from \eqref{eqGammaInt} we obtain: 
\bl \label{lemTopCaseThree}
 Suppose $[\gi:\gip]=3.$   Then the poles of \eqref{eqLogDerTop} are simple and are located at the points $s=n,$ with $n = -1,-2,\dots, $ and residue:
$$ m_n = \begin{cases}
\frac{2}{3}l_\infty - k_\infty, 	&  \text{if $n$ is a multiple of 3}, \\
 \frac{1}{6}l_\infty +\frac{1}{2}k_\infty, 	& \text{else} 
\end{cases}
$$  and at the point $s=0$ with residue $\frac{1}{2}(\tr \smat(0) - k_\infty).$
\el
As an application of the lemma above, we have:
\begin{cor}
Let $\Gamma = \PSL(2,\ZZ[-\frac{1}{2}+\frac{\sqrt{-3}}{2}]),$ and $\chi \equiv 1$ (the trivial representation). Then $Z(s,\Gamma,\chi)$ is not a meromorphic function\footnote{This is the first example that the author is aware of where the Selberg zeta-function is not meromorphic.\\} (it is the 6-th root of a meromorphic function).
\end{cor}
\pf
Since $\QQ(\sqrt{-3})$ has class number one, $\PSL(2,\ZZ[-\frac{1}{2}+\frac{\sqrt{-3}}{2}])$ has one class of cusps (\cite[Chapter 7]{Elstrodt}).  In addition it has $\infty$ for a cusp.  Elementary calculations show that $[\gi:\gip]=3.$   By definition, $\chi \equiv 1$ implies that $k_\infty = l_\infty =1.$  The result follows from Lemma~\ref{lemTopCaseThree}.
\epf
\subsection{The Remaining Cases}
The cases that remain are: $[\gi:\gip] = 4$ and $[\gi:\gip] = 6.$  Using the same ideas as in the other three cases one can show the following:
\begin{thm}
Suppose that $[\gi:\gip] = 4.$  Then for some integer $N,$ $\sz^N $ is a meromorphic function.
\end{thm}
The author does not know a good bound for the integer $N$ ($N$ depends on $\Gamma$ and $\chi).$ 
On the other hand we conjecture the following:
\begin{conj*}
Suppose that $[\gi:\gip] = 6.$ Then for some integer $N,$ $\sz^N $ is a meromorphic function.
\end{conj*}
\section{The Entire Function Associated to the Selberg Zeta-Function}
As Fischer \cite[Chapter 3]{Fischer} observed, it is useful to group the elliptic, parabolic, and identity terms, together with the loxodromic terms to define an entire function associated to the Selberg zeta-function called the \emph{Selberg xi-function} $\Xi(s,\Gamma,\chi).$  
\bd
For $\R(s) > 1 $ set\footnote{Theorem \ref{thmXiFun} shows that we are justified in defining $\Xi(s,\Gamma,\chi)$  by its logarithmic derivative. \\} 
\begin{multline*}
\frac{\Xi^\prime}{\Xi}(s,\Gamma,\chi)  =  \sum_{ \{ T \}\LOX}  \frac{   \tr (\chi(T)) \log N(T_{0})}{m(T)|a(T)-a(T)^{-1}|^{2}}N(T)^{-s} 
 \\ + \frac{l_\infty}{2 \pi [\gi:\gip]} \int_\RR  \left(\frac{2s}{s^2 +w^2}  \right) \frac{\Gamma^\prime}{\Gamma}(1 + iw)~dw \\
+ \frac{1}{2s} \left( \tr \smat(0)
 - \frac{l_\infty}{[\gi:\gip]} \right) 
\\ - \sum_{i=1}^l \frac{\tr \chi(g_i) }{|C(g_i)||1-\epsilon_i^2|^2} \int_0^\infty  \left( \frac{e^{-sx}}{2s}  \right)   \frac{\sinh x}{\cosh x -1 +\frac{|1-\epsilon_i^2|^2}{2} }~dx - E,
\end{multline*}
where $E$ is defined in Theorem \ref{thmFuncEq}.
\ed
\begin{thm} \label{thmXiFun}
The Selberg xi-function can be continued to an entire function $\Xi(s,\Gamma,\chi) $ with 
\begin{multline*}
\frac{1}{2s}\frac{\Xi^\prime}{\Xi}(s,\Gamma,\chi) -  \frac{1}{2B}\frac{\Xi^\prime}{\Xi}(B,\Gamma,\chi)\\ = \sum_{n \in D}  \left(\frac{1}{s^2 - s_n^2} - \frac{1}{B^2 - s_n^2}  \right) 
- \frac{1}{4 \pi} \int_\RR  \left(\frac{1}{s^2 +w^2} - \frac{1}{B^2 + w^2}  \right) \frac{\phi^\prime}{\phi}(i w)~dw.  
\end{multline*}
Here $B > 1$, $B > \R(s).$
\end{thm}
\pf
The second assertion follows from \eqref{eqLogDer}.  For the first asserstion: by Lemma~\ref{lemSpecDiv}, after multiplying through by $2s,$ the residues of the expression on the right side of the equal sign are all integers and positive.  Hence $\Xi(s,\Gamma,\chi) $ is entire.
\epf
An explicit product formula (for $\R(s) > 1$) can be obtained for $\Xi(s,\Gamma,\chi) $ in the case of \mbox{$[\gi:\gip]=1,2,3$} by  integrating, and then exponentiating, the explicit formula  for the contribution of the parabolic and cuspidal elliptic elements to the Selberg zeta-function.

\begin{rem} The Selberg zeta-function for symmetric spaces of rank-one was studied by Gangolli and Warner (\cite{Gangolli}, \cite{GangWarn}).  More specifically, they studied torsion-free cocompact quotients with finite-dimensional unitary representations, and non-cocompact torsion-free quotients with  trivial representations (the scalar case).  
\end{rem}

\bibliography{mthesis}

\providecommand{\bysame}{\leavevmode\hbox to3em{\hrulefill}\thinspace}
\providecommand{\MR}{\relax\ifhmode\unskip\space\fi MR }
\providecommand{\MRhref}[2]{%
  \href{http://www.ams.org/mathscinet-getitem?mr=#1}{#2}
}
\providecommand{\href}[2]{#2}
\begin{thebibliography}{EGM98}

\bibitem[CdV81]{Verdi}
Yves Colin~de Verdi{\`e}re, \emph{Une nouvelle d\'emonstration du prolongement
  m\'eromorphe des s\'eries d'{E}isenstein}, C. R. Acad. Sci. Paris S\'er. I
  Math. \textbf{293} (1981), no.~7, 361--363.

\bibitem[EGM98]{Elstrodt}
J.~Elstrodt, F.~Grunewald, and J.~Mennicke, \emph{Groups acting on hyperbolic
  space}, Springer Monographs in Mathematics, Springer-Verlag, Berlin, 1998,
  Harmonic analysis and number theory.

\bibitem[Fad69]{Faddeev}
L.~D. Faddeev, \emph{The eigenfunction expansion of {L}aplace's operator on the
  fundamental domain of a discrete group on the {L}oba\v cevski\u\i\ plane},
  Transactions of the {M}oscow {M}athematical {S}ociety for the year 1967
  ({V}olume 17) (1969), 357--386.

\bibitem[Fis87]{Fischer}
J{\"u}rgen Fischer, \emph{An approach to the {S}elberg trace formula via the
  {S}elberg zeta-function}, Lecture Notes in Mathematics, vol. 1253,
  Springer-Verlag, Berlin, 1987.

\bibitem[Gan77]{Gangolli}
Ramesh Gangolli, \emph{Zeta functions of {S}elberg's type for compact space
  forms of symmetric spaces of rank one}, Illinois J. Math. \textbf{21} (1977),
  no.~1, 1--41. \MR{MR0485702 (58 \#5524)}

\bibitem[GR65]{Grad}
I.~S. Gradshteyn and I.~M. Ryzhik, \emph{Table of integrals, series, and
  products}, Fourth edition prepared by Ju. V. Geronimus and M. Ju. Ce\u\i
  tlin. Translated from the Russian by Scripta Technica, Inc. Translation
  edited by Alan Jeffrey, Academic Press, New York, 1965.

\bibitem[GW80]{GangWarn}
Ramesh Gangolli and Garth Warner, \emph{Zeta functions of {S}elberg's type for
  some noncompact quotients of symmetric spaces of rank one}, Nagoya Math. J.
  \textbf{78} (1980), 1--44.

\bibitem[Hej83]{Hejhal}
Dennis~A. Hejhal, \emph{The {S}elberg trace formula for psl(2,r). {V}ol. 2},
  Lecture Notes in Mathematics, vol. 1001, Springer-Verlag, Berlin, 1983.

\bibitem[Iwa02]{Iwaniec}
Henryk Iwaniec, \emph{Spectral methods of automorphic forms}, second ed.,
  Graduate Studies in Mathematics, vol.~53, American Mathematical Society,
  Providence, RI, 2002.

\bibitem[Lan76]{Langlands}
Robert~P. Langlands, \emph{On the functional equations satisfied by
  {E}isenstein series}, Springer-Verlag, Berlin, 1976, Lecture Notes in
  Mathematics, Vol. 544.

\bibitem[Lan87]{Lang}
Serge Lang, \emph{Elliptic functions}, second ed., Graduate Texts in
  Mathematics, vol. 112, Springer-Verlag, New York, 1987, With an appendix by
  J. Tate.

\bibitem[Roe66]{Roelcke}
Walter Roelcke, \emph{Das {E}igenwertproblem der automorphen {F}ormen in der
  hyperbolischen {E}bene. {I}, {II}}, Math. Ann. 167 (1966), 292--337; ibid.
  \textbf{168} (1966), 261--324.

\bibitem[Sar87]{Sarnak}
Peter Sarnak, \emph{Determinants of {L}aplacians}, Comm. Math. Phys.
  \textbf{110} (1987), no.~1, 113--120.

\bibitem[Sel56]{Selberg1}
A.~Selberg, \emph{Harmonic analysis and discontinuous groups in weakly
  symmetric {R}iemannian spaces with applications to {D}irichlet series}, J.
  Indian Math. Soc. (N.S.) \textbf{20} (1956), 47--87.

\bibitem[Sel89]{Selberg2}
Atle Selberg, \emph{Collected papers. {V}ol. {I}}, Springer-Verlag, Berlin,
  1989, With a foreword by K. Chandrasekharan.

\bibitem[Sel91]{Selberg3}
\bysame, \emph{Collected papers. {V}ol. {II}}, Springer-Verlag, Berlin, 1991,
  With a foreword by K. Chandrasekharan.

\bibitem[Sie80]{Siegel}
Carl~Ludwig Siegel, \emph{Advanced analytic number theory}, second ed., Tata
  Institute of Fundamental Research Studies in Mathematics, vol.~9, Tata
  Institute of Fundamental Research, Bombay, 1980.

\bibitem[Ven82]{Venkov}
A.~B. Venkov, \emph{Spectral theory of automorphic functions}, Proc. Steklov
  Inst. Math. (1982), no.~4(153), ix+163 pp. (1983), A translation of Trudy
  Mat. Inst. Steklov. \textbf{153} (1981).

\end{thebibliography}
\bibliographystyle{amsalpha}

\end{document}